\documentclass[12pt]{article}
\usepackage{amsmath,amsfonts,amssymb}
\usepackage[all]{xy}
\usepackage{graphicx}
\usepackage{epsf}
\usepackage{cancel}
\usepackage{multirow}
\usepackage{rotating}
\usepackage{arydshln}
\usepackage{mathrsfs}
\usepackage[colorlinks=true,urlcolor=black]{hyperref}
\hypersetup{			
backref = true,			
pagebackref = true,		
hyperindex = true, 		
colorlinks = true, 		
breaklinks = true, 		
urlcolor = blue, 		
linkcolor = blue, 		
bookmarks = true,		
bookmarksopen = true,  
citecolor=red,
}

\advance \textheight by \topskip
 

\newtheorem{theorem}{Theorem}[section]

\newtheorem{corollary}[theorem]{Corollary}

\newtheorem{definition}[theorem]{Definition}
\newtheorem{example}[theorem]{Example}

\newtheorem{lemma}[theorem]{Lemma}

\newtheorem{proposition}[theorem]{Proposition}
\newtheorem{remark}[theorem]{Remark}

\newcommand{\qeed}{\hfill\textrm{QED}\break\null}
\newenvironment{demo}{\noindent\textit{Proof.}~}{\qeed}

\def\NN{\mathbb N}
\def\ZZ{\mathbb Z}

\def\FF{\mathbb F}

\def\KK{\mathbb K}
\def\M{{\cal A}}
\def\G{GL(V)}

\newcommand{\beq}{\begin{equation}}
\newcommand{\eeq}{\end{equation}}
\newcommand{\beqa}{\begin{eqnarray}}
\newcommand{\eeqa}{\end{eqnarray}}
\newcommand{\A}{{\mathrm {Aut}}}

\newcommand{\Q}{{Q_{\text{mod}}}}
\newcommand{\noi}{\noindent}
\newcommand{\nn}{\nonumber}

\def\>{\rangle}
\def\<{\langle}

\providecommand{\keywords}[1]{{{Keywords:}} #1}
\begin{document}

\title{Polynomial invariants and moduli of generic two-dimensional commutative algebras}
{\bf }

\author{
{\sf M. Rausch de Traubenberg}\thanks{e-mail:
michel.rausch@iphc.cnr.fr}$\,\,$${}^{a}$ and
{\sf M. J. Slupinski}\thanks{e-mail:
marcus.slupinski@math.unistra.fr}$\,\,$${}^{b}$
\\
{\small ${}^{a}${\it  Université de Strasbourg, CNRS, IPHC UMR 7178, F-67000 Strasbourg, France 
}}\\
{\small ${}^{b}${\it
Université de Strasbourg, CNRS, IRMA UMR 7501,  F-67000 Strasbourg, France 
}}\\
}

\maketitle
\date
\vskip-1.5cm

\vspace{2truecm}

\begin{abstract} 
Let   $V$ be a two-dimensional vector space over a field $\FF$ of characteristic not $2$  or $3$.
We  show there is  a canonical surjection  $\nu$  from the  set of suitably generic
commutative algebra structures on $V$ modulo the action of $GL(V)$ onto the plane $\FF^2$.
In these coordinates, which are quotients of invariant quartic polynomials,  properties  such as associativity and the existence of  zero divisors
are described by simple algebraic conditions.
The map $\nu$ is a bijection  over  the complement of  a degenerate elliptic curve $\Gamma$ and over $\Gamma$
we give an explicit parametrisation of the fibre in terms of Galois extensions of $\FF$.
Algebras  in
 $\nu^{-1}(\Gamma)$ are  exactly those which admit non-trivial automorphisms.
We show how  $\nu$ can be lifted to  a  map  from   the  $SL(V)-$moduli space  to an algebraic hypersurface $\Gamma'$
in a four-dimensional vector space  whose equation is essentially  the classical Eisenstein equation for the covariants of
a binary cubic. This map is the restriction of  a surjective map from the set of stable  commutative algebras on $V$
modulo the action of $SL(V)$ onto $\Gamma'$.

\end{abstract}

\vspace{2cm}

\keywords{Commutative algebras, Moduli space, Binary cubics,  Galois extensions.}

\newpage
\tableofcontents

\newpage
\section{Introduction} 
If $\FF$ is a field it is well known that  the equivalence classes of separable
quadratic extensions of $\FF$ (including the split extension) are classified
by the Galois cohomology group $H^1(\FF,\mathbb Z_2) \cong \mathbb F^*/\mathbb F^{*2}$ 
\cite{serre}. 
It is easy to check that if $\mathbb K/\mathbb F$ is a quadratic extension, the  associated ``twisted'' cubic $Q(x) = x^2 \wedge x$
with values in $\Lambda^2(\mathbb K)$ has three distinct roots in a splitting field and 
generically,
the ``twisted'' cubic of an arbitrary two-dimensional $\mathbb F-$algebra
will have this property. Such algebras   we will call {\it generic} and
in this paper we first  investigate the 
moduli spaces with respect to various group actions of the
set of generic  two-dimensional {\it commutative} $\mathbb F-$algebras.
A more general class of algebras is the class of stable algebras which, by definition are algebras for which
some  invariant polynomial does not vanish. In the latter part of the paper (Section \ref{sec:compact} and the Appendix) we study
moduli spaces of stable two-dimensional commutative $\FF-$algebras and  their relation with moduli spaces of generic
two-dimensional commutative $\FF-$algebras.
Of course, different moduli spaces with respect to  different group actions  
 reveal different features.
In \cite{bhar} (p. 224)  Bhargava  shows that there is 
a natural group structure on the  $SL(2,\mathbb Z)\times SL(2,\mathbb Z)-$moduli space of 
$2\times 2\times 2$ cubes with ``twofold symmetry'' of fixed nonzero discriminant.
We will prove the analogous  result in the context of commutative two-dimensional algebras
over a field of characteristic not two or three.\\

It is natural to consider functions invariant under  a group action when looking for coordinates
on the corresponding moduli space.
If $\FF$ is algebraically closed of characteristic not two or three, we show that
  coordinates on the $GL(2)-$moduli space of generic, two-dimensional, commutative algebras
are given  by two  ratios of ``twisted'' $GL(2)-$invariant polynomials which  are quartic in the structure constants of the algebra.
This gives a natural parametrisation of the  $GL(2)-$moduli space   by $\FF\times \FF$. A similar result was obtained
by Anan'in and Mironov in \cite{am} for { stable} two-dimensional algebras.
However, if $\FF$ is not algebraically closed, we show that   a point in  $\FF\times \FF$ does not necessarily 
characterise a class of generic algebras. Another (discrete) invariant,  the splitting field of the associated
twisted cubic, may be needed
to distinguish non-isomorphic  algebras but 
we prove that this is only necessary for equivalence classes of algebras with the property that their  moduli lie on a
degenerate elliptic curve $\Gamma_{\text{Cardano}}$ in $\FF \times \FF$ 
whose equation we give explicitly ({\it cf}  Eq.[\ref{eq:card}]).

The  parametrisation above has the advantage that many
intrinsic properties of an algebra are described by simple algebraic conditions on the corresponding
moduli. This is  a first step to understanding the global ``geometry'' of moduli space. 
For example, all equivalence classes of associative algebras 
correspond to a single point. 
This is a version of the celebrated result of B. Peirce  (completed  by his son  C. Peirce) \cite{bp} which says that
there are, up to equivalence,  seven two-dimensional  real associative algebras (not necessarily commutative).
Perhaps more surprisingly,
equivalence classes of algebras which admit non-trivial automorphisms are characterised
by the fact that their moduli lie on the curve $\Gamma_{\text{Cardano}}$ already introduced above.

The curve $\Gamma_{\text{Cardano}}$ appears naturally in yet another context:  if $\FF$ is algebraically closed,
the $SL(2)-$moduli space of generic two-dimensional commutative algebras is a principal  $\FF^\ast-$fibration over its complement in
the $GL(2)-$moduli space of generic two-dimensional commutative algebras.
In fact, we  will realise this fibration as a projective fibration by
 defining a bijective map from the $SL(2)-$moduli
space of {\it stable commutative} algebras onto a  homogeneous algebraic hypersurface in  a four-dimensional vector space.
This map is a lift of the projective embedding of the $GL(2)-$moduli space of stable two-dimensional
commutative algebras defined in \cite{am}.
\\

There is a considerable literature on the subject of two-dimensional algebras,
the vast majority of which  \cite{lu,markus,bu,ah,wal,pet,dz,di}  is concerned with finding an exhaustive  list of normal forms
 for algebras satisfying  varying constraints (commutative, division {\it etc.}).
In most of these articles the base field is $\mathbb R$ or algebraically closed but the results of  \cite{wal,pet} are
proved for more general fields.  Let us recall that the motivation for \cite{markus} was the fact that
classifying systems of two quadratic differential equations in two variables 
is equivalent to classifying  two-dimensional commutative
algebras.

More relevant to this paper is  \cite{am} which,  to the best of the authors' knowledge,
is the  only  article to study the  structure of  moduli 
space. In this paper Anan'in and Mironov  construct explicitly an embedding $\Phi$ of ${\cal M}$ into
the projective space $\mathbb P^8$, where ${\cal M}$ denotes the  $GL(2)-$moduli  space
of stable, two-dimensional (not necessarily commutative) algebras over an algebraically closed field of characteristic
not two or three. The nine affine  coordinates are given by  $SL(2)-$invariant
polynomials, quartic in the structure constants, and the authors  establish that the image of this
embedding is exactly the intersection of six quadrics.
If $\FF$ is algebraically closed there is a  link between this paper and some of our results which we now explain.
Commutative algebras which are generic  in our sense are stable in the sense of \cite{am} (the converse is
not true) so it makes sense to restrict $\Phi$ to the moduli space of generic commutative algebras.
It is easy to see that the image
of this restriction is a $\mathbb P^2\setminus \mathbb P^1$ embedded in $\mathbb P^8$ and, with respect to
a suitable identification of $\mathbb F \times \mathbb F$ with $\mathbb P^2\setminus \mathbb P^1$
their coordinates and  our coordinates are the same.
This allows us to show that, in the case of generic commutative algebras,
the  algebro-geometric singular points of ${\cal M}$ identified by
Anan'in and Mironov are precisely  the points  of the projective completion $\widetilde{\Gamma}_{\text{Cardano}}$ of
${\Gamma}_{\text{Cardano}}$.
 In other words
a generic  commutative algebra defines  a singular point of moduli space in the abstract sense if and only if
it admits non-trivial automorphisms. 
This is also true for any stable commutative algebra as shown in Appendix \ref{sec:non-gen}, and it would be very interesting
to know  whether  this property remains true for an arbitrary two-dimensional algebra.
If $\FF$ is not algebraically closed, let us recall, as mentioned above,  ``discrete'' invariants are needed to give
a complete parametrisation of the $GL(2)-$moduli space.
In Section \ref{sec:mod-sl}  we lift this projective  embedding of the
$GL(2)-$moduli space of generic commutative algebras in $\mathbb P^2$
to an  embedding of the corresponding   $SL(2)-$moduli space in a  four-dimensional vector space as the set of solutions
of an equation of  Eisenstein type.

Finally to complete our survey of the literature, let us mention that in the articles \cite{GBS, wright, wy, wood, wood2, df}
equivalence classes of twisted binary cubics are parametrised
by  Galois extensions in the case of  fields \cite{wright}  or by cubic rings in the case of  $\mathbb Z$ \cite{ wy, wood, wood2, df}.
This is relevant to the present paper since, up to equivalence, certain types of algebras are completely
determined by their associated twisted cubic.
\\

We now give a more detailed account of the results in this paper.
Let  $\FF$ be a field of characteristic not two or three and let $V$ be a two-dimensional $\FF-$vector space. 
We denote by $\M^c$ the set of commutative algebra structures on $V$, {\it i.e.}, the set of symmetric bilinear maps
from $V\times V$ to $V$.
In Sections  \ref{sec:intro} and \ref{sec:poly} we introduce several  ``twisted'' polynomial invariants 
on   $\M^c$ where 
by     ``twisted'' polynomial invariant (TPI)  we mean a $GL(V)-$equivariant polynomial map from
$\M^c$ to  some power of $\Lambda^2(V^*)$
(which from now on we denote by $L$). The decomposition 
$$\M^c \cong V^\ast \oplus\Big({\cal S}^3(V^\ast) \otimes L^{-1}\Big)\ $$
into $GL(V)-$irreducible summands has three important consequences. Firstly, there are no TPIs
of degree less than or equal to three. Secondly,  the space of ``twisted'' binary cubics ${\cal S}^3(V^\ast)\otimes L^{-1}$ appears
naturally and,  as  has already been observed in the literature \cite{wright, GBS},  the description of $GL(V)-$orbits in
this space is simpler than the description of $GL(V)-$orbits in the space ${\cal S}^3(V^\ast)$ of ``untwisted'' binary cubics.
Thirdly it allows us to define analogues of the 
 the four classical covariants of binary cubics for $\M^c$.
Maybe the most important of these is the discriminant Disc$(Q_m)$ of the  fundamental cubic $Q_m$ given by
 $Q_m(v)=m(v,v)\wedge v$ for $m\in \M^c, v\in V$.
The set of  generic algebras  $\M^c_3$ is then defined to be 
the set of  algebras $m$ such that   Disc$(Q_m)$ is non zero,
or equivalently, such that $Q_m$ has three distinct roots in a splitting field.
The other three invariants 
 $\tilde p_2,\tilde p_3$  and Inv  are respectively  quartic, quartic and sextic
 in the structure constants. It essentially follows from the fact that
 the  covariants of a binary cubic satisfy the Eisenstein identity \cite{eisen}, that these
 four invariants  also satisfy 
a kind of  Eisenstein identity.
 Furthermore, it turns out that they  also furnish  coordinates
(Section \ref{sec:moduli})
for the moduli space of generic algebras. 
In fact we have the following theorems (see Theorems \ref{theo:miracle}, \ref{theo:inv-gen1-1},
\ref{theo:inv-non-gener}, \ref{theo:inv-non-gener2},
\ref{theo:assos}, \ref{cor:aut},
Lemma \ref{lem:toto}) 
\begin{theorem}\label{th:th1.1}
Let $m \in \M^c$ (not necessarily generic). The polynomial invariants $\tilde p_3(m), \tilde p_2(m)$, \rm{Disc}$(Q_m)$ and
\rm{Inv}$_m$ satisfy the following identity in 
$ L^3$:
\beqa
27 \mathrm{Disc}(Q_m) \tilde p_3(m)^2+ 4 \tilde p_2(m)^3 = -2^4 \cdot 3^6\  \mathrm{Inv}_m^2 \ . \nn
\eeqa
\end{theorem}

\noi
Note that  this identity ({\it a fortiori} the four TPIs appearing in it) make sense 
over any unital commutative ring ({\it e.g.} the integers). In the case 
of generic $m$ we  give a  second proof of the identity  based on the fact that $\tilde p_3$ and $\tilde p_2$ turn
out to be   symmetric   functions of the roots of a ``modular cubic'' (see Corollary \ref{cor:altern}).
\begin{theorem}
If  $m\in \M_3^c$ we define $p_2(m), p_3(m) \in \FF$ by 
\beqa
p_2(m) = \frac{\tilde p_2(m)}{\mathrm{Disc}(Q_{m})}\ , \ \ \ \ 
p_3(m) = \frac{\tilde p_3(m)}{\mathrm{Disc}(Q_{m})} \  \nn
\eeqa
and  the map $\nu:\M_3^c/GL(V ) \to \FF\times \FF$ by 
\beqa
\nu(m)= (p_3(m),p_2(m)) \ . \nn
\eeqa
The Cardano curve is given by
\beqa
\label{eq:int-card}
\Gamma_{\mathrm{Cardano}} =\Big\{ (p_3,p_2) \in \FF\times \FF \ \ \text{s.t.} \ \ 
27 p_3^2 + 4 p_2^3=0 \Big\} \  , 
\eeqa
and we denote by  $\FF_{Q_m}$   a splitting field of $Q_m$.
\begin{enumerate}
\item[(i)] The map $\nu$ is surjective.
\item[(ii)]  Let $m_1,m_2 \in \M_3^c$ be such that $\nu(m_1),\nu(m_2) \not \in \Gamma_{\mathrm{Cardano}}$. Then
  there exists $g \in GL(V)$ s.t. $g\cdot m_1 =m_2$ iff 
\beqa
\nu(m_1)=\nu(m_2) \ . \nn
\eeqa
\item[(iii)]  Let $m_1,m_2 \in \M_3^c$ be such that $\nu(m_1),\nu(m_2)  \in \Gamma_{\mathrm{Cardano}}$. Then
  there exists $g \in GL(V)$ s.t. $g\cdot m_1 =m_2$ iff 
\beqa
\nu(m_1)=\nu(m_2)  \ \ \text{and}\ \ \FF_{Q_m}\cong\FF_{Q_{m'}} \ . \nn
\eeqa
\item[(iv)]  Let $m \in \M_3^c$.
\begin{enumerate}
\item If $\nu(m)  \in \Gamma_{\mathrm{Cardano}}\setminus\{0\}$ then $[\FF_{Q_m}:\FF]\le 2$. Conversely, if
$(p_3,p_2) \in \Gamma_{\mathrm{Cardano}}\setminus\{0\}$  and  $\KK/\FF$ is
an extension of degree at most two, then there exists  $m\in \M_3^c$  such that $\nu(m)=(p_3,p_2)$ and $\FF_{Q_m}\cong \KK$.
\item If $\nu(m)=(0,0)$ then $[\FF_{Q_m}:\FF]=1, 2, 3$ or $ 6$. Conversely, if  $\KK/\FF$ is a Galois extension of degree at
most six, there exists
  $m\in \M_3^c$  such that $\nu(m)=(0,0)$ and $\FF_{Q_m}\cong \KK$.
\end{enumerate}
\item[(v)] Let $m \in \M_3^c$. Then $m$ is associative iff $(p_3(m),p_2(m))=(16,-12)$. 
\item[(vi)]  The ``split'' automorphism group $\widetilde{\mathrm{Aut}}(m)$ 
(see Corollary \ref{cor:aut}) of  $m \in \M_3^c$ is given by
\beqa
\widetilde{\A}(m) \cong \left\{
   \begin{array}{cl}
    \{\bf 1\}&\text{\rm{if }}   \nu(m) \not \in \Gamma_{\text{Cardano}}\  ,\\
    \mathbb Z_2&\text{\rm{if }}  \ \ \nu(m)  \in \Gamma_{\text{Cardano}} \ \   \text{and} \ \  \nu(m) \neq (0,0)  \ ,   \\
    S_3&\text{\rm{if }}   \nu(m) = (0,0)  \ . \\
   \end{array}\right.\nn
\eeqa

\end{enumerate}
\end{theorem}
The main idea of the proof  is to find a normal form for any generic algebra, which we call 
a fundamental triple,  and  prove  its uniqueness (see Section \ref{sec:triple}).

Division algebras have always played a special r\^ole in mathematics and in  section \ref{sec:app} we obtain the following
characterisation of 
generic two-dimensional commutative division algebras in terms of our moduli.
\begin{theorem}
Let $m\in \M_3^c$.
\begin{itemize}
\item[(i)] If $\big(p_3(m),p_2(m)\big) \not \in \Gamma_{\mathrm{Cardano}}$, then 
 $m$ is 
a division algebra iff 
$$
-3(p_2(m)-p_3(m)+1)(-27p_3^2(m)- 4p_2^3(m))$$ is not a square in $\FF$.
\item[(ii)] If $\big(p_3(m),p_2(m)\big)  \in \Gamma_{\mathrm{Cardano}} \setminus\{0\}$,
write  $\FF_{Q_m}= \FF(\sqrt{a})$ where $a \in \FF$.
 Then $m$ is a division algebra iff
$$-3a(p_2(m)-p_3(m) +1)$$
  is not a square in  $\FF$.
\item [(iii)] If $\big(p_3(m),p_2(m)\big)=(0,0)$ then $m$ is  a division algebra iff $-3$ is not
a square in the splitting field $\FF_{Q_m}$  of $Q_m$.

\end{itemize}
\end{theorem}

In Section \ref{sec:mod-sl} 
 we give  moduli for the action of $SL(V)$ on generic two-dimensional commutative algebras.  
For this we introduce 
functions  invariant under $SL(V)$ but not under $GL(V)$ which distinguish algebras with the same
$GL(V)-$moduli.   In order to simplify presentation 
we state here the results only  in the case when the field 
$\FF$ is algebraically closed but  for the general case, which is considerably more complicated, see Theorem \ref{theo:slmodule}.

\begin{theorem}
Let $\FF$ be an algebraically closed field.
\begin{enumerate} 
\item  Let $\M^c_{3g}= \Big\{ m \in \M_3^c \ \ \mathrm{s.t.} \ \ \nu(m) \not \in \Gamma_{\mathrm{Cardano}}\Big\}$. The map
 $\hat \nu_s: \M^c_{3g}/SL(V)\to \Big(\FF^2 \setminus \Gamma_{\mathrm{Cardano}}\Big)
\times\Big(L\setminus\{0\}\Big)$ given by 
$$\hat \nu_s ( m) = \Bigg(\nu(m),\frac{\mathrm{Inv}_m}{\mathrm{Disc}(Q_m)}\Bigg) $$
 is a bijection.
\item  Let  $\M^c_{3C} = \Big\{ m \in  \M_3^c \ \ \mathrm{s.t.} \ \ \nu(m)  \in \Gamma_{\mathrm{Cardano}} \setminus \{0\}\Big\}$.
The map $\hat{\hat \nu}_s:\M^c_{3C}/SL(V)\to
\Big(\Gamma_{\mathrm{Cardano}}\setminus\{0\} \Big) \times \Big(L^2\setminus\{0\} \Big)$ given by
\beqa 
\hat {\hat \nu}(m)_s= (\nu(m),\mathrm{Disc}(Q_m))\nn
\eeqa 
is a bijection.
\item  Let $\M^c_{3}{}_0= \Big\{ m \in  \M_3^c \ \ \mathrm{s.t.} \ \ \nu(m) =(0,0) \Big\}$.  The map 
$\hat \nu_0 : \M^c_{3}{}_0 /SL(V)
 \to L^2\setminus\{0\}$  given by
\beqa
\hat \nu_0(m) = \text{Disc}(Q_m) \  \nn
\eeqa 
is a bijection.
\end{enumerate}
\end{theorem}

In Section \ref{sec:compact},  using Theorem \ref{th:th1.1}, we construct a surjective (bijective
if $\FF$ is algebraically closed) map from $\M^c_{\text{st}}/SL(V)$ 
 (here $\M^c_{\text{st}}$ is the set of  stable two-dimensional commutative algebras)
to  a closed hypersurface in a four-dimensional vector space. The hypersurface is the set of solutions of what
is essentially the classical Eisenstein equation considered as an equation in four variables.
The main results can be summarised as follows:
\begin{theorem} Let $\FF$ be algebraically closed and let  
\beqa
\Gamma_{\mathrm{Eisenstein}}= \Big\{ (A,B,D,C) \in L^2\times L^2 \times L^2 \times L^3 \ \ \mathrm{s.t.} \ \
27 D A^2 + 4 B^3= -2^4 3^6 C^2 \Big\} \ .\nn
\eeqa
\begin{enumerate}
\item The map $\zeta_s : \M^c_{\text{st}}/SL(V) \to \Gamma_{\mathrm{Eisenstein}}$ given by
\beqa
\zeta_s(m) = (\tilde p_3(m), \tilde p_2(m), \mathrm{Disc}(Q_m), \mathrm{Inv}_m) \ , \nn
\eeqa
is a bijection.
\item  Let  $\FF^\ast$ act on $\Gamma_{\mathrm{Eisenstein}}$ by:
\beqa
\lambda\cdot (A,B,D,C) = \Big(\frac 1 {\lambda^2} A, \frac 1 {\lambda^2} B,\frac 1 {\lambda^2} D,\frac 1 {\lambda^3} C\Big)  \  \ \ \lambda \in \FF^\ast \ , \nn
\eeqa
and  denote by $ \widetilde  {\mathbb P} (\Gamma_{\mathrm{Eisenstein}})$ the quotient of $\Gamma_{\mathrm{Eisenstein}}$ by this action.
Then $\zeta_s$ induces a bijective map $\zeta_g: \M^c_{\text{st}}/GL(V) \to \widetilde  {\mathbb P} (\Gamma_{\mathrm{Eisenstein}})$.
\end{enumerate}
\end{theorem}

\bigskip
The relationship between moduli spaces of {\it generic} two-dimensional commutative algebras and
 {\it stable} two-dimensional commutative algebras can be summarised in the following commutative cube (Theorem
  \ref{theo:c-cube} and  Appendix  \ref{sec:non-gen}).

\beqa
\xymatrix{
&\M^c_{\mathrm{st}}/SL(V)\ar[rr]^{\zeta_s}\ar@{-->}[ddd]_p  &&\Gamma_{\mathrm{Eisenstein}}\ar[ddd]^\pi \\{\color{white} toto}\\
{\color{white} toto}\\
&\M^c_{\mathrm{st}}/GL(V)\ar@{-->}[rr]^{\zeta_g}  &&\widetilde{\mathbb P}(\Gamma_{\mathrm{Eisenstein}})
\\
\M^c_3/SL(V)\ar[rr]^{\zeta_s}\ar@{^{(}->}[uuuur] \ar[ddd]_p&&\Gamma^*_{\mathrm{Eisenstein}} \ar[ddd]^\pi\ar@{^{(}->}[uuuur]
\\ {\color{white} toto}\\
{\color{white} toto}\\
\M^c_3/GL(V)\ar[rr]^{\zeta_g}\ar@{^{(}-->}[uuuur] && \widetilde{\mathbb P}(\Gamma^*_{\mathrm{Eisenstein}}) \ar@{^{(}->}[uuuur]
} \nn
\eeqa
The front face of the cube involves only  two-dimensional {\it generic} algebras, the back face only  two-dimensional
{\it stable} algebras and the maps including
 the former in the latter  can be thought of as ``compactifications''.

\bigskip
In the final section of this paper 
we consider  the  invariant
 obtained by  associating to a commutative algebra  $m$ the quadratic form $D_m$
given by the determinant of (left) multiplication.  The discriminant  of $D_m$ is not just $SL(V)-$invariant but is
invariant under the action of a larger group  $SL(V)\times SL(V)$ in which it is diagonally embedded.
We show that this association establishes a
bijection between 
the set of
commutative algebra structures
 of fixed non-vanishing discriminant $\Delta$ modulo the action of $SL(V)\times SL(V)$, and the
set of quadratic forms of discriminant $\Delta$ modulo the action of $SL(V)$.
By pullback of the Gauss composition law this allows us to endow the set of $SL(V)\times SL(V)-$orbits of fixed discriminant
$\Delta$ with a group law. In the context of commutative two-dimensional algebras over a field of characteristic not
two or three, this result is analogous to a result of Bhargava  \cite{bhar}  who
exhibited a natural group structure on the $SL(2,\mathbb Z)\times SL(2,\mathbb Z)$ moduli space of
$2\times 2 \times 2$ cubes with ``twofold'' symmetry.

\section{Notations and definitions}\label{sec:intro}
Throughout this paper:

\begin{itemize}\setlength {\itemsep }{0.1 mm }
\item[-] $\FF$ is a field of 
characteristic not two or three; 
\item[-] $V$ is  a two-dimensional vector space  over $\mathbb F$;
\item[-]  $V^\ast$ is the  dual of $V$;
\item[-] ${\cal S}^2(V^\ast)$ is the set of symmetric bilinear forms on $V$;
\item[-] $L$ is the set of anti-symmetric bilinear forms on $V$ (this is one-dimensional);
\item[-] $L^{-1}$ is  the dual of $L$;
\item[-]   $L^r \otimes L^s$  is always identified  with  $L^{r + s }$ for   any integers $r,s$.
\end{itemize}

\begin{definition}
\label{def:L-R}
\phantom{marcus} \hskip 3.truecm
\begin{enumerate}
\item Let $\M^c= {\cal S}^2(V^\star) \otimes V$. We think of this as   the set of commutative algebra structures on $V$,
 that is, the set of bilinear maps $m$ from $V\times V$ to $V$ such that 
for all $v,w \in V, m(v,w)=m(w,v)$.
\item Let  ${\cal F} = \{ m \in \M^c \ \text{s.t.} \ m \mathrm{~is~a~field~or~} m \mathrm{~is~a~ split~field}\}$ \
(by `` $m$ is a split field'' we mean that  $m$  has two idempotents $e_1,e_2$ and an identity $1$:
$m(e_1,e_1)=e_1, m(e_2,e_2)=e_2, m(e_1, e_2) =0$ and $1=e_1+e_2$).
\end{enumerate}
\end{definition}

As usual, the left  multiplication operation
$L_v : V \to V$ 
is
\beqa
\label{eq:L-R}
 L_v(m)(w) := m(v,w)   \ \ \forall v,w \in V, m \in \M^c \ .  
\eeqa

\noi
The group $GL(V)$ acts  on $V$ by
\beqa
\label{eq:GLV}
g \cdot m(v,w) = g\Big(m\big(g^{-1}(v),g^{-1}(w)\big)\Big) \ ,\ \ 
\forall v,w \in V, m \in \M^c  \ .
\eeqa

In 1881 Benjamin Peirce introduced the notion of ``idempotent'' \cite{bp} in order to classify two-dimensional associative
 (not necessarily commutative) real algebras.  In fact, one can associate to  any two-dimenional algebra a (twisted) binary cubic on $V$ whose roots correspond to its idempotents (see also for example \cite{dz}.

\begin{definition} \label{def:cubic}
Let $m$ be in $\M^c$. The fundamental binary cubic of $m$ is the function 
$Q_m : V \to L^{-1}$ defined for all $v \in V$ by
\beqa
\label{eq:Q}
Q_m(v):=  m(v,v)\wedge v \ .   
\eeqa
An element $v \in V$ is an idempotent of $m$ iff
\beqa
\label{eq:idem}
Q_m(v) = 0 \ .
\eeqa
\end{definition}

\section{Twisted polynomial invariants and  the twisted Eisenstein identity}\label{sec:poly}

In this section, based on the decomposition of $\M^c$ into $GL(V)-$irreducible summands and inspired by  the classical
theory of covariants of binary cubics, we introduce four invariant polynomial maps from $\M^c$ to powers of $L$.
These four maps satisfy an identity analogous to the classical Eisenstein identify satisfied
by the four classical covariants of a binary cubic. Furthermore, as we shall see in Section \ref{sec:mod-sl},
they also define  an embedding (at least if $\FF$ is algebraically closed) of the $SL(V)-$moduli space of stable two-dimensional commutative algebras onto
a hypersurface in a four-dimensional vector space.
This embedding  turns out to be a lift of the projective embedding of the $GL(V)-$moduli space of stable
two-dimensional commutative algebras 
defined in \cite{am}.

It has been known for a long time that one can associate three covariants to a homogeneous cubic polynomial in two variables
\cite{eisen}.
These are polynomial functions of the cubic  of  which perhaps 
the most important is the discriminant,  and we will now generalise them
to the case of ``twisted'' cubics.

We begin with the discriminant. 
Recall that the classical discriminant of  $P(x,y)=Ax^3 + B x^2 y +C x y^2 + D y^3$ 
is given by
\beqa
P \mapsto 18 ABCD + B^2 C^2 -4 AC^3 -4B^3 D -27 A^2 D^2 \ . \nn
\eeqa
This formula defines a quartic function of $P$ which is invariant under  special linear transformations
but not under a general linear transformation. In order to define a $GL(V)-$invariant ``discriminant'' we
need to consider quartic functions of the cubic $P$ with values in powers  of 
the one dimensional vector space $L$. 
Since  $\lambda$ Id acts on cubic polynomials by  $1/\lambda^3$ and on  $L^{-1}$  by $\lambda^2$, it follows that
the formula
\beqa
\text{Disc}_c(P) = (18 ABCD + B^2 C^2 -4 AC^3 -4B^3 D -27 A^2 D^2)(\epsilon_1 \wedge \epsilon_2)^{\otimes^6}  \nn
\eeqa 
defines a $GL(V)-$equivariant map  Disc$_c: {\cal S}^3(V^\ast) \to L^{6}$, where $\{\epsilon_1,\epsilon_2\}$ is
the dual basis of $\{e_1,e_2\}$.
 Motivated by this definition of the classical discriminant of a cubic  we define 
Disc$:{\cal S}^3(V^\ast)\otimes L^{-1} \to L^2$ by
\beqa
&\hskip -3.truecm \text{Disc}\Big((Ax^3 + B x^2 y +C x y^2 + D y^3)(e_1\wedge e_2) \Big)= \nn\\
&\hskip 1.truecm \big(18 ABCD + B^2 C^2 -4 AC^3 -4B^3 D -27 A^2 D^2\big)
(\epsilon_1\wedge \epsilon_2)^2 \ . \nn
\eeqa 

If $m \in \M^c$, we now give the  explicit formul\ae \ for $Q_m$ and Disc$(Q_m)$ in terms of the structure constants of $m$.
Let  $(e_1,e_2)$ be a basis of $V$ with dual  basis  $(\varepsilon_1,\varepsilon_2)$ and $m: V\times V \to V$ be given by:
\beqa\label{rem:surjQ}
m(e_1,e_2) = a e_1 + b e_2\ , \ \
m(e_2,e_2)=  c e_1 + d e_2 \ , \nn \\\ \
m(e_1,e_2) = e e_1 + f e_2 \ , \ \ 
m(e_2,e_1) = e e_1 + f e_2 \ . 
 \eeqa
Then  $Q_m$  and  Disc$(Q_m)$ are given by
\beqa
\label{eq:decom}
Q_m&=&\Big( -b x^3 + (a -2f) x^2 y + (2e -d) x y^2 + c y^3\Big)(e_1 \wedge e_2) \ , \\
\text{Disc}(Q_m) &=&\Big(-36bcfd+72bcfe+18bcad-36bcae+4f^2d^2-16f^2de+16f^2e^2 -4fad^2\nn \\
&&\hskip .4truecm+16fade-16fae^2+a^2d^2-4a^2de+4a^2e^2-4bd^3+
24bd^2e -48bde^2 +32be^3\nn \\
&&\hskip .4truecm +32cf^3-48cf^2a+24cfa^2-4ca^3-27b^2c^2\ \Big)
 (\varepsilon_1 \wedge \varepsilon_2)\otimes 
(\varepsilon_1 \wedge \varepsilon_2)\ . \nn
\eeqa
It follows from this that given any binary cubic $G$ there exists  $m$ in $\M^c$
such that $Q_m = G e_1 \wedge e_2$. 
Recall also that $G$ is irreducible {\it iff} its splitting field is of degree three or six, and in
this case  the degree is three  {\it iff}
Disc$(G)$ is a square.

\begin{remark} \label{rem:am-disc}
In \cite{am} Anan'in and Mironov (p 4483) attach to any  
two-dimensional algebra  nine ``projective coordinates'':
\beqa
[x_1,x_2,x_3,y_1,y_2,y_3,z_0,z_1,z_2] \ .\nn 
\eeqa
These are all quartic functions of the structure constants and
if the algebra is commutative only $x_1,x_2,z_0$ do not vanish identically.
 One can check that  the link between our invariant Disc$(Q)$ and $z_0$ is  Disc$(Q)=81 z_0 (\epsilon_1 \wedge \epsilon_2)^2$  where
$(\epsilon_1, \epsilon_2)$ is the basis dual to
the basis $(e_1,e_2)$  given by their Eq.[6] (p  4483). 
\end{remark}

The normalisation of the discriminant above  has been chosen so that if $m$ 
is a split field then Disc$(Q_m) \in L^2$ is a perfect
square. In fact using Disc, we have a natural parametrisation of field structures on $V$ up to $SL(V)-$equivalence which extends the
well-known natural parametrisation  by $\FF^*/\FF^{*2}$ of field structures on $V$ up to $GL(V)-$equivalence:

\begin{proposition}

Recall that  ${\cal F}=\{ m \in \M^c \ \text{s.t.} \ m \mathrm{~is~a~field~or~} m \mathrm{~is~a~ split~field}\}$.
Let:
\begin{itemize}
\item  $s: L^2\setminus\{0\}\to \FF^*/\FF^{*2}$  be defined by
\beqa
s\big(\lambda(\omega_0\otimes \omega_0)\big) =[\lambda] \ , \nn
\eeqa
for $\omega_0$ any in $ L$. This is clearly independent of the choice of $\omega_0$.
\item  $\delta: {\cal F}/GL(V) \to \FF^*/\FF^{*2}$ be defined by
\beqa
\delta([m]) = [z^2] \ ,\nn
\eeqa
where $z$ is any (nonzero) purely imaginary element of $m$.
\item  $p: {\cal F}/SL(V)\to {\cal F}/GL(V)$ be the natural projection.
\end{itemize}

\noi
Then the diagram
\beqa
\xymatrix{
{\cal F}/SL(V) \ar[r]^{\mathrm{Disc}} \ar[d]_p& L^2\setminus\{0\}\ar[d]^s\\
{\cal F}/GL(V) \ar[r]^{\delta} &\FF^*/\FF^{*2}
}\nn
\eeqa
 is commutative  and the maps Disc and $\delta$ are bijections.
\end{proposition}

\begin{demo}
The fact that $\delta$ is a bijection is well known (see {\it e.g.} \cite{serre}) and 
one can check easily that the diagram is commutative. 

To prove the surjectivity of Disc, let $x \in L^2\setminus\{0\}$. Then by the surjectivity of $\delta$ there exists
$m \in {\cal F}$ such that $\delta([m]_{GL(V)})=[x]$.  Hence by the definition of $s$ there exists $\lambda \in \FF^*$ such that 
\beqa
\mathrm{Disc}([m]_{SL(V)}) =\lambda^2 x \  \nn
\eeqa 
and from this  it follows that 
\beqa
\mathrm{Disc}(\frac1 \lambda Id \cdot m)= x \ .\nn
\eeqa
This prove that the map Disc is surjective.

To prove the injectivity of Disc, let $m_1, m_2 \in {\cal F}$ be such that
\beqa
\mathrm{Disc}([m_1]_{SL(V)})= \mathrm{Disc}([m_2]_{SL(V)}) \ .\nn
\eeqa
Then
\beqa
\delta([m_1]_{GL(V)})= \delta([m_2]_{GL(V)}) \ , \nn
\eeqa
and by the injectivity of $\delta$ there exists $g \in GL(V)$ such that $m_1 = g\cdot m_2$. Thus
\beqa
\mathrm{Disc}(m_1) = \frac 1 {\det(g)^2} \mathrm{Disc}(m_2)  \  \nn
\eeqa
from which it follows that $\det (g) = \pm 1$.  However, we can always assume that  $\det (g) = 1$ since if not, we compose
with the non-trivial element of the Galois group which is always of determinant equal to minus one. 
This prove injectivity.
\end{demo}

\begin{remark}
Our convention is the following: if   $(e_1,e_2)$ is a basis of $V$ with dual  basis  $(\varepsilon_1,\varepsilon_2)$ 
and $v_i= x_i e_1 + y_i e_2$ are vectors in $V$ then 
\beqa
(\varepsilon_1 \wedge \varepsilon_2)(v_1,v_2) = (x_1 y_2 -x_2 y_1) \ . \nn
\eeqa

\end{remark}

We  now give the decomposition of $\M^c$  into its  $GL(V)-$irreducible components
and for this we need the following
\begin{definition}
\label{def:KH}
Let $m \in \M^c$ and let  $L_v(m)$  be left multiplication
as in  \eqref{eq:L-R}. Then 
$T_{m} \in V^\ast$ and $D_m \in {\cal S}^2(V^\ast)$ are defined by
\beqa
\label{eq:DT}
\begin{array}{ll}
T_m(v) = \mathrm{Tr}(L_v(m)) \ , &
D_m(v) = \det(L_v(m)) \  , \ \forall v \in V \ .  
\end{array}
\eeqa
In terms of the structure constants of $m$  (see \eqref{rem:surjQ}) this gives the formu\ae,
\beqa
\label{eq:TD}
T_m(x e_1 + y e_2)&=& (a+f) x + (e+d)y \ , \nn\\
 D_m(x e_1 + y e_2)&=&(af-be)x^2+(ad-bc)xy+(ed-cf)y^2 \ . 
\eeqa
\end{definition}
\begin{proposition}\label{prop:Ac-irred} 
\phantom{titi} \hskip 1.truecm
\begin{enumerate}
\item As $GL(V)-$representations we have: $\M^c \cong V^\ast \oplus \Big({\cal S}^3(V^\ast) \otimes L^{-1}\Big)$.
\item For $m \in \M^c$  define $m' \in \M^c$ by
\beqa
m'(v,w) = \frac13 \Big(T_m(v) w + T_m(w)v\Big) \ . \nn
\eeqa
Then 
\beqa
m = m' + (m-m') \ , \nn
\eeqa
is the irreducible decomposition of $m$.

\end{enumerate}
\end{proposition}
\begin{demo}
(1): Let $K= \{m \in \M^c \ \ \text{s.t.}\ \  Q_m\equiv 0 \}$ and consider the exact sequence of $GL(V)-$representations:
\beqa
\{0\} \longrightarrow K \longrightarrow \M^c \longrightarrow  {\cal S}^3(V^\ast) \otimes L^{-1} \longrightarrow \{0\} \ .\nn
\eeqa
Since the map $Q$ is not identically zero and $ {\cal S}^3(V^\ast) \otimes L^{-1}$ is irreducible, it follows that $Q$ is surjective
and  that $K$ is a two-dimensional representation of $GL(V)$. Define the $GL(V)-$equivariant map $k:V^\ast \to \M^c$ by
\beqa
k_\alpha(v,w)= \alpha(v) w + \alpha(w) v \ , \ \forall \alpha \in V^*, \forall u,v \in V \ . \nn
\eeqa
Since  $Q_{k_\alpha}(v)= 2 \alpha(v) v \wedge v =0$, we have Im$(k_\alpha) \subset K$ and,
since $V^\ast$ is irreducible of the same dimension as $K$, the map $k$ is an isomorphism of $GL(V)-$representations.

(2): For all $m \in \M^c$, it is evident that $m'= k_{\frac13 T_m}$ and easily checked that $T_{m'}=T_m$.
\end{demo}

\begin{remark}
The space of twisted binary cubics $B={\cal S}^3(V^\ast) \otimes L^{-1}$ appearing in the above decomposition has already
been studied in the literature. In both \cite{GBS,wright} it was explicitly pointed out that
the crucial  advantage of $B$ over  the space of {\it untwisted} binary cubics ${\cal S}^3(V^\ast)$ is that
$\lambda\cdot$Id$_V$ acts by  $\frac1 \lambda$Id$_B$ instead of by $\frac1 {\lambda^3}$Id$_{{\cal S}^3(V^\ast)}$
(recall that  in \cite{GBS} the authors are working over $\mathbb Z$ whereas in \cite{wright} the author is
working over over a number field).
\end{remark}

\begin{remark}
Let $q \in {\cal S}^3(V^*) \otimes L^{-1}$ be a ``twisted'' binary cubic.
It follows from the proposition that there
  exists a unique  algebra $m$ such that $Q_m=q$ and $T_m=0$. Explicitly,
if   $(e_1,e_2)$ is a basis of $V$ and  $q(x e_1 + y e_2)= (\alpha x^3 + \beta x^2 y + \gamma x y^2 + \delta y^3) e_1 \wedge e_2$
in this basis, it is easy to see that:
\beqa
m(e_1,e_1) = \frac13\beta  e_1 -\alpha e_2 \ , \ \
m(e_2,e_2) = \delta e_1 -\frac13 \gamma  e_2 \ , \ \
m(e_1,e_2) =  \frac13 \gamma e_1 -\frac13 \beta e_2 \ . \nn
\eeqa
\end{remark}

In the rest of this paper an important r\^ole will be played by twisted polynomial invariants (TPI)  on $\M^c$.
A basic example of a TPI   is the quartic  polynomial   Disc$: \M^c \to L^2$  and this motivates the following definition:
\begin{definition}
A twisted polynomial invariant (TPI) of degree $(2k,k)$  on $\M^c$   is a degree $2k$ polynomial map $P : \M^c \mapsto  L^k$ which is 
$GL(V)-$equivariant, {\it i.e.}, $P$ is an element   of ${\cal S}^{2k}(\M^{\ast c}) \otimes L^k$ which is fixed by the action
of $GL(V)$. 
\end{definition}

Using Proposition \ref{prop:Ac-irred}  we can adapt the classical theory of polynomial invariants of binary cubics to construct
TPIs as follows. Setting
$B=  {\cal S}^3(V^\ast) \otimes L^{-1}$,
the decomposition $\M^c = V^\ast \oplus  B$  implies the $GL(V)-$equivariant decomposition of degree $(2k,k)$ polynomials on
$\M^c$: 
\beqa
\label{eq:bideg}
{\cal S}^{2k}(\M^{ c *}) \otimes L^k  = \bigoplus\limits_{m+n=2k} {\cal S}^m(V) \otimes  {\cal S}^n(B^\ast) \otimes L^k \ . \nn
\eeqa
Recall that there is an $GL(V)-$equivariant isomorphism
\beqa
V\cong V^\ast \otimes L^{-1}  \nn
\eeqa
given by: 
\beqa
\xymatrix{
  i_\alpha b&\ar@{|->}[l] \alpha \otimes b 
} \ \ \ \ \forall \alpha \in V^\ast \ , \ b \in L^{-1}  \ ,
\nn
\eeqa
which induces    an $GL(V)-$equivariant isomorphisms
\beqa
 {\cal S}^m(V)  \otimes  {\cal S}^n(B^\ast) \otimes L^k \cong {\cal S}^m(V^\ast) \otimes L^{k-m}  \otimes  {\cal S}^n(B^\ast)  \ . \nn
\eeqa
and
\beqa
{\cal S}^m(V)  \otimes  {\cal S}^n(B^\ast) \otimes L^p \cong {\cal S}^m(V^\ast) \otimes L^{p-m}  \otimes  {\cal S}^n(B^\ast)  \ .\nn
\eeqa     
Hence TPIs of degree $(2k,k)$ are described by the $GL(V)-$equivariant isomorphism:
\beqa
\label{prop:TPI}
{\cal S}^{2k}(\M^{ c *}) \otimes L^k  = \bigoplus\limits_{m+n=2k} {\cal S}^m(V^*) \otimes {\cal S}^n(B^\ast)  \otimes L^{k-m} \ . 
\eeqa

Now we recall that the algebra of  classical covariants  of  the space of binary cubics ${\cal S}^3(V^\ast)$ is the algebra of
$GL(V)-$equivariant 
polynomial maps $f: {\cal S}^3(V^\ast) \to {\cal S}(V^\ast)$, {\it i.e.}, the $GL(V)-$invariant elements of
\beqa
{\cal S}\big(({\cal S}^3(V^\ast))^\ast \otimes {\cal S}(V^\ast) \big) =
{\cal S}(B^\ast) \otimes {\cal S}(L^{-1})  \otimes {\cal S}(V^\ast)  \ .
\nn
\eeqa
Hence by  \eqref{prop:TPI} each classical covariant suitably adapted to our situation
 will give rise to a TPI. For this we need the following definition:

\begin{definition}
Let $E$ be a vector space,
let $\alpha \in V^\ast$ and let $P \in {\cal S}^k(V^\ast) \otimes E$. Then we denote by $P(\alpha)$ the element of 
$L^k \otimes E$ defined by:
\beqa
\label{eq:P}
P(\alpha)=P(v)\otimes \omega^k \ , 
\eeqa
where $v\in V$ and $\omega \in L$ are chosen so that $\alpha= i_v \omega$. One easily checks that this definition is
independent of the choice of such  pair $(v,\omega)$.
\end{definition}

Let us  illustrate this with: $\alpha= T_m$, $P=f(Q_m)$ where $m\in\M^c$ and $f :  B \mapsto S^k(V^\ast)\otimes L^p$ 
is a polynomial   function. 
 In the basis of  \eqref{rem:surjQ} we have
\beqa
T_m(x e_1 + y e_2) = (a+f) x + (e+d)y \ , \nn
\eeqa  
in other words,
\beqa
T_m = (a+f)\epsilon_1 + (e+d) \epsilon_2  \ . \nn
\eeqa
Hence 
\beqa
T_m = i_{\big((e+d)e_1 - (a+f) e_2\big)} \epsilon_1 \wedge \epsilon_ 2 \ , \nn
\eeqa
and
\beqa
\label{eq:TPI-bas}
f(Q_m)(T_m) = f(Q_m)\big((e+d)e_1 - (a+f) e_2\big) \otimes (\epsilon_1 \wedge \epsilon_ 2)^{\text{deg}(f)} \ . 
\eeqa
Note that  the function $f(Q_m)(T_m)$ is of bidegree $(k,\text{deg}(f))$ (see \eqref{eq:bideg}).
We now give explicit formul\ae \ for the TPIs on $\M^c$ which are obtained by taking for $f$ 
in \eqref{eq:TPI-bas}
 the classical covariants of a binary cubic \cite{eisen}.

The simplest example is obtained by taking  for $f$ the covariant
Disc$:B \to L^2$ (corresponding to $n=4,m=0,k=2$ in \eqref{prop:TPI}) which gives rise to Disc$(Q_m)(T_m) : \M^c \to L^2$ of bidegree $(4,2)$ (see \eqref{eq:decom}
for the explicit formula).

The second example is  obtained by taking  for $f$  the covariant Id$: B \to S^3(V^\ast)\otimes L^{-1}$
(corresponding to $m=3,n=1,k=2$ in \eqref{prop:TPI}):
\begin{definition}
If $f : B \to S^3(V^\ast)\otimes L^{-1}$  is the identity then we define $\tilde p_3: \M^c \to L^2$ 
of bidegree $(4,2)$
by 
\beqa
\tilde p_3(m)=- 8 Q_m(T_m) \ .\nn
\eeqa
In the basis   \eqref{rem:surjQ} we have 
\beqa
\label{eq:pt3}
  \tilde p_3(m)&=&\Big(-32fade+24bde^2+24bd^2e+24cfa^2+24cf^2a +8a^2de\nn\\
&& \hskip .34truecm-40fae^2+8fad^2 -40f^2de-8f^2d^2-32f^2e^2+16a^2d^2 \\
&&\hskip .34truecm -8a^2e^2+8bd^3+8be^3+
8cf^3+8ca^3 \Big) (\epsilon_1 \wedge \epsilon_2)^2 \ . \nn
\eeqa

\end{definition}

\begin{remark}\label{rem:am-p3}
In terms of the  coordinates of \cite{am}  we have
$\tilde p_3 = 81 x_2 (\epsilon_1\wedge \epsilon_2)^2$ where
$(\epsilon_1, \epsilon_2)$ is the basis dual to
the basis $(e_1,e_2)$  given by their Eq.[6] (p  4483). 
\end{remark}

We now recall the definition of the the two remaining classical  covariants of a twisted binary cubic $Q \in B$.

\newpage
\begin{definition}
Let   $(e_1,e_2)$ be a basis of $V$,
let $x,y$ be the coordinates of $v\in V$ in this basis and 
 let $Q =Q_1 e_1 \wedge e_2\in B$. 
\begin{enumerate}
\item The map $\mu: B \to S^2(V^\ast)$ is  defined  by taking the Hessian of $Q_1$:
\beqa
\label{eq:hess}
\mu(Q) = \frac{\partial^2 Q_1}{\partial x^2} \frac{\partial^2 Q_1}{\partial y^2}- 
\Big(\frac{\partial^2 Q_1}{\partial x \partial y}\Big)^2 \ . \nn
\eeqa
One checks that the quadratic form  $\mu(Q)$ does not depend on the choice of $(e_1,e_2)$.
In the basis of  \eqref{rem:surjQ} this gives
\beqa
\label{eq:mu}
\mu(Q) (x e_1 + y e_2)&=&\phantom{+}(-16f^2+16fa-4a^2+12bd-24eb)x^2\nn \\
&&+(-36bc-8fd+4ad-8ae+16fe)yx  \ \nn \\ 
&&+ (16de-16e^2+12ac-24cf-4d^2)y^2 \ , 
\eeqa
\item The map $G:B\to {\cal S}^3(V^\ast)$ is defined by taking   the Poisson bracket of $Q_1$ and $\mu$:
\beqa
G(Q)= \frac{\partial Q_1}{\partial x} \frac {\partial \mu}{\partial y} -  \frac{\partial Q_1}{\partial y} \frac{\partial \mu}{\partial x} \ . \nn
\eeqa
One checks that the binary cubic  $G(Q)$ does not depend on the choice of $(e_1,e_2)$.
In the basis of  \eqref{rem:surjQ} this gives
\beqa
\label{eq:G}
&&G(Q)(x e_1 + y e_2)=\nn\\
&&\phantom{+}\Big(-144bfe+72bae+96f^2a+72bfd-36bad-48fa^2+108b^2c-64f^3+8a^3\Big)x^3\nn\\
&&+\Big(48fad-288bde+96f^2e+288be^2-12a^2d+24a^2e\nn\\
&&\hskip 2.truecm -108bac+216bcf-48f^2d-96fae+72bd^2\Big)yx^2\\
&&+\Big(72a^2c-12ad^2+48ade-48ae^2+96fe^2-96fde\nn\\
&&\hskip 2.truecm -288fac+288cf^2+24fd^2-108bcd+216bce\Big)y^2x\nn\\
&&+
\Big(72dcf+72eac-144ecf-36dac+96de^2+108c^2b+8d^3-64e^3-48d^2e\Big)y^3\nn
\eeqa
\end{enumerate}
\end{definition}

Note that these covariants are invariant under the action of $GL(V)$ whereas the classical covariants of (untwisted) binary cubics
are only invariant under the action of $SL(V)$. This means that that we have to slightly adapt the Eisenstein identity
satisfied by the classical covariants to our situation.
\begin{proposition}
\label{prop:eisen}
The values of the covariants  Disc, Id,  $\mu$, $G$ at any $v$ in $V$  satisfy the Eisenstein identity
\beqa
16 \times 27 \mathrm{Disc}(Q) Q(v)^2 + \mu(Q)(v)^3= - G(Q)(v)^2 \  . \nn
\eeqa
\end{proposition}

\begin{demo}
A straightforward consequence of the classical Eiseinstein equation \cite{eisen}.
\end{demo}

Taking for $f$ in \eqref{eq:TPI-bas} the covariants $\mu$ and $G$ we obtain 
 two new TPIs.
\begin{definition}
\phantom{Quaspar} \hskip 1.truecm
\begin{enumerate}
\item  If $f: B \to S^2(V^*)$ is the covariant $\mu$ we define $\tilde p_2 : \M^c \to {\cal S}^2(V^\ast)$ of bidegree $(4,2)$ by
\beqa
\tilde p_2(m) =\mu(Q_m)(T_m) \ . \nn
\eeqa
In the basis of \eqref{rem:surjQ} this gives
\beqa
\label{eq:pt2}
\tilde p_2(m)  &=&\Big(36bcfd+36bcfe+36bcad+36bcae+60fade-36bde^2\nn\\
&&\hskip -.1truecm -36cf^2a+12a^2de -24fae^2+12fad^2-24f^2de-12f^2d^2\\
&&\hskip -.1truecm  -48f^2e^2-12a^2d^2-12a^2e^2+12bd^3-24be^3
-24cf^3 +12ca^3 \Big) (\epsilon_1 \wedge \epsilon_2)^2  \ . \nn
\eeqa
\item  If $f: B \to S^6(V^*)$ is the covariant $G$ we define Inv$: \M^c \to{\cal S}^3(V^\ast)$ of bidegree $(6,3)$ by
\beqa
\mathrm{Inv}_m = \frac1 {54} G(Q_m)(T_m) \ . \nn
\eeqa
In the basis of  \eqref{rem:surjQ} this gives
\beqa
\label{eq:Invd}
\mathrm{Inv}_m&=&\Big(
-2a^3cdf-4a^3cef-6a^2bc^2f+6b^2cde^2-6a^2de^2f\nn\\
&& \hskip .4truecm +4bd^3ef -12bde^3f+12acef^3-8de^2f^3+6a^2bce^2\nn\\
&& \hskip .4truecm +2acdf^3+4a^2e^3f-2abd^4-4abe^4+
2ad^3f^2+8ae^3f^2\nn\\
&& \hskip .4truecm +2b^2cd^3 +2b^2ce^3-2bc^2f^3 -8be^4f+4cdf^4+8cef^4 \nn\\
&& \hskip .4truecm -4d^2ef^3-2abde^3+6abd^2e^2-2a^2d^3f -2a^3de^2
+2a^4cd\nn\\
&& \hskip .4truecm   -2a^3bc^2+2a^3d^2e+6b^2cd^2e+2abd^3e-6bcd^2f^2+6ad^2ef^2 \\
&&\hskip .4truecm -6abc^2f^2 -6bcdef^2+6abce^2f-6abcd^2f+6a^2bcde-6a^2cdf^2\Big)\times  (\epsilon_1 \wedge \epsilon_2)^3 \ . \nn
\eeqa
\end{enumerate}
\end{definition}

\begin{remark} \label{rem:am-p2}
In terms of the  coordinates of \cite{am}  we have
$\tilde p_2 = 81 x_1 (\epsilon_1\wedge \epsilon_2)^2$ where
$(\epsilon_1, \epsilon_2)$ is the basis dual to
the basis $(e_1,e_2)$  given by their Eq.[6] (p  4483). 
\end{remark}

\begin{remark}
In contrast to Disc$(Q), \tilde p_3$ and $\tilde p_2$ which are quartic in the structure constants and
correspond to  coordinates in \cite{am}, the invariant Inv is  sextic
  and is not one of the  coordinates in \cite{am}. It turns out that these four invariants
satisfy a twisted  Eisentein identity (see below) which 
is essential to our
 parametrisation of the  $SL(V)-$moduli space of commutative algebras in Section \ref{sec:mod-sl}.  
\end{remark}

We now prove the  twisted Eisenstein identity mentioned above. This will be a consequence of 
Proposition \ref{prop:eisen}.

\begin{theorem}\label{theo:miracle}
Let $m \in \M^c$.  Then as elements of $\Lambda^2(V^*)^6$ we have the identity:
\beqa
\begin{array}{ll}
27 \mathrm{Disc}(Q_m) \tilde p_3(m)^2+ 4 \tilde p_2(m)^3 = -(4 \times 27)^2\  \mathrm{Inv}_m^2 \ .
\end{array}
\eeqa

\end{theorem}

\begin{demo}
Since
\beqa
\tilde p_3(m) &=& -8Q_m(T_m) \ , \nn \\
\tilde p_2(m) &=& \mu(Q_m)(T_m) \ , \nn \\ 
\mathrm{Inv}_m&=& \frac1{54} G(Q_m)(T_m) \ ,  
 \nn
\eeqa
the theorem follows from Proposition \ref{prop:eisen}.

\end{demo}

To conclude this section we introduce one more  TPI on $\M^c$  which, although it can
 be expressed in terms of
the polynomial invariants above (Proposition  \ref{prop:miracle}),  is of interest in its own right and will be useful later.
\begin{definition}\label{def:D-bha}           
Let $m \in \M^c$ and let  $L_v$ be the left multiplication by $v \in V$.  Define the quartic twisted polynomial invariant
 \rm{Disc}$\circ D:
\M^c \to L^2$ by  
\beqa
\mathrm{Disc}\circ D_m(v) = \mathrm{Disc}(\det(L_v)) \ . \nn
\eeqa
In the basis of   \eqref{rem:surjQ}  we have 
\beqa
\label{eq:dD}
\mathrm{Disc}(D_m)=
&&\Big(a^2d^2-2bcad+b^2c^2-4fade+4cf^2a+4bde^2-4bcfe\Big) (\epsilon_1 \wedge \epsilon_2)^2 
\eeqa 
\end{definition}

The relationship between Disc$(D)$ and the TPIs  appearing in the twisted Eisenstein identity is given by
the following proposition for which we omit the proof.
\begin{proposition}
\label{prop:miracle}
Let $m \in \M^c$. Then as elements of $L^2$ we have the identity
\beqa
27 \mathrm{Disc}(D_m)= \tilde p_3(m) - \tilde p_2(m) - \mathrm{Disc}(Q_m) \ . \nn
\eeqa
\end{proposition}

\begin{remark}
All of the TPIs appearing in the  identities of Theorem \ref{theo:miracle} and Proposition \ref{prop:miracle}
are defined over $\mathbb Z$ and hence the identities
themselves are valid over $\mathbb Z$.
\end{remark}

Finally let us point out that one can always generate ``new'' TPIs from ``old'' TPIs
 by applying to them
 natural operations  on $\M = V^* \otimes V^* \otimes V$, {\it i.e.}, elements 
 of $GL(\M)$ which commute with the action of
$GL(V)$.  For example in this way, the invariant $D$  will generate  three  (in general)
linearly independent invariants, say $D, D',D''$.
These  are essentially the three  quadratic forms introduced
by Bhargava in \cite{bhar}. 
If  the multiplication is commutative it turns out there are only three possiblities:   (i) $D=D'=D''$; (ii) $D=D'=-D''$; (iii)
$D=D'$ and $D''$ independent of $D$.
Explicitly, in the basis of  \eqref{rem:surjQ} with $v = x e_1 + y e_2$, we have the following formul\ae
\beqa
D_m(v) &=& (af-be)x^2+(ad-bc)xy+(ed-cf)y^2 \ , \nn\\
D''_m(v)&=& (-f^2+bd)x^2+(2ef -ad-bc)xy+(-e^2+ac)y^2 \ . \nn
\eeqa
Notice that if $T_m=0$ $a+f=d+e=0$ and hence $D_m=D''_m$.
Bhargava observed
an analogous property for cubes with ``twofold symmetry'' or ``threefold symmetry'' \cite{bhar}, p. 225.

In a forthcoming paper we will investigate invariants of  arbitrary  two-dimensional algebras in more detail.

\section{Moduli space of generic algebras  with respect to the action of $GL(V)$  }\label{sec:moduli}
\subsection{R\'esum\'e and comparison with  the results of Anan'in and Mironov}\label{sec:res}
In this section we consider the set of generic commutative algebras 
\beqa
\M^c_3 = \Big\{ m \in \M^c\  \text{s.t.} \ \  \text{Disc}(Q_m) \ne 0 \Big\} \ , \nn
\eeqa
the set of stable commutative algebras \cite{am}
\beqa
\label{eq:stable}
\M^c_{\text{st}} = \Big\{ m \in \M^c \ \ \text{s.t.} \  \exists \ \ SL(V)-\mbox{invariant polynomial}\ \ F \not \equiv 0 \ \ \mbox{on}\ \ \M^c \  \text{s.t.} \ \ F(m) \ne 0 \Big\} 
\eeqa\
and especially their moduli spaces
\beqa
{\cal M}_3^c = \M^c_3 /GL(V) \ , \ \ 
{\cal M}_{\text{st}}^c = \M^c_{\text{st}}/GL(V) \ . \nn
\eeqa
Commutative algebras which are generic are stable  but the converse is not true (Appendix \ref{sec:non-gen}).
Recall that $\tilde p_2(m),\tilde p_3(m)$ and $\text{Disc}(Q_m)$ are obtained by
evaluation  of the three classical covariants
of the fundamental cubic $Q_m$ of $m$  at  the ``trace'' $T_m$ (see Definition \ref{def:KH}) and hence the 
functions $p_2, p_3$
\beqa
\label{eq:p2p3}
p_3(m) = \frac{\tilde p_3(m)}{\text{Disc}(Q_m)} \ , \ \  p_2(m) = \frac{\tilde p_2(m)}{\text{Disc}(Q_m)} \ ,
\eeqa
are well defined on the moduli space ${\cal M}^c_3$ and take their values in $\FF$.
For $m \in \M_3^c$ we introduce the following notations:
\begin{enumerate}
\item $\hat \nu (m) = (p_3(m),p_2(m))$;
\item $\FF_{Q_m}$ is the  splitting field of $Q_m$ (this is well defined up to equivalence of field extensions);
\item $\Gamma_{\text{Cardano}}= \Big\{(a,b)\in \FF \times \FF \ \ \text{s.t.} \ \ 27 a^2 + 4 b^3=0 \Big\}$;
\item ${\cal E}_k$ denotes the set of equivalences classes of extensions of degree at most $k$ which are isomorphic
to splitting
fields of cubic polynomials.
\end{enumerate}
The essential result of this section can now be summarised in the

\begin{theorem}\label{theo:gose}
 Let $m,m'$ be in $\M^c_3$.
\begin{enumerate}
\item[(i)] There exists a non-trivial automorphism of $m$ {\it iff} $\hat \nu(m) \in \Gamma_{\text{Cardano}}$.
\item[(ii)]
Suppose $\hat\nu (m),\hat \nu (m') \not \in\Gamma_{\text{Cardano}}$. Then there exists $g\in GL(V)$ such that $m' = g \cdot m$ iff 
$\hat\nu (m)=\hat \nu (m') $ and $\hat \nu$   establishes  a bijection of  ${\cal M}^c_3 \setminus
\hat \nu^{-1}(\Gamma_{\text{Cardano}})$ with $\FF^2\setminus\Gamma_{\text{Cardano}}$.
\item[(iii)] Suppose $\hat\nu (m),\hat \nu (m')  \in\Gamma_{\text{Cardano}}$. Then exists  $g\in GL(V)$ such that $m' = g \cdot m$ iff
$\hat\nu (m)=\hat \nu (m') $ and $\FF_{Q_m} \cong \FF_{Q_m'}$.
\begin{enumerate}
\item[(iii.1)]The map  $(\hat \nu,\FF_{Q_{\bullet}})$ establishes a bijection between $\hat \nu^{-1}\big(\Gamma_{\text{Cardano}}\setminus\{0\}\big)$
and $\big(\Gamma_{\text{Cardano}}\setminus\{0\}\big) \times {\cal E}_2$;
\item[(iii.2)] The map $\FF_{Q_{\bullet}}$ establishes a bijection between $\hat \nu^{-1}\{0\}$
and ${\cal E}_6$;
\end{enumerate}
\end{enumerate}
\end{theorem}

\begin{remark}
Note that $27 p_3(m)^2 + 4 p_2(m)^3$ is minus the discrimant of the cubic $x^3 + p_2(m) x - p_3(m)$. By Theorem \ref{theo:miracle}
this equals $-2^4\cdot 3^6 \cdot \text{Inv}_m^2/\text{Disc}(Q_m)^3$.
\end{remark}

If $\FF$ is an {\it algebraically closed} field, Anan'in and Mironov \cite{am} constructed an embedding
$\Phi:{\cal M}\to \mathbb P^8$ of ${\cal M}$,  the $GL(V)-$moduli space of {\it stable} two-dimensional algebras, in
the projective space $\mathbb P^8$.
By Remarks  \ref{rem:am-disc}, \ref{rem:am-p3} and \ref{rem:am-p2} the restriction  of $\Phi$ to ${\cal M}^c_3$ is given by
\beqa
\label{eq:id1}
\Phi(m)=[\tilde p_2(m),\tilde p_3(m),0,0,0,\text{Disc}(Q_m),0,0] \  
\eeqa
which clearly defines an embedding (also denoted $\Phi$) of $ {\cal M}^c_3$ in  $\mathbb P^2$: 
\beqa
\label{eq:id2}
\Phi(m)= [\tilde p_2(m),\tilde p_3(m),\text{Disc}(Q_m)]  \ .
\eeqa
It follows from the definition of $\M^c_3$ that
\beqa
\Phi\Big({\cal M}^c_3\Big) = \Big\{[a,b,c] \in \mathbb P^2 \ \ \text{s.t.} \ \ c \ne 0 \Big\} \ . \nn
\eeqa
Hence  we can think of the map  $\hat \nu : {\cal M}^c_3 \to \FF \times \FF$ of Theorem   \ref{theo:gose}  given by
\beqa
\hat \nu(m)=\Bigg(\frac{\tilde p_3(m)}{\text{Disc}(Q_m)}, \frac{\tilde p_2(m)}{\text{Disc}(Q_m)}\Bigg) =(p_3(m),p_2(m))\ , \nn
\eeqa
as providing  two affine coordinates on the moduli space  ${\cal M}^c_3$.



In order to prove  the ``projective'' versions Theorem \ref{theo:gose} we need to consider
 the homogenisation of $\Gamma_{\text{Cardano}}$:
\beqa
\widetilde{\Gamma}_{\text{Cardano}}= \Bigg\{ [a,b,c] \in \mathbb P^2 \ \ \text{s.t.} \ \  27 c a^2+ 4 b^3 = 0 \Bigg\} \ . \nn
\eeqa
We will show in Appendix \ref{sec:non-gen} that if $\Phi(m) \not\in \widetilde{\Gamma}_{\text{Cardano}}$ then $\Phi(m)$ uniquely
  determines the equivalence class  of the {\it stable} commutative algebra $m$.
What is more, 
we will give a complete set of discrete invariants  which distinguish stable commutative algebras
$m,m'$ such that $\Phi(m) =\Phi(m')\in \widetilde{\Gamma}_{\text{Cardano}}$.


Finally to complete this section let us point out that, as we shall see in Sections \ref{sec:mod-sl} and \ref{sec:compact},
the map $\zeta : \M^c_{\text{st}}/SL(V) \to W$ (here $W$ is a suitable four-dimensional vector space) given by
\beqa
\zeta(m)=(\tilde p_3(m), \tilde p_2(m), \text{Disc}(Q_m), \text{Inv}_m) \ , \nn
\eeqa
defines an embedding (at least if $\FF$ is algebraically closed)
of the  $SL(V)-$moduli space of stable, commutative, two-dimensional algebras
onto  the Eisenstein hypersurface
\beqa
\Gamma_{\text{Eisenstein}} = \Bigg\{(a,b,c,d) \in W \ \ \text{s.t.}\ \  27 ca^2 + 4 b^3 +2^4 3^6 d^2 =0 \Bigg\}  \ .\nn
\eeqa
Since $( \tilde p_3(m),  \tilde p_2(m), \text{Disc}(Q_m)) \ne (0,0,0)$, one can define  the map
$\pi: \Gamma_{\text{Eisenstein}} \to \mathbb P^2$  given by
\beqa
\pi\Big(\tilde p_3(m),  \tilde p_2(m),\text{Disc}(Q_m), \text{Inv}_m\Big)=[\tilde p_3(m),  \tilde p_2(m),\text{Disc}(Q_m)] \ ,\nn
\eeqa
 (some care must be taken to make sense of this) and then
 \beqa
\pi^{-1} \Big(\widetilde{\Gamma}_{\text{Cardano}}\Big) = \Bigg\{ (a,b,c,d) \in \Gamma_{\text{Eisenstein}} \ \ \text{s.t.} \ \
d=0 \Bigg\} \ . \nn
\eeqa

\subsection{Fundamental triples and the invariants $p_2,p_3$}\label{sec:triple}

The main tool in the proof of Theorem \ref{theo:gose} will be the fact that
 one can associate to  $m$ in $\M^c_3$ 
a unique  triple  of vectors in $V_{\FF_{Q_m}} $ ({\it i.e.} $V$  tensored by
the splitting  field $\FF_{Q_m}$  of $Q_m$) together with a unique triple of elements of $\FF_{Q_m}$.
Very surprisingly, it turns out that the two polynomial invariants $p_2$ and $p_3$ in \eqref{eq:p2p3} are essentially the
elementary symmetric functions of these  three numbers.
Recall that since $\text{Disc}(Q_m) \ne 0$,  the polynomial $Q_m$ has three distinct roots in its splitting field.

\begin{remark}
\label{rem:splitting}
Since $Q_m$ is a cubic, the degree of a splitting field is either one, two, three or six.
\end{remark}

\begin{remark}
\label{rem:gal-perm}
If $\FF_{Q_m}$ is a splitting field of $Q_m$ we can write $Q_m =F_1 F_2 F_3 \otimes b$ with $b \in \Lambda^2(V_{\FF_{Q_m}})$
and $F_i \in (V_{\FF_{Q_m}})^*$ such that $(F_1,F_2,F_3)$
are pairwise independent. The Galois group $\mathrm{Gal}(\FF_{Q_m}/\FF)$  permutes
the elements of the set 
${\cal R}= \big\{\mathrm{Ker}(F_1),\mathrm{Ker}(F_2),\mathrm{Ker}(F_3)\big\}$
and this gives  a homomorphism to  the permutation group $S_3({\cal R})$.
Note that although the factorisation of $Q_m$ is not unique, the set ${\cal R}$ is
uniquely defined by $Q_m$.
\end{remark}

\begin{theorem} 
\label{theo:can}
Let $m \in \M^c_3$, let $Q_m$ be its fundamental cubic and let $V_{\mathbb \FF_{Q_m}} = V \otimes_\mathbb F \FF_{Q_m}$, where $\FF_{Q_m}$ is
the splitting field of $Q_m$.
\begin{enumerate}
\item[(i)] There exists a triple $((f_1,\gamma_1),(f_2,\gamma_2),(f_3,\gamma_3)) \in
(V_{\FF_{Q_m}} \times \mathbb \FF_{Q_m})^3$ such that
 \beqa
\begin{array}{ll}
a)& (f_1,f_2,f_3) {\rm\  are \ pairwise \ linearly \ independent}; \\
b)& f_1+f_2+f_3=0 \ ,\\
c)& f_i^2= \gamma_i f_i \ ,\\
d)&\gamma_1+\gamma_2 +\gamma_3=-1 \ .
\end{array} \nn
\eeqa
\item[(ii)] The triple of a) is unique up to permutation.
\end{enumerate} 
\end{theorem}

\begin{demo}
By hypothesis there exist three non-zero linear forms $F_1,F_2,F_3 \in V_{\FF_{Q_m}}^*$ which are
pairwise independent and  $b \in \Lambda^2(V_{\FF_{Q_m}})$ such 
$Q_m(v)= F_1(v) F_2(v) F_3(v)  b$ for all $v$ in $V$.
The vector space 
\beqa
B= \Big\{(f_1,f_2,f_3) \in {\rm Ker}(F_1) \times {\rm Ker}(F_2) \times {\rm Ker}(F_3) :
f_1 + f_2 + f_3 =0 \Big\} \ ,\nn
\eeqa  
is one-dimensional.  Since $Q_m(v)= m(v,v)\wedge v$ and $Q_m(f_i)=0$,
 there
exist three linear forms $\gamma_i : \text{Ker}(F_i) \to \FF_{Q_m}$ such that
\beqa
m(f_i,f_i)= \gamma_i(f_i) f_i \ . \nn
\eeqa
We define $\gamma : B \to \FF_{Q_m}$  by
\beqa
\gamma(f_1,f_2,f_3) = \gamma_1(f_1) + \gamma_2(f_2) + \gamma_3(f_3) \ . \nn
\eeqa
\begin{lemma}
\label{lem:gamma}
The map $\gamma$ is non-zero.
\end{lemma}
\begin{demo}
Suppose for contradiction that $\gamma\equiv 0$. Then
\beqa
m(f_1,f_2) &=& \frac12\big(m(f_1+f_2,f_1+f_2) - m(f_1,f_1) - m(f_2,f_2)\big)\nn \\ 
&=&\frac12\big(m(f_3,f_3) - m(f_1,f_1) - m(f_2,f_2)\big)=\frac12
\big(\gamma_3(f_3) f_3 - \gamma_1(f_1) f_1-\gamma_2(f_2) f_2\big)\nn \\ 
&=&\frac12 \gamma_2(f_2) f_1 + \frac12 \gamma_1(f_1) f_2 \ . \nn
\eeqa
But now for any vector $v=x f_1 + y f_2$ we have
\beqa
m(v,v)&=& x^2 \gamma_1(f_1) f_1 + y^2 \gamma_2(f_2) +
x y (\gamma_2(f_2) f_1 + \gamma_1(f_1) f_2)=(x \gamma_1(f_1) + y \gamma_2(f_2) v \nn \ .
\eeqa
From this $Q_m(v)=0$ for all $v \in V_{\FF_{Q_m}}$ which is a contradiction. 
\end{demo}

By Lemma  \ref{lem:gamma} there exists a unique $(f_1,f_2,f_3) \in B$ such that
$\gamma(f_1,f_2,f_3)=-1$ and this completes the proof of the first part of the
theorem.

To prove  unicity assume now that one can find another triple
 $((f_1',\gamma_1'),(f_2',\gamma_2'),(f_3',\gamma_3')) \in
(V_{\FF_{Q_m}} \times \mathbb \FF_{Q_m})^3$ satisfying
\beqa
\begin{array}{ll}
i)& (f_1',f_2',f_3') {\rm\  pairwise \ linearly \ independent},\\
ii)& f_1'+f_2'+f_3'=0 \ ,\\
iii)& f_i'{}^2= \gamma_i' f_i' \ ,\\
iv)&\gamma_1'+\gamma_2' +\gamma_3'=-1 \ .
\end{array} \nn
\eeqa
Since the $f_i$ are idempotent and  pairwise linearly  independent, 
renumbering if necessary, there exist non-zero
$\lambda_1, \lambda_2 , \lambda_3 \in \mathbb K$ such that
\beqa
f_i'=\lambda_i f_i \ . \nn
\eeqa
Substituting in ii) gives
\beqa
\lambda_1 f_1 + \lambda_2 f_2 + \lambda_3 f_3 =0 \ , \nn 
\eeqa
and since $f_3 = - f_1 - f_2$ and $(f_1,f_2)$ is a basis of $V_{\FF_{Q_m}}$, we get
$\lambda_1=\lambda_2=\lambda_3=\lambda$. 

Substituting $f_i'{}= \lambda f_i$ in iii) gives $\gamma_i'=\lambda \gamma_i$ and 
finally, substituting in iv), gives $\lambda=1$. 
\end{demo}

\begin{remark}
\label{rem:galois-compat}
The Galois group of the extension $\FF_{Q_m}/\FF$ acts on $V_{\FF_{Q_m}}$ by automorphism and 
hence by  Theorem \ref{theo:can} (ii), each element  induces a permutation of the
triple  $((f_1,\gamma_1),(f_2,\gamma_2),(f_3,\gamma_3))$.
\end{remark}

As an application of Theorem \ref{theo:can} we now characterise the  possible automorphism groups
of $m \in \M^c_3$. For this we need the notion ``idemvalue'' which is well defined by  Theorem \ref{theo:can}

\begin{definition}
The fundamental triple of $m \in\M^c_c$ is the triple $((f_1,\gamma_1),(f_2,\gamma_2),(f_3,\gamma_3))$
associated to $m$ by Theorem \ref{theo:can}. The idemvalues of $m$ are $\gamma_1,\gamma_2,\gamma_3$ and
the reduced idemvalues of $m$ are $\Delta_i = 3 \gamma_i +1$. Note that 
\beqa
\Delta_1 + \Delta_2 + \Delta_3 =  0 \ .\nn
\eeqa
\end{definition}

\begin{theorem}
\label{theo:Delp}
Let $m \in \M^c_3$ and let $\Delta_1,\Delta_2,\Delta_3$ be the reduced idemvalues. Then the polynomial invariants
$p_2(m), p_3(m)$ are given by
\beqa
p_2(m) &=& \Delta_1 \Delta_2 + \Delta_2 \Delta_3 + \Delta_3 \Delta_1 \ , \nn\\
p_3(m)& =& \Delta_1 \Delta_2 \Delta_3 \ . \nn 
\eeqa
\end{theorem}
\begin{demo}
In the basis $(f_1,f_2)$ the multiplication table is given by
\beqa
\label{eq:f}
m(f_1,f_1) = \frac{\Delta_1-1}3 f_1\ , \ \ 
m(f_2,f_2) = \frac{\Delta_2-1}3 f_2\ , \ \ 
m(f_1,f_2) = \frac{2+\Delta_2}6 f_1 + \frac{2+\Delta_1}6 f_2 \ .
\eeqa

A straightforward but long computation gives (see  \eqref{eq:decom}, \eqref{eq:dD}  and \eqref{eq:Invd})
\beqa
\label{eq:Inv}
\text{Disc}(Q_m)(v_1,v_2,v_3,v_4)&=& v_{12} v_{34} \ , \nn \\
\text{Disc}(D_m)(v_1,v_2,v_3,v_4)&=&\frac1{27} (\Delta_1-1)(\Delta_2-1)(\Delta_3-1)v_{12} v_{34}\ , \nn \\
\text{Inv}_m(v_1,v_2,v_3,v_4,v_5,v_6)&=& -\frac1{108} (\Delta_1- \Delta_2)  (\Delta_2 -\Delta_3) (\Delta_3 -\Delta_1)
v_{12} v_{34} v_{56}  \ ,
\eeqa 
where $v_{ij}= x_i y_j - x_j y_i$ and hence (see  \eqref{eq:pt3} and \eqref{eq:pt2})
\beqa
p_2(m) &=& \Delta_1 \Delta_2 + \Delta_2 \Delta_3+ \Delta_3 \Delta_1 \ , \nn \\
p_3(m) &=& \Delta_1 \Delta_2  \Delta_3 \ . \nn 
\eeqa
\end{demo}

\begin{corollary}\label{cor:altern}
(Alternative proof of Theorem \ref{theo:miracle} for a generic algebra.)
Let $m$ be in $\M^c_3$.  Then as elements of $L^6$ we have the identity 
\beqa
27 \mathrm{Disc}(Q_m) \tilde p_3(m)^2+ 4 \tilde p_2(m)^3 = -2^4 \cdot 3^6\  \mathrm{Inv}_m^2 \ . \nn
\eeqa

\begin{demo}
Since $\mathrm{Disc}(Q_m) \ne 0$, the identity above in
$L^6$
is equivalent to the identity in $\FF$
\beqa
27 p_3(m)^2 + 4 p_2(m)^3 = -2^4\cdot 3^6 \frac{\mathrm{Inv}_m^2}{\mathrm{Disc}(Q_m)^3}\ , \nn
\eeqa
which, using \eqref{eq:Inv},  reduces to
\beqa
27 p_3(m)^2 + 4 p_2(m)^3 =
-\Big((\Delta_1 -\Delta_2) (\Delta_2 -\Delta_3) (\Delta_3 -\Delta_1)\Big)^2\ . \nn
\eeqa
This is the classical formula \ for the discriminant of the cubic polynomial $x^3+p_2(m) x -p_3(m)$ in terms of its roots.
\end{demo}
\end{corollary}
\begin{remark}
It is important to emphasise that the bivectors $f_1\wedge f_2, f_2\wedge f_3,f_3\wedge f_1$ although {\it a priori}
elements of $\Lambda^2(V_{\FF_{Q_m}})$ are in fact elements of  $L^{-1}$. 
This means that in order to deduce the expression for the invariants above in a basis $(e_1,e_2)$ of $V$ from their expressions in the 
basis $(f_1,f_2)$ of $V_{\FF_{Q_m}}$, one has to express $f_1\wedge f_2$ in terms of $e_1\wedge e_2$.
When the $\Delta_i$'s are distinct this is the formula \eqref{eq:e-f} below.
\end{remark}

\noi
A simple characterisation of the automorphism group of a generic algebra can be given in terms of the  invariants Inv$_m$ and 
$p_2(m),p_3(m)$.  (For alternative calculations  of the automorphism group of  real division algebras see \cite{dz,d},
Proposition 4/Proposition 2.8.) 
  
\begin{corollary}
\label{cor:aut}
\phantom{marcus}\hskip 2.truecm
\begin{enumerate}
\item
Let $m \in \M^c_3$, let $\FF_{Q_m}$ be a splitting field of $Q_m$ and let  
$\widetilde{\A}(m)= \{ g \in GL(V\otimes \FF_{Q_m})$ such that $g\cdot m = m\}$.
Then, 
\beqa
\widetilde {\A}(m) \cong \left\{
   \begin{array}{cl}
    \{\bf 1\}&\text{\rm{if }}  {\mathrm{Inv}}_m\ne 0 \  ,\\
    \mathbb Z_2&\text{\rm{if }}  {\mathrm{Inv}}_m=0 \ \ \text{and} \ \ (p_2(m), p_3(m)) \neq (0,0)  \ ,   \\
    S_3&\text{\rm{if }}  p_2(m) = p_3(m) =0 \ , \\
   \end{array}\right.\nn
\eeqa
where $S_3$ is the group  of permutations  of three objects. Recall that Inv$_m=0$ {\it iff} $27 p_3(m)^2 + 4 p_2(m)^3=0$.
\item Let $m \in \M^c_3$,  and let  $\A(m)= \{ g \in GL(V)$ such that $g\cdot m = m$\}
\beqa
 {\A}(m) \cong \left\{
   \begin{array}{cl}
    \{\bf 1\}&\text{\rm{if }}  {\mathrm{Inv}}_m\ne 0 \  ,\\
    \mathbb Z_2&\text{\rm{if }}  {\mathrm{Inv}}_m=0 \ \ \text{and} \ \ (p_2(m), p_3(m)) \neq (0,0)  \ ,   \\
    S_3&\text{\rm{if }}  p_2(m) = p_3(m) =0  \ \text{and} \ [\FF_{Q_m}:\FF]=1 , \\
    \ZZ_2&\text{\rm{if }}  p_2(m) = p_3(m) =0  \ \text{and} \ [\FF_{Q_m}:\FF]=2 , \\
      \ZZ_3&\text{\rm{if }}  p_2(m) = p_3(m) =0  \ \text{and} \ [\FF_{Q_m}:\FF]=3 , \\
      \{1\}&\text{\rm{if }}  p_2(m) = p_3(m) =0  \ \text{and} \ [\FF_{Q_m}:\FF]=6 \ .
   \end{array}\right.\nn
\eeqa
\end{enumerate}
\end{corollary}
\begin{demo}
(1): 
From the formul\ae \  \ref{eq:Inv} and $\Delta_1 + \Delta_2 + \Delta_3=0$  it follows that 
\beqa
\begin{array}{ll}
    \hskip .3truecm {\mathrm {Inv}}_m \neq 0 &\text{\rm{if }}  m \text{\rm{ has  three  distinct  reduced idemvalues}},\\ \\
   \left\{\begin{array}{l}
   \hskip -.2truecm  {\mathrm{Inv}}_m=0 \\ \hskip -.2truecm  (p_2(m), p_3(m)) \neq (0,0) 
\end{array} \right.
 &\text{\rm{if }}  m \text{\rm{ has  two  distinct reduced   idemvalues}}, \\ \\
  \hskip .4truecm  p_2(m) = p_3(m) =0 &\text{\rm{if }}  m \text{\rm{ has  one  reduced idemvalue}}. \\
   \end{array}\nn
\eeqa
The result follows from this and the fact that by Theorem \ref{theo:can}
any  automorphism is equivalent to  a permutation of the
triple  $((f_1,\Delta_1),(f_2,\Delta_2),(f_3,\Delta_3))$.

(2):  This is a consequence of  Propositions 2.1 and 2.3 in  \cite{wright}.
\end{demo}

\newpage
\begin{corollary}
\label{cor:detminus}
Let   $m \in \M_3^c$ be such that Inv$_m=0$ (or, equivalently, $m \in \Gamma_{\mathrm{Cardano}}$ in the terminology of  Eq.[\ref{eq:card}]). 
Then there exists $g$ in $\widetilde{\text{Aut}}(m)$ such 
\begin{enumerate}
\item $g \cdot m = m $;
\item $\det(g) =-1$.
\end{enumerate}
\end{corollary}

\begin{demo}
First note that $-$Id$_{V\otimes \FF_{Q_m}}$ cannot be an automorphism of  $m$ in $\M_3^c$. If so, we would have
\beqa
m(-v,-w) = - m(v,w) \ , \forall v,w \in V \ , \nn
\eeqa
which is impossible. Secondly, by Corollary \ref{cor:aut}{\color{blue}(1)},
 the group of automorphisms of $m$ is isomorphic to either $\mathbb Z_2$ or $S_3$. Hence $m$ has a non-trivial
automorphism $g$ of square one which 
 cannot be $-$Id$_{V\otimes \FF_{Q_m}}$ by the above remark. Thus $g\not =  $Id$_{V\otimes \FF_{Q_m}}$, $g\not =-$Id$_{V\otimes \FF_{Q_m}}$  and $\det(g) =-1$.
\end{demo}
\begin{remark}
In \cite{am} the authors identify explicitly the image
of  the singular locus  Sing$({\cal M})$ under the embedding 
$\Phi: {\cal M} \to\mathbb P^8$.
Furthermore, they  show that the image of Sing$($Sing$({\cal M}))$  is a single  point which they also identify.
Restricting to generic commutative algebras their results imply
\beqa
\hat \nu\Big(\text{Sing}({\cal M}) \cap {\cal M}^c_3\Big) &=& {\Gamma}_{\text{Cardano}} \ , \nn\\
\hat \nu\Big(\text{Sing}(\text{Sing}({\cal M})) \cap {\cal M}^c_3\Big)&=&(0,0) \in {\Gamma}_{\text{Cardano}} \ , \nn
\eeqa
and restricting to stable commutative algebras 
\beqa
\Phi\Big(\text{Sing}({\cal M})\cap {\cal M}^c_{\text{st}}\Big) &=& \widetilde{\Gamma}_{\text{Cardano}} \ , \nn\\
\Phi\Big(\text{Sing}(\text{Sing}({\cal M}\cap {\cal M}^c_{\text{st}})\Big)&=&[0,0,1]  \in \widetilde{\Gamma}_{\text{Cardano}} \ .\nn
\eeqa

The significance of the Corollary above is that it gives 
 a ``geometric'' interpretation of the singular locus and of the singular-singular locus.
Namely, if $m$ is a generic commutative two-dimensional algebra then:\\
\beqa
\begin{array}{lllll}
\big[m\big]&\in & \text{Sing}({\cal M})& {\it iff} & \text{Aut}(m) \ne \{\text{Id}\}\ ; \\ \\
\big[m\big] &\in & \text{Sing}({\cal M}) \setminus  \text{Sing}(\text{Sing}({\cal M}))&{\it iff} &\widetilde{\text{Aut}}(m) =
\mathbb Z_2 \ ; \\ \\
\big[m\big]&\in&  \text{Sing}(\text{Sing}({\cal M})) & {\it iff} & \widetilde{\text{Aut}}(m) = S_3 \ ;
\end{array}
\nn
\eeqa
 We will give the proof  of  the ``stable'' version of this result  in the Appendix.
\end{remark}

\subsection{Moduli space}

In this section we give a  simple parametrisation of the  $GL(V)-$moduli space of
generic  commutative algebras.  
When the field $\mathbb F$                                                                       
is algebraically closed the parameter space turns out to be $\mathbb F^2$. 
If $\FF$ is not algebraically closed,  a point in the plane $\FF^2$ does not necessarily 
characterise a class of algebras, and another (discrete) invariant may be needed
to distinguish non-isomorphic  algebras which become isomorphic 
after tensoring with an  extension of $\FF$.

In order to ease the notation, from now on we write $(p_3,p_2)$ instead of $(p_3(m),p_2(m))$ whenever 
 the choice of $m$ is clear from the context.

Our strategy will be to consider the map $\nu: \M^c_3 \to \FF^2$ defined by 
\beqa
\nu(m) = (p_3,p_2) \ . \nn
\eeqa  
This factors to a map $\hat \nu: {\cal M}^c_3\to \FF^2$  ({\it c.f. \eqref{eq:p2p3}}). 
Note that in terms of the reduced idemvalues of $m$, we have
$\nu(m)=(\Delta_1 \Delta_2 \Delta_3,\Delta_1 \Delta_2 + \Delta_2 \Delta_3+\Delta_3 \Delta_1)$ and so the degenerate elliptic curve
\beqa
\label{eq:card}
\Gamma_{\mathrm{Cardano}}=\Big\{ (p_3,p_2) \in \FF^2: \ -27 p_3^2 -4 p_2^3=0
\Big\} \ ,
\eeqa
defined by the vanishing of the discriminant of the cubic  
$Q_{\text{mod}}(x)=x^3 + p_2 x - p_3$ 
will play a special r\^ole. This is because the roots of $Q_{\text{mod}}$ are $\Delta_1,\Delta_2,\Delta_3$ and we will
call $Q_{\text{mod}}$ the modular cubic of the algebra $\M^c$. Note that by Corollary \ref{cor:aut}: 
$\nu(m) \in \Gamma_{\mathrm{Cardano}}$ {\it iff} the algebra $\M^c$ has a non-trivial automorphism.

We first show that over $\FF^2 \setminus \Gamma_{\mathrm{Cardano}}$ the map 
 $\hat\nu$ is a bijection. Then we will show that over $\Gamma_{\mathrm{Cardano}}\setminus\{(0,0)\}$ there is a point in the moduli
space for every  class of extensions of $\FF$ of degree at most two and   and finally that over  $(0,0)$  
there is a point in the moduli
space for  every class of Galois extensions of $\FF$ of degree one, two, three and six.

\subsubsection{Surjectivity of $\nu$}

\begin{theorem}
\label{theo:inv-gen}
\phantom{toto}\hskip 12truecm
\begin{itemize}
\item[(i)] $\nu(g\cdot m) = \nu(m), \forall m \in \M_3^c, \  \forall g \in \G$.
\item[(ii)] $\nu: \M_3^c\to \mathbb F^2$ is surjective.
\end{itemize}
\end{theorem}

\begin{demo}
Part (i) is a consequence of \eqref{eq:p2p3}. 

To prove (ii), let $(p_3,p_2) \in \FF^2$. First we suppose that
 $(p_3,p_2) \not \in \Gamma_{\text{Cardano}}$, {\it i.e.},
$-27 p_3^2 - 4 p_2^3\ne 0$. Since $-27 p_3^2 -4 p_2^3$ is the discriminant of
the   cubic $P(x)=x^3 +p_2 x -p_3$, this means $P$  
has three distinct roots $\Delta_1,\Delta_2$ and $\Delta_3$ 
in a splitting field $\KK_1$. 
Since the characteristic of $\mathbb F$ is zero or greater than three and $P$ is of degree three, $\mathbb K_1/ \mathbb F$ is
a Galois extension (see {\it e.g.} Artin \cite{ar}).
\begin{lemma} \label{lem:base}
Let $V_{\KK_1} = V \otimes_\FF \KK_1$. Then there exists a triplet $f_1,f_2,f_3 \in V_{\KK_1}$ such that
\begin{itemize}
\item[(a)] $(f_1,f_2,f_3)$ are pairwise linearly independent;
\item[(b)] $f_1 + f_2 + f_3 =0$;
\item[(c)] $\Big\{f_1,f_2,f_3\Big\}$ is stable under the action of the Galois group
${\mathrm{ Gal}}(\KK_1/\FF)$;
\item[(d)] There exists a group homomorphism  $\gamma : {\mathrm{ Gal}}(\KK_1/\FF) \to S_3$ 
 such that for all $g \in {\mathrm{ Gal}}(\KK_1/\FF)$ and for $i=1,2,3$,
\beqa
g\cdot (f_i,\Delta_i) = (f_{\gamma(g) \cdot i}, \Delta_{\gamma(g) \cdot i}) \ . \nn
\eeqa
\end{itemize}
\end{lemma}
\begin{demo}
Let $e_1, e_2$ be a basis of $V$ and define the following elements of $V_{\KK_1}$:
\beqa
\label{eq:prod-gene}
f_1&=& \Delta_1 e_1 + (\Delta_2 \Delta_3 -\frac{p_2}3)e_2 \ , \nn\\
f_2&=& \Delta_2 e_1 + (\Delta_3 \Delta_1 -\frac{p_2}3)e_2 \ , \\
f_3&=& \Delta_3 e_1 + (\Delta_1 \Delta_2 -\frac{p_2}3) e_2\ . \nn
\eeqa
If $\KK_1=\FF$, (c) and (d) are evident. If $\mathbb K_1 \ne \mathbb F$ 
the Galois group acts  by permutations on the indices of the roots and so (c) and
(d) are obviously satisfied. Property (b) follows since $\Delta_1, \Delta_2$ and
$\Delta_3$ are the roots of the polynomial $P(x) = x^3 + p_2 x - p_3$.

We now prove (a). This follows from
\beqa
\label{eq:e-f}
f_i\wedge f_j =(\Delta_j -\Delta_i) \Big( (\Delta_i + \Delta_j)^2 + \frac{p_2}3\Big) 
e_1\wedge e_2\ = \frac13(\Delta_1 -\Delta_2)(\Delta_2-\Delta_3)(\Delta_3-\Delta_1)
e_1 \wedge e_2  
\eeqa
 and the fact that the  $\Delta_i$ are
distinct by hypothesis. Hence, we have proved (a) for $(f_1,f_2,$ $f_3)$.
\end{demo}

Since
$f_1,f_2 \in V_{\KK_1}$ is a basis we can  define a multiplication $m$ on $V_{\KK_1}$ by,
\beqa
\label{eq:prod}
m(f_1,f_1)= \frac{\Delta_1 -1}{3} f_1\ , \ \
m(f_2,f_2)=\frac{\Delta_2-1}{3} f_2 , \ \ 
m(f_1,f_2)= \frac{2+\Delta_2}6 f_1 + \frac{2+\Delta_1}6 f_2 \ .  
\eeqa  
Then a direct calculation shows that in fact we have for all $i,j=1,2,3$ such that $i \neq j$,
\beqa
\label{eq:ff}
m(f_i,f_i)=\frac{\Delta_i-1}{3} f_i , \ \
m(f_i,f_j)= \frac{2+\Delta_j}6 f_i+ \frac{2+\Delta_i}6 f_j\ .   
    \eeqa
Since for any $g \in \text{Gal}(\KK_1/\FF)$ we have
\beqa
g\cdot f_i= f_{\gamma(g) \cdot i }\ , \ \ 
g\cdot \Delta_i= \Delta_{\gamma(g) \cdot i} \ , \nn
\eeqa
it follows from \eqref{eq:ff}  that 
\beqa
g\cdot m(f_i,f_j)= m(g \cdot f_i,g \cdot f_j)  \nn
\eeqa 
and so 
$$g\cdot m(v,w)= m(g \cdot v,g \cdot w)\  , \forall v,w \in V_{\mathbb K_1} \ . $$
In particular, since $V \subseteq V_{\KK_1}$ is the fixed  point set of the Galois group and since $\KK_1/\FF$ is Galois,
this implies that
$m(v,w) \in V$ if $v,w \in V$.

It remains  to prove that the fundamental cubic of $m$ has three distinct roots.
This follows since
\beqa
m(f_i,f_i)= \frac{\Delta_i -1}{3} f_i\ , \nn 
\eeqa
so by  \eqref{eq:idem}
  there are three pairwise independent   idempotents  
and the fundamental cubic (not to be confused with $P(x) !$) indeed has three distinct roots.

To recap we have now shown that if   $(p_3,p_2) \not \in \Gamma_{\text{Cardano}}$ then
there exists $m \in \M_3^c$ such that $\nu(m)=(p_3,p_2)$. To complete the proof of (ii) we now
consider the case $(p_3,p_2)  \in \Gamma_{\text{Cardano}}$. This means  the cubic  $P(x)$
has a multiple root in $\FF$ and we have
\beqa
P(x) = x^3 + p_2 x - p_3 = \left\{ \begin{array}{l}
(x-\frac{3 p_3}{2 p_2})^2 (x + \frac{3p_3}{p_2}) \ \  \text{if~} p_2 \ne 0 \\
x^3 \ \  \text{if~} p_2 = 0 \ . \nn
\end{array} \right. \nn
\eeqa
Let $(\Delta,\Delta,-2 \Delta)$ be the roots of $P(x)$. Choose $(e_1,e_2)$ a basis of $V$ and
set,
\beqa
f_1 = e_1 -2 e_2 \ , \ \ 
f_2= -2 e_1 + e_2 \ , 
f_3=e_1 + e_2 \ . \nn
\eeqa
It is clear  that $(f_1,f_2,f_3)$ are pairwise linearly independent and that $f_1 + f_2 + f_3=0$.
Define the multiplication $m: V \times V \to V$  by
\beqa
m(f_1,f_1) = \frac{\Delta -1}3 f_1\ , \ \ 
m(f_2,f_2) = \frac{\Delta -1}3 f_2\ \ ,
m(f_1,f_2) = \frac{2 + \Delta}6 \big(f_1 + f_2) \ .\nn
\eeqa
One checks that
\beqa
m(f_3,f_3) = \frac{-2 \Delta -1}3 f_3 \ , \ \ 
m(f_i,f_3) =  \frac{2-2  \Delta }6 f_i + \frac{2 + \Delta}6 f_3 \ , i=1,2 \ . \nn
\eeqa
Thus $m$ has three idempotents $f_1,f_2,f_3$ and by definition,
\beqa
\nu(m)=( -2\Delta^3, -3 \Delta^2 )= (p_3,p_2) \ . \nn 
\eeqa
This complete the proof of (ii).
\end{demo}
\begin{remark}
In the proof we showed that the product \eqref{eq:prod} defined {\it a priori}
on $V_{\KK_1}$  in fact defines  a product on $V$ because it is compatible
with the action of the Galois group $\text{Gal}(\KK_1/\FF)$. In order to get explicit
formul\ae \ for the products of the $e_i$,  inverting the system \eqref{eq:prod-gene} gives
\beqa
e_1 = \frac1{d}
\Big((3\Delta_3 \Delta_1 -p_2) f_1 + (p_2 -3\Delta_2 \Delta_3)f_2\Big) \ ,
e_2=-\frac3{d}
\Big(\Delta_2 f_1 - \Delta_1 f_2\Big) \ ,\nn
\eeqa
where $d=(\Delta_1-\Delta_2)(\Delta_2 -\Delta_3)(\Delta_3-\Delta_1)$ and then
(tedious) calculation shows that
\beqa
\label{eq:mul-gen}
m(e_1,e_1)&=& \frac{1}{3d^2}\Big[\big(-4p_2^3-27p_3^2-9p_2p_3\big) e_1 +
\big(2p_2^3+27 p_3^2\big)e_2\Big] \ ,\nn \\
m(e_2,e_2)&=& \frac{1}{d^2}\Big[9 p_3 e_1 +2 p_2^2 e_2 \Big] \ ,  \\
m(e_1,e_2)&=&\frac3{d^2}\Big[-\frac23 p_2^2 e_1+
\big(-\frac29 p_2^3 -\frac{3}2 p_3^2 +  p_2 p_3\big)e_2 \Big] \ . \nn 
\eeqa
It should be noted that
\beqa
d^2= -27p_3^2 -4 p_2^3 \ , \nn
\eeqa
so  that all of the structure constants are rational functions in the moduli $(p_3,p_2)$.
Note further that these formul\ae \ are valid in characteristic not two or three but
not over $\mathbb Z$.
\end{remark}
\subsubsection{Moduli over  $\FF^2\setminus \Gamma_{\mathrm{Cardano}}$}
We saw in Theorem  \ref{theo:inv-gen} that $\nu(m) \in \FF^2$ is an invariant of the
equivalence classe of $m$. We now show that if $\nu(m)$ is not on a $\Gamma_{\text{Cardano}}$
it completely determines the class of $m$.
\begin{lemma}\label{lem:Q-P}
Let $m  \in  \M^c_3$ with reduced idemvalues $(\Delta_1,\Delta_3,\Delta_3)$
and let $\KK$ be a splitting field of $Q_{\mathrm{mod}}(x)=(x-\Delta_1)(x-\Delta_2)(x-\Delta_3)$.
If  the idemvalues are distinct,  then $\KK$ is also a splitting field of the
fundamental binary cubic $Q_m$ ({\it cf} \ref{def:cubic}).
\end{lemma}
\begin{demo}
Let $\mathbb F_{Q_m}$ be a splitting field of $Q_m$.
By definition, it is always true that the polynomial $Q_{\text{mod}}$ splits over any splitting
field of $Q_m$ so   without loss of generality we can suppose $\KK \subset \FF_{Q_m}$.
We now show that  if the idemvalues are distinct the converse inclusion is also
true.

Let $G=\text{Gal}(\FF_{Q_m}/\KK)$ be the Galois
group of the extension $\FF_{Q_m}/\KK$. By definition, there  exist three idempotents 
$f_1,f_2,f_3 \in V \otimes_{\FF} \FF_{Q_m}$ and three idemvalues $\Delta_1,\Delta_2,\Delta_3 \in \KK$
such that
\beqa
m(f_i,f_i)=\frac{\Delta_i -1}3 f_i \ . \nn
\eeqa
The Galois group $G$ acts by automorphism on $V \otimes_{\FF} \FF_{Q_m}$ and maps
idempotents to idempotents but fixes the  idemvalues $(\Delta_1,\Delta_2,\Delta_3)$.
Since the idemvalues are distinct, the Galois group must fix each 
idempotent because it act by automorphisms.
 However, the extension $\FF_{Q_m}/\KK$ is Galois 
(being a splitting  field of the separable polynomial $Q_m$), so the fundamental theorem 
of Galois theory
implies that the fixed point set of $G$ in $V \otimes_\FF \FF_{Q_m}$
is exactly $V\otimes_\FF \KK$. Hence $f_1,f_2,f_3 \in V \otimes \KK$ and by 
\eqref{eq:idem}, $Q_m(f_1)=Q_m(f_2)=Q_m(f_3)=0$  and $Q_m$ splits 
over $\KK$.
\end{demo}

Perhaps surprisingly, the next theorem shows that just as in the case when $\FF$ is algebraically closed \cite{am}, the invariant $\nu(m)$ completely determines the $GL(V)-$equivalence classe of $m$ even 
when $\FF$ is not algebraically closed.

\begin{theorem}
\label{theo:inv-gen1-1}
 Let $\Gamma_{\mathrm{Cardano}} = \Big\{ (p_3,p_2) \in \mathbb F^2: 27p_3^2 +4p_2^3=0\Big\}$,
and let $\M_3^c{}_g= \{ m \in \M_3^c : \nu(m) \not \in \Gamma_{\text{Cardano}}\}$.
If $m,m'\in \M_3^c{}_g$ satisfy $\nu(m) = \nu(m')$ then there exists $g \in GL(V)$ such that $m'=g \cdot m$.

\end{theorem}
\begin{demo}
Suppose $\nu(m)=\nu(m') \in \FF^2 \setminus \Gamma_{\text{Cardano}} 
$ with $m,m'\in {\cal A}_3^c$. Let $(\Delta_1,\Delta_2,\Delta_3)$ and
$(\Delta'_1,\Delta'_2,\Delta'_3)$ be the idemvalues of  $m$ and $m'$ respectively. 
Since  $\nu(m)=\nu(m')$ we have  
\beqa
Q_{\text{mod}}(x)=(x-\Delta_1)(x-\Delta_2)(x-\Delta_3)=(x-\Delta'_1)(x-\Delta'_2)(x-\Delta'_3) \ , \nn
\eeqa
and therefore  $Q_m$ and $Q_{m'}$ split over $\FF_{Q_\text{mod}}$ by lemma \ref{lem:Q-P}.
Without loss of generality we can suppose $\Delta_i=\Delta_i'$.
By Theorem \ref{theo:can} one can find fundamental triples 
$(f_1,f_2,f_3)$ and
$(f'_1,f'_2,f'_3)$
in $V \otimes_\FF \FF_{Q_{\text{mod}}}$,
 such that
\beqa
m(f_i,f_i)=\frac{\Delta_i-1}3 f_i\ , \ \ 
m'(f'_i,f_i')=\frac{\Delta_i-1}3 f'_i \ .\nn
\eeqa
The map $h: V_{\FF_{Q_\text{mod}}} \to V_{\FF_{Q_\text{mod}}} $ given by $h(f_i)= f_i'$ clearly defines a $\FF_{Q_{\text{mod}}}-$algebra isomorphism.
By Remark \ref{rem:galois-compat},  if $g \in \text{Gal}(\FF_{Q_\text{mod}}/\FF)$ and
$g \cdot \Delta_i = \Delta_j$ then $g \cdot f_i=f_j$ and  $g \cdot f_i'=f_j'$.
Hence  $h\circ g  = g \circ h$ and $h$ maps $V$ to $V$.
\end{demo}

\subsubsection{Moduli over  $\Gamma_{\mathrm{Cardano}}\setminus\{(0,0)\}$}
\label{sec:car} 
Theorem \ref{theo:inv-gen1-1}  gives a complete set of invariants for isomorphism classes
of algebras $m$ such that $\nu(m) \not \in \Gamma_{\text{Cardano}}$.
As we remarked above, this set of invariants is the same whether $\FF$ is algebraically closed or not.
This is not the case if $\nu(m) \in  \Gamma_{\text{Cardano}}$.

Suppose now that   $\nu(m) \in \Gamma_{\text{Cardano}}$. In this case the idemvalues
of $\Q$  are  not distinct and hence 
are
in $\FF$. 
Consequently, Lemma \ref{lem:Q-P}  does not apply, {\it i.e.},
the splitting field of the fundamental cubic
$Q_m$ is not determined by the splitting field of the cubic polynomial $\Q$. 

 Given a binary cubic $Q$ we denote by
$\FF_Q$ a splitting field of $Q$. 

\begin{lemma}
\label{lem:toto}
Let $m \in \M_3^c$.
\begin{itemize}
\item[(i)] If $\nu(m) \in \Gamma_{\text{Cardano}} \setminus \big\{(0,0)\big\}$ then 
$[\FF_{Q_m}:\FF] =1$ or $2$;
\item[(ii)]  If $\nu(m)=(0,0)$ then $[\FF_{Q_m}:\FF] =1,2,3$ or $6$.
\end{itemize}
\end{lemma}
\begin{demo}
(i): Let $\FF_{Q_m}$ be a splitting field of the fundamental cubic, and let
$f_1,f_2,f_3 \in V\otimes_\FF \FF_{Q_m}$ be such 
\beqa
m(f_1,f_1) = \frac{-2\Delta-1}3 f_1 \ , \ \ 
m(f_2,f_2) = \frac{\Delta-1}3 f_2\ , \ \ 
m(f_3,f_3) = \frac{\Delta-1}3 f_3 \ .\nn
\eeqa
If $g\in \text{Gal}(\FF_{Q_m}/\FF)$ then $g$ permutes 
$\Big((f_1,-2\Delta), (f_2,\Delta),( f_3,\Delta)\Big)$ by Remark \ref{rem:galois-compat}.
Since $-2 \Delta \ne \Delta$ we must have $g\cdot f_1=f_1$ and either 
$g\cdot f_2=f_2,\ g\cdot f_3=f_3$ or $g\cdot f_2=f_3,\ g\cdot f_3=f_2$. 
This means that $\text{Gal}(\FF_{Q_m}/\FF)$ has at most two elements and hence $[\FF_{Q_m}:\FF] \le 2$.
Part two follows from Remark \ref{rem:splitting}
since the only possible degree for a splitting field of $Q_m$ is
$1,2,3$ or $6$.
\end{demo}

Given a field $\FF$ and an integer $ k \in \NN^*$, 
we denote by ${\cal E}_k(\FF)$ the set of equivalence classes of
extensions of $\FF$ of degree at most $k$ which are isomorphic to a splitting field of a degree
three polynomial with distinct roots. 
It is well known that there is a bijection
\beqa
\label{rem:E2}
{\cal E}_2(\FF) \cong \FF^*/\FF^{*2} \ ,
\eeqa
given by
\beqa
\FF(\sqrt{a}) \mapsto [a] \ . \nn
\eeqa
(Under this bijection the field $\FF$ itself corresponds to the identity element.)

\begin{theorem}
\label{theo:inv-non-gener}
 Let $\M_3^c{}_C=  \nu^{-1}\big(\Gamma_{\text{Cardano}} 
\setminus \big\{(0,0)\big\}\big)$. 
Define $\tilde \nu : \M_3^c{}_C
\to \big(\Gamma_{\text{Cardano}} \setminus \big\{(0,0)\big\} \big)\times {\cal E}_2(\FF)$
by
\beqa
\tilde \nu(m) = (\nu(m), [\FF_{Q_m}]) , \ \ 
\forall m \in  \M_3^c{}_C \ . \nn 
\eeqa
\begin{enumerate}
\item[(i)] Let $m,m'$ be in $\M^c_{3C}$. Then
\beqa
\tilde \nu(m) = \tilde \nu(m') \ \ \Rightarrow \ \ \exists g \in \mbox{GL}(V) \ s.t.\
m'=g\cdot m \ . \nn
\eeqa
\item[(ii)] The map $\tilde \nu$ is  surjective.
\end{enumerate}

\end{theorem}
\begin{demo}
(i): Since $\tilde \nu(m) = \tilde \nu(m')$ we have $\nu(m) = \nu(m')$ and
$\FF_{Q_m} \cong \FF_{Q_{m'}}$. Hence $m$ and $m'$ have the same idemvalues
and without loss of generality we can suppose $\FF_{Q_m} = \FF_{Q_{m'}}$.
Let $((f_1,-2\Delta),(f_2,\Delta),(f_3,\Delta))$ and $((f'_1,-2\Delta'),(f'_2,\Delta'),(f'_3,\Delta'))$ be two fundamental triples in
 $V \otimes_\FF \FF_{Q_m}$ s.t. 
\beqa
\begin{array}{lll}
m(f_1,f_1) = \frac{-2\Delta-1}3 f_1\ , &
m(f_2,f_2) = \frac{\Delta-1}3 f_2\ ,& 
m(f_3,f_3) = \frac{\Delta-1}3 f_3\ ;\\ \\
m'(f_1',f_1') = \frac{-2\Delta-1}3 f_1'\ , &
m'(f_2',f_2') = \frac{\Delta-1}3 f_2'\ ,& 
m'(f_3',f'_3) = \frac{\Delta-1}3 f_3'\ .
\end{array} \nn
\eeqa
Define  $h \in \text{GL} \big(V \otimes_\FF \FF_{Q_m}\big)$ by
$h(f_i)=f'_i$ then $h$ is an isomorphism of $\FF_{Q_m}-$algebras.
There are two possibilities: either $\FF_{Q_m}=\FF$ or $\FF_{Q_m}$ is a quadratic
extension of $\FF$ (see Lemma \ref{lem:toto} {\color{blue}(i)}).
 In the first case $h$ is an isomorphism of $\FF-$algebras and we are finished. In the second case
to conclude, we have to show that $h$ commutes with the action of $\text{Gal}(\FF_{Q_m}/\FF)
\cong \mathbb Z_2$. The Galois group has to fix $f_1$ and $f_1'$   and has to permute 
$f_2,f_3$ and $f_2',f_3'$. This says exactly that $h$ commutes with the action of the
Galois group.\\

(ii): We now show the surjectivity of $\tilde \nu$. Let $(p_3,p_2) \in \Gamma_{\text{Caradano}}
\setminus \{ (0,0)\}$ and let   $\KK/\FF$ be an extension of degree at most two. 
By hypothesis there exists $\Delta \in \FF^\star$ s.t.
\beqa
x^3 + p_2 x - p_3=(x+2 \Delta) (x -\Delta)^2 \ . \nn
\eeqa
Since $\KK$ is of degree at most two there exists $a \in \FF\setminus\{0\}$ s.t. $\KK$ is a 
splitting field of $x^2 -a$.
Let   $\{e_1,e_2\}$ be a basis of $V$, pick a square root $\sqrt{a}$ in $\KK$  and define the following elements of $V_\KK$:
\beqa
\label{eq:lem}
f_1 &= &\sqrt{a} e_1 + \frac13 a e_2 \ , \nn\\
f_2 &= &-\sqrt{a} e_1 + \frac13 a e_2 \ , \\
f_3&=&-\frac23 e_2 \ . \nn
\eeqa
Define the multiplication $m: V_\KK \times V_\KK \to V_\KK$  by 
\beqa
m(f_1,f_1) = \frac{\Delta -1}3 f_1\ , \ \ 
m(f_2,f_2) = \frac{\Delta -1}3 f_2\ \ ,
m(f_1,f_2) = \frac{2 + \Delta}6 \big(f_1 + f_2) \ .\nn
\eeqa
One then checks that
\beqa
\label{eq:mult-card}
m(e_1,e_1)=-\frac16  e_2 \ , \ \ 
m(e_2,e_2) = \frac{2\Delta+1}{2a} e_2\ , \ \ 
m(e_1,e_2)=\frac{\Delta -1}{2a} e_1 \ . 
\eeqa
From this it follows that  $\nu(m) = (-2 \Delta^3, -3 \Delta^2)=(p_3,p_2)$. It remains to show
that the splitting field of the fundamental cubic $Q_m$ is $\KK$. Up to a constant,
$Q_m=F_1F_2F_3 \otimes e_1 \wedge e_2$ where
\beqa
F_1(x e_1+ y e_2) =\frac13 a x -\sqrt{a} y \ , \ \ 
F_2(x e_1+ y e_2)=\frac13 a x +\sqrt{a} y \ , \ \ 
F_3(x e_1+ y e_2 )=-\frac23a x  \ 
\eeqa
since $F_i(f_i)=0, i=1,2,3$.
A short calculation then gives
\beqa
Q_m(x e_1+ y e_2)=-\frac{2}{27}a^3 x(x^2 -\frac1a y^2) \otimes e_1 \wedge e_2\ , \nn
\eeqa
which shows that $x^2-a$ and $Q_m$ have isomorphic splitting fields.
Hence  
$\tilde \nu(m) = ((-2 \Delta^3, -3 \Delta^2), \KK)=((p_3,p_2),\KK)$.

\end{demo}

\subsubsection{Moduli over  $\{(0,0)\}$:  exceptional algebras}

  Theorems \ref{theo:inv-gen1-1} and \ref{theo:inv-non-gener} give  a complete set of invariants 
for isomorphism classes
of algebras $m$ such that $\nu(m)\not = (0,0)$. We now consider the final case which corresponds
to the case of  algebras $m$ satisfying  $\nu(m)=(0,0)$ which we call exceptional.
Recall that, by definition, ${\cal E}_6(\FF)$ is the set of equivalence classes of
extensions of $\FF$ of degree at most six and which are isomorphic to splitting fields of cubic polynomials over
$\FF$. Recall also that the splitting field of a cubic polynomial is of degree three or six
{\it iff} it is irreducible.

Exceptional algebras have a very simple characterisation in terms of the representation theory of $GL(V)$.
In fact, recall that as a $GL(V)$ representation there is an isomorphism (see Proposition \ref{prop:Ac-irred})
\beqa
\M^c \cong V^\ast \oplus B  \nn
\eeqa
(here $B={\cal S}^3(V^\ast) \otimes L^{-1}$) and we now show that $m \in \M^c_3$ is exceptional {\it iff} its component along $V^\ast$ vanishes.
 \begin{lemma}\label{eq:very-simple}
Let $m$ be in $\M_3^c$ and recall that  $T_m: V \to \FF$ is  defined by $T_m(v)=\mathrm{Tr}(L_v)$ for all $v \in V$.
Then $T_m(v)=0, \forall v \in V$ iff $m$ is an exceptional algebra.
\end{lemma}
\begin{demo}
Let $((f_1,\Delta_1),(f_2,\Delta_2),(f_3,\Delta_3)) \in (V_{\FF_{Q_m}} \otimes_\FF \FF_{Q_m})^3)$ be a
reduced fundamental triple:
\beqa
m(f_i,f_i)=\frac{\Delta_i-1}3 f_i\ , \ \ 
m(f_i,f_j)= \frac{\Delta_j+2} 6 f_i + \frac{\Delta_i+2}6 f_j \ , i \ne j \ . \nn
\eeqa
Without loss of generality, we can suppose that $\FF_{Q_m} = \FF$ since
$\text{tr}(L_v \otimes_\FF \text{Id}_{\FF_{Q_m}})= \text{tr}(L_v)$. Then,
\beq  
L_{x f_1 + y f_2} = \begin{pmatrix}
\frac{\Delta_1-1}3 x + \frac{\Delta_2+2}6 y&\frac{\Delta_2+2}6x\\
\frac{\Delta_1+2}6y&\frac{\Delta_2-1}3 y + \frac{\Delta_1+2}6 x \nn
\end{pmatrix}
\eeq
and $T_m(x f_1 +y f_2)= \frac12(\Delta_1 x + \Delta_2 y)$.
Since by defintion $m$ is exceptional  {\it iff} $\nu(m)=0$ {\it iff} $\Delta_1= \Delta_2=\Delta_3=0$, this proves the result.
\end{demo} 

\begin{corollary}\label{cor:wright}
\phantom{marcus} \hskip 1.truecm
\begin{enumerate}
\item Let $m_1,m_2 \in \M^c_3$ be exceptional and let $g \in GL(V)$. 
Then   $g\cdot m_1 = m_2$ iff $g \cdot Q_{m_1} = Q_{m_2}$.
\item If $Q \in B$  is non-zero, there exists a unique $m \in \M^c_3$ such that $m$ is exceptional and $Q_m = Q$.
\end{enumerate}
\end{corollary}
\begin{demo}
This is an immediate consequence of the lemma above and   Proposition \ref{prop:Ac-irred}.
\end{demo}

\begin{theorem}
\label{theo:inv-non-gener2}
Let  $\M_3^c{}_0= \nu^{-1}\big((0,0)\big)$. Define 
 $\nu_0:  \M_3^c{}_0\to {\cal E}_6(\FF)$ by 
\beqa
\nu_0(m)= [\FF_{Q_m}] \ ,   \ \ \forall m \in  \nu^{-1}\big((0,0)\big)  \ . \nn
\eeqa
\begin{enumerate}
\item [(i)] Let
$m,m' \in\M_3^c{}_0$.  Then
\beqa
 \nu_0(m) =  \nu_0(m') \ \ \Rightarrow \ \ \exists g \in \mbox{GL}(V) \ s.t.\
m'=g\cdot m \ . \nn
\eeqa
\item[(ii)] The map  $\nu_0$ is surjective.
\end{enumerate}
\end{theorem}
\begin{demo}
(i): By Proposition 2.1(iv)  in  \cite{wright}, if $Q_1, Q_2 \in B$  are non-degenerate then they
 are in the same $GL(V)-$orbit {\it iff} their splitting
fields are isomorphic. Combining this with the previous corollary it follows that two generic exceptional 
algebras $m_1,m_2$ are in the same $GL(V)-$orbit {\it iff} the splitting fields of $Q_{m_1}$ and $Q_{m_2}$ are
isomorphic.

\noi
(ii): By definition ${\cal E}_6(\FF)$ is the set of equivalence classes of  extensions of $\FF$  of degree at
most six which are isomorphic to splitting fields of cubic polynomials over $\FF$.
Hence by Corollary  \ref{cor:wright}  (2) $\nu_0$ is  surjective.
\end{demo}

\begin{example}
Let $\KK/\FF$ be an extension of $\FF$
such that $[\KK] \in {\cal E}_6(\FF)$.
We now give explicitly the multiplication table of an exceptional algebra $m$ such that $\nu_0(m) = [\KK]$.
Without loss of generality we can assume that $\KK$ is the
 splitting field of $x^3 + d_2 x -d_3 \in \FF[X]$ and that the roots
 of this polynomial are distinct.  If
\beqa
\label{eq:rho}
x^3 + d_2 x -d_3 =(x-\rho_1)(x-\rho_2)(x-\rho_3) \ \ \text{in} \ \  \KK,
\eeqa  
 we define $f_1,f_2,f_3 \in V_\KK$ by
\beqa
\label{eq:fff}
f_1&=& \rho_1 e_1 + ( \rho_2  \rho_3 -\frac13 d_2) e_2 \ , \nn \\
 f_2&=&  \rho_2 e_1 + (  \rho_3 \rho_1 -\frac13 d_2) e_2 \ ,  \\
f_3&=& \rho_3 e_1 + (  \rho_1  \rho_2 -\frac13 d_2) e_2 \ . \nn
\eeqa
Let $m$ be the unique multiplication on $V_\KK$ s.t.
\beqa
\label{eq:Q3}
m(f_i,f_i) = -\frac{1}3 f_i\ , \ \ 
m(f_i, f_j) =\frac{1}3 f_i + \frac{1}3 f_j \ , \text{for~} i \ne j \ . 
\eeqa
One can check, after lengthy but straightforward  calculation, that
the product of two elements in $V$ is also in $V$ and 
in the basis $e_1,e_2$ we have
\beqa
\label{eq:muf}
m(e_1,e_1)&=&-\frac{1}{ 3 D^2} \Big( 9 d_3 d_2 e_1 -(27d_3^2+2d_2^3) e_2\Big) \ , \nn \\
m(e_2,e_2) &=& \frac{9}{D^2} \Big(d_3 e_1 +\frac29 d_2^2 e_2\Big) \ ,  \\
m(e_1,e_2)&=& \frac1{D^2} \Big(-2d_2^2 e_1 + 3 d_3 d_2 e_2 \Big) \ , \nn
\eeqa
where 
\beqa
\label{eq:DDiscc}
D=(\rho_1 - \rho_2) ( \rho_2 - \rho_3) ( \rho_3 -\rho_1) \ .
\eeqa
Notice that $D^2=-27 d_3^2 -4 d_2^3$  is
 the discriminant of the cubic $x^3 +d_2 x-d_3$.
It is clear  that $\nu(m)=(0,0)$ from \eqref{eq:Q3} and to see that $\nu_0(m)=\KK$
we have to show that  $\KK$ is a splitting field of $Q_m$.  Computation gives
\beqa
\label{eq:DDDDQ}
Q_m = -\frac9 {D^2}F_1 F_2  F_3 \otimes e_1 \wedge e_2 \ \,
\eeqa
where
\beqa
\label{eq:FFFF}
  F_i(xe_1 + y e_2)=( \rho_j \rho_k  -d_2/3) x - \rho_i y \ , \ \ \{i,j,k\}=\{1,2,3\} \ .
\eeqa 
Since $\rho_i$  is a root of $x^3 + d_2 x -d_3$ 
whose splitting field is  $\KK$, the  splitting field of
$Q_m$ is also $\KK$.
After some (tedious) calculation, we get
\beqa
Q_m (xe_1 + y e_2)= \Big((d_3^2 + \frac2{27} d_2^3) x^3 +  d_3 d_2 x^2 y + \frac23 d_2^2 xy^2 -d_3 y^3 \Big)\otimes e_1\wedge e_2\ . \nn
\eeqa

 Note that it follows from \eqref{eq:muf} that
\beqa
\label{eq:cube}
m(m(x e_1 + y e_2,x e_1 + y e_2),x e_1 + y e_2)=
\frac{1}{3 D^4}(d_2^2 x^2 -9d_3 xy -3 d_2 y^2)
(x e_1 + y e_2) \ , \nn
\eeqa
in other words the cube of an element in an exceptional algebra is proportional to itself.\\

Note also that  equations \eqref{eq:mult-card} define a algebra $m$ such that
$\tilde\nu(m)=((-2\Delta^3,-3\Delta^2), \FF(\sqrt{a})$. If we set $\Delta=0$,\
$m$ becomes the algebra \eqref{eq:muf}  with $d_3=0,d_2=-a$ and $D^2=4 a^3$.
This means that in some sense exceptional algebras  corresponding to extensions of degree at most
two are limits of non-exceptional algebras over the Cardano curve.

\end{example}

\begin{remark}  If $\KK / \FF$ is an extension of degree six we can describe the action of the Galois group on
the six dimensional representation $\KK$ as follows. Without loss of generality we can suppose that $\KK$ is
the splitting field of $x^3+d_2 x -N $. 
Then  we can find $\rho,\Delta \in \KK$ such that 
$\KK=\FF(\rho,{\Delta})$ where $\rho,\Delta$ are defined by
\beqa
\begin{array}{llll}
x^3+d_2 x -N &= (x-\rho)(x^2 + \rho x +\frac{N}{\rho})& (\mathrm{in~} \FF(\rho))  \\ 
&=(x-\rho)(x-\frac{-\rho +\Delta}{2})(x-\frac{-\rho -\Delta}{2})&(\mathrm{in~} \KK) \\
&=(x-\rho)(x-\bar\rho)(x-\bar {\bar \rho})
\ . \end{array} \nn
\eeqa
From $\rho \bar \rho + \bar \rho \bar {\bar \rho} + \bar {\bar \rho} \rho = d_2$
it follows that  $\Delta^2 = -3\rho^2 -4 d_2$. Since $\Delta = \bar \rho -\rho$ it follows
immediately that  $\bar \Delta = 1/2(3 \rho -\Delta)$. The action of
the Galois group  $S_3$  on $\KK$  is then given by
\beqa
(23) \longrightarrow\left\{ \begin{array}{clc} 
                             \rho &\to& \rho \\
                              {\Delta}&\to &-{\Delta}
                             \end{array} \right.\ ,  \ \ \ 
{}^-=(123)  \longrightarrow\left\{ \begin{array}{clc}
                             \rho &\to& \frac12\big(-\rho -{\Delta}\big) \\
                              {\Delta}&\to &\frac12\big(3\rho -{\Delta}\big) 
                             \end{array} \right.\ ,  \ \ \nn 
\eeqa 
and the decomposition of $\KK$  into irreducible  $\FF-$representations by
\beqa
\KK = \big< N\big> \oplus \big<D\big> \oplus \big<\rho, \bar \rho, \bar{\bar \rho}\big>
\oplus \big<\rho \bar \rho -\frac13 d_2, \bar\rho \bar{\bar \rho} -\frac13 d_2,
\bar {\bar \rho} \rho- \frac13 d_2\big> \ , \nn
\eeqa
where $D=(\rho-\bar \rho)({\bar \rho} - \bar{\bar\rho})(\bar{\bar \rho} -\rho)$.
These representations as representations of $S_3$ are respectively isomorphic to: ${\bf 1}, {\bf 1}_\epsilon, {\bf 2}$ and
${\bf 2}$. 
\end{remark}

\section{Applications}\label{sec:app}

\subsection{Division algebras}

In this section we give a criterion  for division algebras
in terms of the invariants   defined above. 
Recall that
\begin{definition}
The algebra  $m \in \M^c_3$ is a division algebra {\it iff}  for all
$v \in V \setminus \{0\}, L_v$ is invertible.
\end{definition}

\begin{lemma}\label{lem:div} (If $\FF = \mathbb R$ see \cite{ak}.) 
An  algebra $m \in \M^c_3$ is a division algebra  {\it iff} $ D_m$ is an anisotropic
quadratic form.
\end{lemma}

\begin{demo}
Let $v \in  V \setminus\{0\}$ then $v$ is $D_m-$isotropic {\it iff} $D_m(v)=0$ {\it iff}
$\det( L_v(m))=0$. This is true {\it iff} there exists  a non-zero $u \in V$ such 
$L_v(m)(u)=0$, {\it i.e.}, such that  $m(u,v)=0$. Hence $D_m$ is anisotropic
{\it iff} $m$ defines a division algebra.  
\end{demo}

We now show that the moduli $p_3,p_2$ enable us to detect when an algebra is a division algebra.

\begin{theorem}
\label{theo:crit-div}
Let $m\in \M_3^c$.
\begin{itemize}
\item [(i)] If $\nu(m)=(0,0)$ then $m$ is  a division algebra {\it iff} $-3$ is not
a square in the splitting field $\FF_{Q_m}$  of $Q_m$.
\item[(ii)] If $\nu(m) \in \Gamma_{\text{Cardano}}\setminus\{(0,0)\}$, let
 $\tilde\nu(m)=({p_3,p_2}, [\FF_{Q_m}])$.
Then $m$ is a division algebra {\it iff}
$-3(p_2 -p_3+ 1)$ is not a square in  $\FF$ when  $\FF_{Q_m} \cong \FF$ or {\it iff}
$-3a(p_2-p_3 +1)$  is not a square in  $\FF$ when $[\FF_{Q_m}:\FF]=2$ and
 $\FF_{Q_m} \cong\FF(\sqrt{a})$.  
\item[(iii)] If $\nu(m) \not \in \Gamma_{\text{Cardano}}$, let $\nu(m)=(p_3,p_2)$. Then $m$ is 
a division algebra {\it iff} $
-3(p_2-p_3+1)(-27p_3^2- 4p_2^3)$ is not a square in $\FF$.
\end{itemize}
\end{theorem}
\begin{demo}
(i): Since $\nu(m)=(0,0)$ the equivalence class of $m$
is completely determined by  $\nu_0(m)=\KK$. Let $x^3+d_2 x - d_3 \in \FF[X]$ be a polynomial 
with three distinct roots $\rho_1,\rho_2,\rho_3$
 whose splitting field is isomorphic to $\FF_{Q_m}$. By 
\eqref{eq:muf} there exists a basis $\{e_1,e_2\}$  
of $V$ such that the matrix of the left multiplication by $x e_1+ye_2$ is given by
\beqa
L_{xe_1 + y e_2} =\begin{pmatrix}
-\frac{3x}{D^2} d_3 d_2  -\frac{2y}{D^2}d_2^2&-\frac{2x}{D^2}d_2^2  + \frac{9y}{D^2} d_3  \\
\frac{x}{D^2}(9 d_3^2 + \frac23 d_2^3) + \frac{3y}{D^2} d_3 d_2  
& \frac{3x}{D^2} d_3 d_2 + \frac{2y}{D^2} d_2^2 
\end{pmatrix} \nn 
\eeqa
where $D^2=(\rho_1-\rho_2)^2(\rho_2-\rho_3)^2(\rho_3-\rho_1)^2= -27 d_3^2 -4 d_2^3$.
By Lemma \ref{lem:div}, $m$ defines a division algebra
{\it iff} $D_m(x,y)=\det(L_{xe_1+ye_2})$ is an anisotropic quadratic form. Computing the determinant we get
\beqa
D_m(x,y)=\frac{27 d_3^2 + 4 d_2^3}{3 D^4} \Big(
d_2^2 x^2 -9 d_3 x y - 3 d_2 y^2 \Big) \nn
\eeqa
and the discriminant reduces to
$
-3\frac1{ (3D)^2}. 
$
Since $D \in \FF_{Q_m} \setminus \FF$ this proves (i).

(ii-a): When $\FF_{Q_m}\cong \FF$ there exists a reduced fundamental triple 
$((f_1,\Delta_1=\Delta),(f_2,\Delta_2=\Delta),(f_3,\Delta_3=-2\Delta)) \in (V\times \FF)^3.$
The fundamental quadratic $D_m(x f_1+yf_2)$ is given by
\beqa
 D_m(x f_1+yf_2)=\frac19 x^2(\Delta_1-1)(1+\frac12\Delta_1)+\frac19 xy(\Delta_1-1)(\Delta_2-1) +
\frac19 y^2(\Delta_2-1)
(1 + \frac12\Delta_2) \ , \nn
\eeqa
and its discriminant by
$
\frac1{27}(p_3-p_2-1).
$
This proves (ii-a).

(ii-b).  By \eqref{eq:mult-card} there exists a basis $\{e_1,e_2\}$ of $V$ such that the matrix of left multiplication
by $xe_1 + y e_2$ is given by
\beqa
L_{xe_1 + y e_2} = \begin{pmatrix}\frac{\Delta-1}{2a} y & \frac{\Delta-1}{2a} x\\
                     -\frac16 x & \frac{2 \Delta +1}{2a} y 
\end{pmatrix} \ . \nn
\eeqa
By Lemma  \ref{lem:div} $m$ defines a division algebra {\it iff} $\det( L_{xe_1 + y e_2})$ is a anisotropic quadratic form.
A short calculation shows that
\beqa
\det( L_{xe_1 + y e_2}) = \frac{p_3-p_2-1}{12 a^3} \ , \nn
\eeqa
and this proves (ii-b).\\
(iii): If $\nu(m)=(p_3,p_2)$ then by \eqref{eq:mul-gen} there is a basis $\{e_1,e_2\}$
 of $V$ such that the matrix of the left multiplication by $x e_1+ye_2$ is given by
\beqa
L_{xe_1+ye_2}=\begin{pmatrix}
\frac{x}{3d^2}\Big(-4p_2^3-27p_3^2-9p_2p_3\Big) -\frac{2y}{d^2}  p_2^2&
-\frac{2x}{d^2} p_2^2 + \frac{9y}{d^2} p_3\\
\frac{x}{3d^2}\Big(2p_2^3+27 p_3^2\Big) + \frac{3y}{d^2}\Big(-\frac29p_2^3-\frac32p_3^2+p_2p_3\Big)&
\frac{2x}{d^2}\Big(-\frac29p_2^3-\frac32p_3^2+p_2p_3\Big) + \frac{2y}{d^2} p_2^2
\end{pmatrix}\nn
\eeqa
with $d^2=(\Delta_1-\Delta_2)^2(\Delta_2-\Delta_3)^2(\Delta_3-\Delta_2)^2 
= - (27 p_3^2 + 4 p_2^3)$. By Lemma \ref{lem:div}, $m$ defines a division algebra
iff $D_m(x,y)=\det(L_{xe_1+ye_2})$ is an anisotropic quadratic form. Computing the determinant we get
\beqa
D_m(x,y) = \frac{4p_2^3 + 27 p_3^2}{18d^4} \Bigg\{
 \Big(27p_3^2 -9p_2 p_3 + 6 p_2^2 + 4 p_2^3\Big) x^2-
6\Big(9p_3 +2 p_2^2\Big) x y +
9\Big(3p_3 -2p_2\Big) y^2  \Bigg\}\ , \nn
\eeqa
and then the discriminant of $D_m(x,y)$ reduces to
$
-\frac19\times3(p_2-p_3+1)(27p_3^2 + 4p_2^3) \ . \nn
$
This proves (iii).
\end{demo}

Part one of this Theorem implies that there is (up to isomorphism) only one two-dimensional real commutative
division algebra whose automorphism group is $S_3$. In fact one can show that any
real division algebra whose automorphism group is $S_3$ must be commutative \cite{dz}.

\subsection{Associativity}
In this section we show that generic commutative associative algebras correspond to a
single point in the Cardano plane.

\begin{theorem}\label{theo:assos} 
Let $m\in \M_3^c$. Then
$
m \text{~is~associative~iff~} \nu(m)=(16,-12).
$
Note that $\nu(m) \in \Gamma_{\text{Cardano}}$.
\end{theorem}
\begin{demo}
Let $\FF_{Q_m}$ be a splitting field of $Q_m$ and let $f_1, f_2 , f_3 \in V_\KK$
be the fundamental triple:
\beqa
m(f_1,f_1) =\frac{\Delta_1 -1}3 f_1\ ,
m(f_2,f_2) =\frac{\Delta_2 -1}3 f_2\ ,
m(f_1,f_2) = \frac{2+ \Delta_2}6 f_1 +  \frac{2+ \Delta_1}6 f_2 \ . \nn
\eeqa
Among the eight associators $A(f_i,f_j,f_k), i,j,k=1,2$ one checks that
only  two do not trivially vanish:
\beqa
m(m(f_2,f_1),f_1)- m(f_2,m(f_1,f_1))\ ,     \ \
m(m(f_1,f_2),f_2)- m(f_1,m(f_2,f_2)) \ . 
\eeqa
Calculation shows that
\beqa
m(m(f_2,f_1),f_1)- m(f_2,m(f_1,f_1))&=& \frac{2 + \Delta_1}{36}
\Big((2+\Delta_2) f_1 + (4-\Delta_1) f_2\Big) \ , \nn \\ 
m(m(f_1,f_2),f_2)- m(f_1,m(f_2,f_2))  &=& \frac{2 + \Delta_2}{36}\Big(
(4-\Delta_2) f_1 + (2+\Delta_1) f_2 \Big) \ . \nn
\eeqa 
It follows from this  that  $m$ is associative iff 
 $(\Delta_1,\Delta_2)=(-2,-2),(-2,4)$ or $(4,-2)$ and hence  iff  $\nu(m)=(16,-12)$. 
\end{demo}

We have already  seen  in Theorem \ref{theo:inv-non-gener} 
that equivalence classes of 
algebras such that $\nu(m)=(16,-12)$ are completely determined by the class
of the field extension $\FF_{Q_m}/\FF$. Moreover,
since $(16,-12) \in \Gamma_{\text{Cardano}}\setminus\{(0,0)\}$, $\FF_{Q_m}$ can be any
extension of $\FF$ of degree at most two.
\begin{corollary}
Let $m \in \M_3^c$ be  associative and let $\tilde \nu(m)=((16,-12),[\KK])$ where
$\KK /\FF$ is an extension of degree one or two. 
\begin{itemize}
\item [(i)] If $\KK \cong \FF$, then the $\FF-$algebra $(V,m)$   is isomorphic
to $\FF \times \FF$.
\item[(ii)] If $\KK$ is a quadratic extension, then the $\FF-$algebra $(V,m)$  
 is isomorphic to $\KK$.
\end{itemize}
\end{corollary}
\begin{demo}
(i): If $\KK \cong \FF$, the splitting field of $Q_m$ is $\FF$. Hence, as we saw
in the proof of Theorem \ref{theo:assos}  one can suppose $\Delta_1=\Delta_2=-2,\Delta_3=4$
and hence there exist
$f_1, f_2 \in V$ such that
\beqa
m(f_1,f_1)=-f_1 \ , \ 
m(f_2,f_2)=-f_2 \ ,
m(f_1,f_2)=0\ , \nn
\eeqa
and the algebra $(V,m)$ is  the product of the ideals
$\text{Vect}(f_1) \times \text{Vect}(f_2)$. The identity element 
is $f_3= -f_1 -f_2$.

(ii): The splitting field  $\FF_{Q_m}$ of $Q_m$ is an extension of degree two so
$\text{Gal}(\FF_{Q_m}/\FF) \cong \mathbb Z_2$ and there exist $\alpha \in V^*$, $\beta \in V^*_{\FF_{Q_m}}$ and $b \in L^{-1}$
such that $Q_m=\alpha \beta \bar \beta \otimes b$.  One can choose 
$f \in \text{Ker}(\beta)$ such that  $((f,-2), (\bar f,-2),(-f-\bar f ,4))$ is a reduced
fundamental triple (see Theorem \ref{theo:can}), thus 
\beqa\label{eq:prod2}
m(f,f)=-f \ , \ 
m(\bar f,\bar f)=-\bar f \ ,
m(f,\bar f)=0\ , 
\eeqa
Since $\FF_{Q_m}$ is a quadratic extension of $\FF$ there exists
$a \in \FF^* \setminus (\FF^*)^2$ ({\it i.e.}, $a$ is not a square in $\FF^*$)
such that $\FF_{Q_m}=\FF(\sqrt{a})$.
In the  basis $e_1= -(f+\bar f), e_2=\sqrt{a}(f-\bar f)$ of $V$ the multiplication
defined in \eqref{eq:prod2} reduces to
\beqa
m(e_1,e_1) =  e_1\ ,\ \ 
m(e_2,e_2) = a e_1 \ , \ \ 
m(e_1,e_2)=  e_2 \ . \nn
\eeqa
Here the map $e_1 \to 1, e_2 \to \sqrt{a}$ defines an  $\FF-$algebra isomorphism of
 $(V,m)$ with $\FF_{Q_m}$.
\end{demo}

\begin{remark}
In case (ii) of the corollary  above,
 the algebra $(V,m)$ is isomorphic to a field and in particular is a division algebra. This means that
the criterion of Theorem \ref{theo:crit-div} (ii) must be satisfied which is true since 
\beqa
-3a(p_2 -p_3 +1) = 81 a \  \nn
\eeqa
is not a square in $\FF$.
\end{remark}

\section{Moduli space of generic algebras with respect to the action of $SL(V)$}\label{sec:mod-sl}
In this section we consider the set of generic algebras
\beqa
\M^c_3 = \Big\{ m \in \M^c\ , \ \  \text{s.t.} \ \  \text{Disc}(Q_m) \ne 0 \Big\} \ , \nn
\eeqa
modulo the action of $SL(V)$. Since $SL(V)$ is a subgroup of $GL(V)$ there is a natural map from the $SL(V)-$moduli space 
$\M^c_3/SL(V)$ to  the $GL(V)-$moduli space ${\cal M}^c_3$.
It turns out that the description of the fibre of this map is quite different depending on whether or not we are working
over the Cardano curve.
We will show that the invariant
$\frac{\text{Inv}_m}{\text{Disc}(Q_m)}$ distinguishes $SL(V)-$moduli corresponding to the same $GL(V)-$moduli
if the latter  do not  belong
 to the Cardano curve. If they do belong to the Cardano curve,  given by the vanishing of Inv (see Theorem \ref{theo:miracle}),
we will show that the invariant Disc$(Q_m)$  distinguishes the corresponding $SL(V)-$moduli unless Disc$(Q_m)$ is irreducible
(see (3.2) of Proposition \ref{theo:inv-ls-1}).

\begin{proposition}
\label{theo:inv-ls-1}
Let $m,m' \in \M_3^c$ and let $\Gamma_{\mathrm{Cardano}} = \Big\{ (p_3,p_2) \in \mathbb F^2: 27p_3^2 +4p_2^3=0\Big\}$.
\begin{enumerate}
\item If $\nu(m), \nu(m') \not \in \Gamma_{\mathrm{Cardano}}$
and  $\nu(m) = \nu(m')$,
then  there exists $g \in SL(V)$ such that $m'=g \cdot m$ {\it iff}
 $\frac{\mathrm{Inv}_m}{\mathrm{Disc}(Q_m)}=\frac{\mathrm{Inv}_{m'}}{\mathrm{Disc}(Q_{m'})}$.
\item If $\tilde\nu(m), \tilde \nu(m')  \in \Gamma_{\mathrm{Cardano}}\setminus\{0\}$
and  $\tilde \nu(m) = \tilde\nu(m')$,
then  there exists $g \in SL(V)$ such that $m'=g \cdot m$ {\it iff}
 $\mathrm{Disc}(Q_m)=\mathrm{Disc}(Q_{m'})$.
\item If $ \nu(m)=\nu(m')=(0,0) $
and  $\nu_0(m) =  \nu_0(m')$,
\begin{enumerate}
\item [3.1] if $Q_m$ is reducible (equivalently if $[\FF_{Q_m}:\FF]=1$ or $2$)
then  there exists $g \in SL(V)$ such that $m'=g \cdot m$ {\it iff}
 $\mathrm{Disc}(Q_m)=\mathrm{Disc}(Q_{m'})$, 
 \item[3.2]  if $Q_m$ is irreducible (equivalentlyf if  $[\FF_{Q_m}:\FF]=3$ or $6$) then Disc$(Q_m)$ determines
 the $SL(V)-$classes up to a $\mathbb Z_2-$factor.
 \end{enumerate}
\end{enumerate}
\end{proposition}
\begin{demo}
The implications ``$\Rightarrow$'' are immediate for (1),(2),(3), so we now prove the three converse implications.

(1):   Suppose $m,m'\in \M_3^c$ are such that
$\nu(m) ,\nu(m') \not \in \Gamma_{\mathrm{Cardano}}$ and 
$$\frac{\mathrm{Inv}_m}{\mathrm{Disc}(Q_m)}=\frac{\mathrm{Inv}_{m'}}{\mathrm{Disc}(Q_{m'})}\ , \ \
\nu(m)=\nu(m')\ .
$$
By Theorem \ref{theo:inv-gen1-1} there exists $g$ in $GL(V)$ such that $m'=g \cdot m$. 
Hence 
\beqa
\frac{\mathrm{Inv}_{m'}}{\mathrm{Disc}(Q_{m'})}= (\det g)^{-1}\frac{\mathrm{Inv}_{m}}{\mathrm{Disc}(Q_{m})} 
= \frac{\mathrm{Inv}_{m}}{\mathrm{Disc}(Q_{m})} \ , \nn
\eeqa 
which means that $g$ is in $SL(V)$.\\

(2): Suppose $m,m'\in \M^c$ are such that   $\tilde \nu(m), \tilde \nu(m') \in  \Gamma_{\mathrm{Cardano}}\setminus\{0\}$ and
 $$\mathrm{Disc}(Q_m)=\mathrm{Disc}(Q_{m'})\ , \ \ \tilde\nu(m)= \tilde \nu(m')\ .$$ 
 By Theorems \ref{theo:inv-non-gener} and \ref{theo:inv-non-gener2} there exists $g$ in $GL(V)$ such that $m'=g \cdot m$. 
Hence  
\beqa
\mathrm{Disc}(Q_{m'})= (\det g)^{-2}\mathrm{Disc}(Q_{m})
= \mathrm{Disc}(Q_{m}) \ , \nn
\eeqa 
whence $\det(g) = \pm 1$. Without loss of generality we can assume $\det(g)=1$ since by Corollary \ref{cor:aut}
there exists an automorphism of $m'$ of determinant minus one.

The proof of (3.1) is  the same as the proof of (2) (see corollary \ref{cor:aut}).

(3.2): By Theorem \ref{theo:inv-non-gener2}, $m$ and $m'$ are $GL(V)-$equivalent {\it iff} the twisted binary cubics
$Q_m$ and $Q_{m'}$ are $GL(V)-$equivalent. It follows that $m$ and $m'$ are $SL(V)-$equivalent {\it iff}
$Q_m$ and $Q_{m'}$ are $SL(V)-$equivalent. As an $SL(V)-$representation the space of twisted binary cubics and the
space of binary cubics are  isomorphic so we are led to the problem of classifying irreducible (untwisted) binary
cubics under the action of $SL(V)$.

Let us first consider the particular case where in some basis $\{e_1,e_2\}$ of $V$ we can write 
\beqa
q_m(x,y)=a x^3 + b y^3 \ , \ \
q_{m'}(x,y)=a' x^3 + b' y^3 \ , \nn
\eeqa
with $Q_m = q_m e_1 \wedge e_2, Q_{m'} = q_{m'} e_1 \wedge e_2$.
Then one can show that
\beqa
\text{Disc}(Q_m) = -3\big(3ab \otimes(\epsilon_1\wedge \epsilon_2)^2\big) \ , \ \ \text{Disc}(Q_{m'}) =
-3\big(3a'b' \otimes(\epsilon_1\wedge \epsilon_2)^2\big) \ ,\nn
\eeqa
where $\{\epsilon_1,\epsilon_2\}$ is the basis dual to $\{e_1,e_2\}$.
It is shown in  \cite{ss} that
$q_{m}$ and $q_{m'}$ are in the same $SL(V)-$orbit {\it iff} 
\beqa
\left\{\begin{array}{ccc}
\big[\frac{a}{b}\big]_{\FF^*/\FF^{* 3} }&=& \big[\frac{a'}{b'}\big]_{\FF^*/\FF^{* 3}}   \nn\\
ab&=&a' b'
\end{array}
\right. \ \text{or} \ \ \ \
\left\{\begin{array}{ccc}
\big[\frac{a}{b}\big]_{\FF^*/\FF^{* 3} }&=& \big[\frac{b'}{a'}\big]_{\FF^*/\FF^{* 3}}   \nn\\
ab&=&-a' b'
\end{array}\right.
\ . \nn
\eeqa
Now the fact that $Q_m$ and $Q_{m'}$ are in the same $GL(V)-$orbit implies that there exists $\lambda$ in $\FF^*$ such that
$\lambda(ab) = a'b'$, and 
the fact that Disc$(Q_m)$=Disc$(Q_{m'})$ implies that  $ab = \pm a'b'$, {\it i.e.},  $\lambda^2=1$.
If $\lambda=1$ we must have
\beqa
\left\{\begin{array}{ccc}
\big[\frac{a}{b}\big]_{\FF^*/\FF^{* 3} }&=& \big[\frac{a'}{b'}\big]_{\FF^*/\FF^{* 3}}   \nn\\
ab&=&a' b'
\end{array}
\right. \ \text{or} \ \ \ \
\left\{\begin{array}{ccc}
\big[\frac{a}{b}\big]_{\FF^*/\FF^{* 3} }&=& \big[\frac{b'}{a'}\big]_{\FF^*/\FF^{* 3}}   \nn\\
ab&=&-a' b'
\end{array}\right.
\ , \nn
\eeqa
and  if $\lambda=-1$ we must have
\beqa
\left\{\begin{array}{ccc}
\big[\frac{a}{b}\big]_{\FF^*/\FF^{* 3} }&=& \big[\frac{a'}{b'}\big]_{\FF^*/\FF^{* 3}}   \nn\\
ab&=&-a' b'
\end{array}
\right. \ \text{or} \ \ \ \
\left\{\begin{array}{ccc}
\big[\frac{a}{b}\big]_{\FF^*/\FF^{* 3} }&=& \big[\frac{b'}{a'}\big]_{\FF^*/\FF^{* 3}}   \nn\\
ab&=&a' b'
\end{array}\right.
\ . \nn
\eeqa
Hence either $Q_{m'}$ and $Q_m$ are in the same orbit or $Q_{m'}$ is in the $SL(V)-$orbit corresponding to the parameters
$([a/b], -ab)$.
This proves (3.2) in the special case where $Q_m$ and $Q_{m'}$
can be written as above. In fact, as shown in \cite{ss},   this is the general situation, and this completes the proof of part (3.2).
\end{demo}

\begin{remark}
Since
\beqa
\mathrm{Disc}(D_m)=\frac1{27}(p_3 - p_2 -1) \mathrm{Disc}(Q_m) \ , \nn
\eeqa
the invariant Disc$(D)$ also distinguishes $SL(V)-$modules with the same $GL(V)-$moduli on the Cardano curve
{\it if} $p_3 - p_2 -1 \not = 0$ if $[\FF_{Q_m}/\FF] \le 2$. Note that the only point on the Cardano curve such that 
$p_3 - p_2 -1 =0$ is the ``associative'' point $(16,-12)$ (see Theorem \ref{theo:assos}).
\end{remark}

Before stating the next theorem we need the

\begin{lemma}
The natural action  of the group $\FF^{*2}$ on $L^2\setminus\{0\}$ defines a 
principal  $\FF^{*2}-$fibration  $s: L^2\setminus\{0\} \to \FF^*/\FF^{*2}$.  
\end{lemma}
\begin{demo}
It is clear that 
\beqa
\omega\otimes \omega \mapsto (\lambda \omega)\otimes    (\lambda \omega) \ , \ \ \forall \lambda \in \FF^*, \omega\in
L \ , \nn
\eeqa
is a well defined  principal $\FF^{*2}-$action, {\it i.e.}, the only element of $\FF^{*2}$ with a fixed point is the identity. 
The point of the Lemma is to identify the quotient space with $\FF^*/\FF^{*2}$. For this  pick $\{\omega_0\}$  of $L$
and define
\beqa
\begin{array}{lclc}
s: &L^2 \setminus\{0\} &\to &\FF^*/\FF^{*2}\\
&(\lambda\omega_0)\otimes(\lambda \omega_0)&\mapsto& [\lambda^2] \ .\nn
\end{array}
\eeqa
One easily checks that  $s$ is independent of the choice of basis $\{\omega_0\}$ and that
$s(\omega_1\otimes \omega_1) = s(\omega_2\otimes \omega_2)$ {\it iff} there exists $\lambda \in \FF^{*}$ such that
$\omega_2 \otimes \omega_2  = \lambda^2 \omega_1 \otimes \omega_1$. This proves the lemma.

\end{demo}

\noi
This lemma enables us to associate a half-line (more precisely a $\FF^{*2}-$orbit) in $L^2\setminus \{0\}$ to
any quadratic extension of $\FF$ since $\FF^*/\FF^{*2}$ is naturally identified with the set of at most quadratic extensions of $\FF$.

\begin{definition}
Let $[a]$ be an element of $\FF^*/\FF^{*2}$. We define the half-line $L_{[a]} \subset L^2\setminus \{0\}$ by
\beqa
L_{[a]} = s^{-1}([a]) \ . \nn
\eeqa

\end{definition}

In part (3) of  the following theorem we use the  notation:\\

\begin{itemize}

\item \ \ $ \mathbf E_i$ is  the set of equivalence classes of Galois
extensions of $\mathbb F$ of degree  equal to
 $i$.

\item \ \  ${}^i \hskip -.15truecm \M^c_{3}{}_0 = \Big\{ m \in \M^c_{3}{}_0 \ \ \text{s.t.} \ \  \nu_0(m) \in \mathbf E_i\Big\}$.


\item \ \ $\mathbf \Delta_{6,2} = \Big\{(\tilde \FF_6, \omega) \in {\mathbf E}_6 \times (L^2 \setminus\{0\})
 \ \ \text{s.t.} \ \  s(\omega)   \subset \tilde \FF_6 \Big\}$\ ,
\end{itemize}
 where $i$ takes the values $1,2,3$ or $6$.
 
\begin{theorem}\label{theo:slmodule}\hskip 1.truecm \phantom{toto} 

\begin{enumerate} 
\item The map $\hat \nu_s: \M^c_{3g}/SL(V)\to \Big(\FF^2 \setminus \Gamma_{\mathrm{Cardano}}\Big)
\times(L\setminus\{0\})$ given by 
$$\hat \nu_s ( m) = \Bigg(p_3,p_2,\frac{\mathrm{Inv}}{\mathrm{Disc}(Q)}\Bigg) $$
 is a bijection.
\item The map $\hat{\hat \nu}_s:\M^c_{3C}/SL(V)\to
\Big(\Gamma_{\mathrm{Cardano}}\setminus\{0\} \Big) \times (L^2\setminus\{0\} )$ given by
\beqa 
\hat {\hat \nu}_s(m)= (p_3,p_2,\mathrm{Disc}(Q))\nn
\eeqa 
is a bijection.
\item  If $\nu(m)=(0,0)$ there are three cases: 
\begin{enumerate}
\item[3.1] The map 
${}_1\hat \nu_0 :  \Big({}^1 \hskip -.15truecm \M^c_{3}{}_0 \cup {}^3 \hskip -.15truecm \M^c_{3}{}_0\Big)\Big/SL(V)
\to \Big({\mathbf E_1} \cup {\mathbf E_3}\Big) \times L_{[1]}$  given by
\beqa
{}_1\hat \nu_0(m) = (\nu_0(m),\text{Disc}(Q_m)) \ , \nn
\eeqa 
is a bijection.
 \item[3.2] The map
${}_2\hat \nu_0 : {}^2 \hskip -.15truecm \M^c_{3}{}_0/SL(V) \to 
(L^2 \setminus \{0\}) \setminus L_{[1]}$ given by 
\beqa
{}_2\hat \nu_0(m) = \text{Disc}(Q_m) \ , \nn
\eeqa 
is a bijection. 
\item[3.3]  The map
${}_6\hat \nu_0 : {}^6 \hskip -.15truecm \M^c_{3}{}_0/SL(V)  \to  {\mathbf \Delta}_{6,2} $ given by 
\beqa
{}_6\hat \nu_0(m) = (\nu_0(m),\text{Disc}(Q_m)) \ . \nn
\eeqa 
\end{enumerate}

\end{enumerate}
\end{theorem}

\begin{demo}
(1): By (1) of Proposition \ref{theo:inv-ls-1}, the map is injective so it remains to  prove  surjectivity. Taking $e=f=b=0$ and
$d=1$ in the formul\ae \  \eqref{eq:decom} and \eqref{eq:Invd}, we get
\beqa
\frac{\mathrm{Inv}_m}{\mathrm{Disc}(Q_m)}= \frac{2a^2 c}{1-4ac}  (\epsilon_1 \wedge \epsilon_2) \ . \nn
\eeqa
Let $X\in \FF\setminus\{0\}$. Then one checks immediately that if we take the multiplication $m$ corresponding to
\beqa
a= 4 X \ , \ \ 
c= \frac 1{48 X} \ , \nn
\eeqa
then  Inv$_m/$Disc$(Q_m)= X \epsilon_1 \wedge \epsilon_2$ and this proves
 surjectivity.\\

(2): To prove  surjectivity in this case, we start with two lemmas.

 Since $(p_3,p_2) \in \Gamma_{\mathrm{Cardano}}\setminus\{0\}$ we know that the splitting field
$\FF_{Q_m}$ of $Q_m$ is of degree at most two  (see Lemma \ref{lem:toto}).
\begin{lemma}
\label{lem:disc-ext}
 Write  $\FF_{Q_m}= \FF(\sqrt{a})$ where $a \in \FF$ is a square  in $\FF$ if $[\FF_{Q_m}:\FF]=1$ and
 not a square in $\FF$  if $[\FF_{Q_m}:\FF]=2$.
Then 
\beqa
s(\mathrm{Disc}(Q_m))=[a] \ . \nn
\eeqa
\end{lemma}
\begin{demo}
Since $\tilde \nu(m) = (p_3,p_2,[\FF_{Q_m}])$, as we showed in \eqref{eq:lem},   there is a basis 
of idempotents $\{f_1,f_3\}$  of $V\otimes \FF_{Q_m}$ 
such that 
\beqa
f_2= -\sqrt{a} e_1 +\frac13 a e_2 \ ,\ \  f_3=-\frac23 e_2 \ . \nn
\eeqa
The dual basis is
\beqa
\varphi_2= -\frac 1 {\sqrt{a}} \epsilon_1 \ , \ \ 
\varphi_3= -\frac {\sqrt{a}} 2 \epsilon_1 -\frac32 \epsilon_2 \ , \nn
\eeqa
and hence from  
 \eqref{eq:Inv}
\beqa
\mathrm{Disc}(Q_m) = (f_2 \wedge f_3)^2=\frac 9{4 a} (\epsilon_1 \wedge \epsilon_2)^2  \ .\nn  
\eeqa
By definition
\beqa
s(\mathrm{Disc}(Q_m)) = \Big[\frac 9 {4a}\Big]= [a]    \nn
\eeqa
and this proves the lemma.
\end{demo}

Recall now (see Theorem \ref{theo:inv-non-gener})
that $\tilde \nu: \M^c_{3C}/GL(V) \mapsto \Gamma_{\mathrm{Cardano}}\setminus\{0\} \times {\cal E}_2(\FF)$ is  a bijection.
In the following Lemma, by abuse of notation, we have identified ${\cal E}_2(\FF)$ with $\FF^*/\FF^{*2}$ as in \eqref{rem:E2}
\begin{lemma}
The following diagram is commutative
\beqa
\xymatrix{
\hskip 1.truecm \M^c_{3C}/SL(V) \phantom{totot} \ar[r]^{\hskip -1.truecm \hat {\hat \nu}_s}\ar[d]_{p} & 
\Big(\Gamma_{\mathrm{Cardano}}\setminus\{0\}\Big) \times (L^2 \setminus\{0\})\ar[d]^{\mathrm{Id}\times s}\\
\hskip .5truecm \M^c_{3C}/GL(V)\ar[r]_{\tilde \nu}\phantom{totot}
 \phantom{to}&\hskip .2truecm
\Big(\Gamma_{\mathrm{Cardano}}\setminus\{0\} \Big)\times \FF^*/\FF^{*2} \\
 }
\nn
\eeqa
\end{lemma}
\begin{demo}
Let $m \in \M^c_{3C}$. We denote by  $[m]_{SL}$ and  $[m]_{GL}$ respectively, the corresponding equivalence classes with respect to 
the actions of  $SL(V)$ and  $GL(V)$. By definition, 
\beqa
\tilde \nu\big([m]_{GL}) = (p_3,p_2,\FF_{Q_m})  \nn
\eeqa
where $\FF_{Q_m}$ is the splitting field of $Q_m$. 
Identifying ${\cal E}_2(\FF)$ with $\FF^*/\FF^{*2}$ (see \eqref{rem:E2}) we can write this as
\beqa
\tilde \nu\big([m]_{GL}) = (p_3,p_2,\FF_{Q_m})=  (p_3,p_2,[a]) \ ,\nn
\eeqa
if $\FF_{Q_m} = \FF(\sqrt{a})$.
By Lemma \ref{lem:disc-ext} this is equivalent to
\beqa
\tilde \nu\big([m]_{GL}) = (p_3,p_2,s(\text{Disc}(Q_m))\ \nn
\eeqa
and  by definition the RHS of this equation is equal to 
\beqa
 (\text{Id}\times s)\circ \hat{\hat \nu}_s \big([m]_{SL}) \ .\nn 
\eeqa
This shows that the diagram is commutative.
\end{demo}

We now show that $\hat {\hat \nu}_s$ is surjective. Let $I=\text{Im}(\hat{\hat \nu}_s) \subset 
(\Gamma_{\mathrm{Cardano}}\setminus\{0\}) \times (L^2 \setminus\{0\})$. On the one hand,
Id$\times s(I)$ is equal to $(\Gamma_{\mathrm{Cardano}}\setminus\{0\}) \times \FF^*/\FF^{*2}$ since 
$\tilde\nu$ is surjective and the above diagram commutes. On the other hand,
$I$ is stable under the action of $\FF^{*2}$ since 
\beqa
\hat{\hat \nu}_s\big([m]_{SL}\big)&=& (p_3,p_2, \text{Disc}(Q_m))
\nn\\
\hat{\hat \nu}_s\big([\lambda \text{Id} \cdot m]_{SL}\big)&=& (p_3,p_2,\lambda^2 \text{Disc}(Q_m)) \ . \nn
\eeqa
 for any $\lambda \in \FF^*$. It follows that $I=(\Gamma_{\mathrm{Cardano}}\setminus\{0\}) \times (L^2 \setminus\{0\})$
and $\hat{\hat \nu}_s$ is surjective.\\

(3): We will need the following  lemma.
\begin{lemma}
\label{lem:marcus}
Let $m$ be an element of $\M^c_{3}{}_0$,  let $\FF_{Q_m}$ be a splitting field of $Q_m$ and let Disc$(Q_m)$  be the discriminant of $Q_m$.
Then 
\beqa
[\FF_{Q_m}:\FF] =  \left\{
\begin{array}{cll}
 1 \ \ \mathrm{or} \ \  3 &\Rightarrow& s\big(\mathrm{Disc}(Q_m)\big) = [1] \ ,\\
 2   &\Rightarrow& s\big(\mathrm{Disc}(Q_m)\big) =[\FF_{Q_m}] \  , \\
 6  &\Rightarrow& s\big(\mathrm{Disc}(Q_m) \big)= [\tilde\FF] \ , \ \mathrm{where} \ \  \tilde \FF \ \mathrm{is \ the \ unique} \\
&& \mathrm{quadratic \ extension\ of}\  \FF \ \mathrm{contained\ in}\  \FF_{Q_m} \ . \nn
\end{array}\right.
\eeqa
\end{lemma}
\begin{demo}
We can always find a polynomial of the form $P(x)=x^3+d_2 x-d_3$ in $\FF[X]$ such that $\FF_{Q_m}$ is the splitting field of
$P$.
By Theorem \ref{theo:inv-non-gener2} there exist $\rho_1,\rho_2,\rho_3 \in \FF_{Q_m}$ and $f_1,f_2,f_3$ in $V\otimes\FF_{Q_m}$
such that the  multiplication table of $m$ in the basis $\{f_1,f_2\}$ is given by  \eqref{eq:fff}.
In particular  by \eqref{eq:DDDDQ} and \eqref{eq:FFFF}
\beqa
Q_m&=&\frac{-1}{3D^2}\Big( (\phantom{-}2\rho_1\rho_2-\rho_2\rho_3 -\rho_3 \rho_1)x - \rho_3 y\Big)\nn\\
&&\hskip .85truecm \Big((-\rho_1\rho_2 +2\rho_2\rho_3 -\rho_3 \rho_1)x - \rho_1 y\Big) \nn \\
&&\hskip .85truecm \Big((-\rho_1\rho_2-\rho_2\rho_3 +2\rho_3 \rho_1)x - \rho_2 y\Big) (f_1\wedge f_2) \ ,\nn
\eeqa
where $D=(\rho_1-\rho_2)(\rho_2-\rho_3)(\rho_3-\rho_1)$. Note that $D^2$ is in $\FF$.
  By \eqref{eq:Inv} we know that 
\beqa
\text{Disc}(Q_m)= (\varphi_1 \wedge \varphi_2)^{\otimes^2} \ , \nn
\eeqa
with $\{\varphi_1,\varphi_2\}$ the dual basis of $f_1,f_2$.
On the other hand it follows from  \eqref{eq:fff} that 
\beqa
 (\varphi_1 \wedge \varphi_2)= \frac3 { (\rho_1-\rho_2) (\rho_2-\rho_3) (\rho_3-\rho_1)}\;(\epsilon_1 \wedge \epsilon_2) \ .  \nn
\eeqa
Hence  we have the formula
\beqa
\text{Disc}(Q_m) =  \frac 9{ (\rho_1-\rho_2)^2 (\rho_2-\rho_3)^2 (\rho_3-\rho_1)^2}\; (\epsilon_1 \wedge \epsilon_2)^{\otimes^2}  \ . \nn
\eeqa

\noi
If $[\FF_{Q_m}:\FF]=1$ then $\rho_1,\rho_2,\rho_3 \in \FF$ and hence $D \in \FF$ which means that $s\big(\text{Disc}(Q_m)\big)=[1]$.
If  $[\FF_{Q_m}:\FF]=3$  then $\rho_1,\rho_2,\rho_3 \in \FF_{Q_m}\setminus \FF$ but $D$, being  invariant under the Galois group
$\mathbb Z_3$, is an element of $\FF$. Hence  $s\big(\text{Disc}(Q_m)\big)=[1]$.

\noi
If  $[\FF_{Q_m}:\FF]=2$  then  Gal$(\FF_{Q_m}/\FF)= \mathbb Z_2$ and the roots of $P(x)$ are $\rho_1 \in \FF, \rho_2 \in\FF_{Q_m}\setminus\FF$
and $\rho_3 = \bar \rho_2 \in \FF_{Q_m}\setminus\FF$. Hence 
\beqa
D= (\rho_1-\rho_2)(\rho_2 -\bar \rho_2) (\bar \rho_2 -\rho_1)  \nn
\eeqa
which is not  invariant under the Galois group. 
Thus $\FF(D)$ is a non-trivial quadratic extension of $\FF$ contained in $\FF_{Q_m}$ which means that
  $s\big(\text{Disc}(Q_m)\big)=[\FF_{Q_m}]$. Finally if $[\FF_{Q_m}:\FF]=6$, then  Gal$(\FF_{Q_m}/\FF)\cong S_3$
which acts on $\rho_1,\rho_2,\rho_3 \in \FF_{Q_m}\setminus\FF$ by permutation. In particular $D$ is not invariant under
odd permutations and hence $D\not \in \FF$. Thus $\FF(D)$ is a non-trivial quadratic extension of $\FF$ contained in 
$\FF_{Q_m}$. By the fundamental theorem of Galois theory there is only one such quadratic extension since the group
$ S_3$ has only one  subgroup of index two. Hence $s\big(\text{Disc}(Q_m)\big)$ must be  this quadratic extension.
\end{demo}

We now prove (3.1) $-$ (3.3) case by case. \\

(3.1): If $[\FF_{Q_m}:\FF]=1$ or $3$  it follows from Lemma \ref{lem:marcus} that Disc$(Q_m)\in L_{[1]}$ and  since 
Disc$(Q_{\lambda\text{Id}\cdot m}) =\lambda^2$Disc$(Q_m)$, it is clear that ${}_1\hat \nu_0$ is surjective.\\

(3.2): If $[\FF_{Q_m}:\FF]=2$   it follows from Lemma \ref{lem:marcus} that Disc$(Q_m)= [\FF_{Q_m}]$.
As in the proof of (3.1),
Disc$(Q_{\lambda\text{Id}\cdot m}) =\lambda^2$Disc$(Q_m)$ and
hence  all points of  the half-line $L_{[\FF_{Q_m}]}$ are in the image of ${}_2\nu_0(m)$.
Since  $\FF_{Q_m}$ can be an arbitrary non-trivial quadratic extension of $\FF$,  this implies that 
the image of ${}_2\nu_0$ is $(L^2\setminus\{0\})\setminus L_{[1]}$.\\

 (3.3): If $[\FF_{Q_m}:\FF]=6$   it follows from Lemma \ref{lem:marcus} that Disc$(Q_m)= [\tilde \FF_m]$ where 
$\tilde \FF_m$ is the unique quadratic extension of $\FF$ such that $\FF \subset \tilde \FF_m\subset \FF_{Q_m}$.
 As in the proof of (3.1),
Disc$(Q_{\lambda\text{Id}\cdot m}) =\lambda^2$Disc$(Q_m)$  and hence
all points in the half-line $L_{[\tilde \FF_m]}$ are in the image of 
 ${}_6\hat \nu_0$.  
By the very definition of $\mathbf \Delta_{6,2}$ this is enough to prove that 
${}_6 \nu_0$  is surjective.
\end{demo}

\begin{remark}
We do not know how to characterise quadratic extensions which are contained in Galois extensions of order six.
\end{remark}

\section{A linear embedding of  stable $SL(V)-$moduli space}\label{sec:compact}

In the two previous sections we gave a parametrisation of   the moduli spaces
of generic  two-dimensional commutative algebras
 under the action of $GL(V)$ and $SL(V)$.
In this section, we $\FF$ is algebraically closed, we will construct 
 an embedding  $\zeta_s$ of the $SL(V)-$moduli space of stable
commutative two-dimensional algebras onto a hypersurface $\Gamma_{\mathrm{Eisenstein}}$ in a four-dimensional vector
space. The equation defining this hypersurface is essentially the classical Eisenstein  identity for
the covariants of a  binary cubic.
Furthermore, the ``projectivisation'' of  $\zeta_s$ gives exactly the embedding of the stable $GL(V)-$moduli space
onto $\mathbb P^2$ of \cite{am}.

With this in mind we make the following definitions. Recall that $\M^c_{\text{st}}$ denotes the set of {\it all} stable
commutative algebra structures on $V$ and  $L=\Lambda^2(V^*)$.
\begin{definition}\label{def:zeta}  \hskip 1.truecm \phantom{marcus}\\
\begin{enumerate}
\item Let $\Gamma_{\mathrm{Eisenstein}} = \Big\{(A,B,D,C) \in L^2\times L^2 \times L^2 \times L^3 \setminus\{0\}\ \ \text{s.t.} \ \
27 D A^2 + 4 B^3 + 2^4 3^6  C^2=0 \Big\}$.

\item
Let $\zeta: \M^c_{\text{st}} \to L^2\times L^2 \times L^2 \times L^3$ be defined by
\beqa
\zeta(m)=(\tilde p_3(m),\tilde p_2(m),\mathrm{Disc}(Q_m),\mathrm{Inv}_m) \ , \ \ \forall m \in \M^c_{\text{st}} .\nn
\eeqa             
Note it follows from  Theorem \ref{theo:miracle}
and the  definition
of $ \M^c_{\text{st}}$ (cf. \eqref{eq:stable})
that $\mathrm{Im}(\zeta) \subset  \Gamma_{\mathrm{Eisenstein}}$. 
\item  The group $\FF^*$ acts on $L^2\times L^2 \times L^2 \times L^3$ by 
\beqa
\lambda\cdot(A,B,D,C)= \Big(\frac1 {\lambda^2} A,\frac1 {\lambda^2} B,\frac1 {\lambda^2} D,\frac1 {\lambda^3}C\Big) \ . \nn
\eeqa
We set
\beqa
\widetilde{{\mathbb P}}(L^2\times L^2 \times L^2 \times L^3) = (L^2\times L^2 \times L^2 \times L^3)/\FF^* \ . \nn
\eeqa
Note that $\Gamma_{\mathrm{Eisenstein}}$ is stable under this action and we will write 
$\widetilde{{\mathbb P}}(\Gamma_{\mathrm{Eisenstein}})=  \Gamma_{\mathrm{Eisenstein}}/\FF^*$. If
$U_3= \big\{[A,B,D,C] \in  \widetilde{{\mathbb P}}(\Gamma_{\mathrm{Eisenstein}}) \ \ \text{s.t.} \ \ D\ne 0 \big\}$,
we define the canonical affine chart $\psi_3: U_3\to \FF^2 $ and the canonical projective chart $\tilde \psi_3:
 \widetilde{{\mathbb P}}(\Gamma_{\mathrm{Eisenstein}})\to \mathbb P^2$ 
by
\beqa
\Psi_3([A,B,D,C])=\Big(\frac A D, \frac B D\Big) \ , \ \ \
\tilde \Psi_3([A,B,D,C])=[A, B,  D] \ .
\nn
\eeqa
\end{enumerate}
\end{definition}

\begin{remark}
Note that the action of $\FF^*$ on the set $\{(A,B,D,C) \in \Gamma_{\mathrm{Eisenstein}}$ s.t. $C\ne 0\}$  is principal.
It follows from the definition of $\zeta$ that $\zeta(m)$ is in this set iff $\nu(m) \not \in \Gamma_{\mathrm{Cardano}}$.
\end{remark}
\medskip
Consider the natural map $p: \M^c_{\text{st}}/SL(V) \to \M^c_{\text{st}}/GL(V)$.  There is a $GL(V)/SL(V) \cong\FF^\ast-$action on
$ \M^c_{\text{st}}/SL(V)$ given by  
\beqa
\label{eq:lam-action}
\lambda\cdot[m]_{SL(V)}= [ g\cdot m] \ \ \text{where} \ \ g \in GL(V) \ \ \mbox{is\ \ s.t.~} \det(g) = \lambda \ 
\eeqa
which satisfies
\beqa
p\Big(\lambda\cdot [m]_{SL(V)}\Big) = p\Big([m]_{SL(V)}\Big) \ . \nn
\eeqa

Since the coordinates of $\zeta$ are $GL(V)-$equivariant ``twisted'' polynomials, the map $\zeta$ factors to 
define a map $\zeta_s: \M^c_{\text{st}}/SL(V) \to \Gamma_{\mathrm{Eisenstein}}$.
It  follows  that $\zeta_s$ is $\FF^\ast-$equivariant with respect to the natural $\FF^\ast-$actions defined
above ({\it c.f.}  Definition  \ref{def:zeta} {\color{blue} (3)} and  Equation \ref{eq:lam-action})  and
factors to define a map  $\zeta_g: \M^c_{\text{st}}/GL(V) \to\widetilde{\mathbb P}(\Gamma_{\mathrm{Eisenstein}})$.
We have the commutative diagram


\beqa
 \xymatrix{
\M_3^c\ar[d]&\hskip -.5truecm\subset&\M^c_{\text{st}} \ar[dr]^{\hskip -.1truecm \zeta}\ar[d]&\\
\M_3^c/SL(V)\ar[d]&\hskip -.5truecm\subset&\M^c_{\text{st}}/SL(V) \ar[r]_{\zeta_s}\ar[d]_p&\Gamma_{\mathrm{Eisenstein}}\ar[d]^\pi \\
\M_3^c/GL(V)&\hskip -.5truecm\subset&\M^c_{\text{st}}/GL(V) \ar[r]_{\zeta_g}& \widetilde{{\mathbb P}}(\Gamma_{\mathrm{Eisenstein}}) 
&\hskip -.5truecm \supset& 
\hskip -.5truecm U_3\ar[r]_{\hskip -.4truecm \Psi_3}&\FF^2\\
}\nn
\eeqa
 and now  determine the isotropy of the two $\FF^\ast-$actions defined above.

\begin{proposition}
\phantom{marcus} \hskip 1.truecm \phantom{marcus}  
\begin{enumerate}
\item Let $(A,B,D,C) \in L^2\times L^2 \times L^2 \times L^3$ and
$I_{(A,B,D,C)} = \Big\{ \lambda \in \FF^\ast \ \ \mathrm{s.t.} \  \ \lambda\cdot(A,B,D,C) = (A,B,D,C) \Big\}$.
Then
$$
I_{(A,B,D,C)} = \left\{
\begin{array}{lll}
\{1\}&\mathrm{if}&C \ne 0 \\
\{1,-1\} &\mathrm{if}&C = 0  \ .
\end{array}\right.
$$
\item Let $[m]_{SL(V)} \in \M^c_{\text{st}}/SL(V)$ and
$I_{ [m]_{SL(V)}}= \Big\{ \lambda \in \FF^\ast \ \mathrm{s.t.} \ \lambda\cdot [m]_{SL(V)} = [m]_{SL(V)} \Big\}$.
Then
$$
I_{ [m]_{SL(V)}}  \left\{
\begin{array}{llll}
=&\{1\}&\mathrm{if}&\mathrm{Inv}_m \ne 0 \\
\subset&\mathbb \{1,-1\} &\mathrm{if} & \mathrm{Inv}_m = 0\ .
\end{array}\right.
$$
\end{enumerate}
\end{proposition}

\begin{demo}
(1) is immediate and (2) follows from the fact that $I_{ [m]_{SL(V)}} \subset I_{\zeta_s([m])}$ since $\zeta_s$ is
$\FF^\ast-$equivariant.
\end{demo}

The following theorem shows that many of the solutions 
of the Eisenstein equation are in the image of $\zeta$.
In fact, as we will show in Corollary  \ref{cor:surj-gene}, it turns out that 
all solutions are in the image of $\zeta$.
\begin{theorem}\label{theo:zeta-psurj}
The image of $\zeta_s :  \M^c_3/SL(V) \to \Gamma_{\mathrm{Eisenstein}}$ is given by
\beqa
\mathrm{Im}(\zeta_s)_{|_ {\M^c_3/SL(V)}} = \Gamma^\ast_{\mathrm{Eisenstein}}= \Big\{(A,B,D,C)  \in \Gamma_{\mathrm{Eisenstein}}  \ \ \mathrm{s.t.}  \ \
D\ne 0  \Big\} \nn \ .
\eeqa
\end{theorem}

\begin{demo}

We split  the proof into two lemmas corresponding to whether  $C=0$ or $C\ne 0$.

\begin{lemma}
Let $(A,B,D,0) \in \Gamma^*_{\mathrm{Eisenstein}}$.  Then there exists $m \in \M^c_3$ such that
$\zeta_s([m]_{SL(V)})=(A,B,D,0)$.
\end{lemma}
\begin{demo}
By  Theorem \ref{theo:slmodule} {\color{blue} (2)}, {\color{blue} (3.1)} and {\color{blue} (3.2)}  there exists $m \in \M^c_3$ such that
\beqa
\frac{\tilde p_3}{\text{Disc}(Q_m)} = \frac A D \ , \ \ 
\frac{\tilde p_2}{\text{Disc}(Q_m)} = \frac B D \ , \ \
\text{Disc}(Q_m) = D  \ . \nn
\eeqa
It follows that $\zeta_s([m]_{SL(V)}) = (\tilde p_3,\tilde p_2 , \text{Disc}(Q_m), 0) = (A, B ,D,0)$.
\end{demo}
\begin{lemma}
Let $(A,B,D,C) \in \Gamma^*_{\mathrm{Eisenstein}}$  be  s.t. $C\ne 0$.  Then there exists $m \in \M^c_3$ such that
$\zeta_s([m]_{SL(V)})=(A,B,D,C)$.
\end{lemma}
\begin{demo}
By  Theorem \ref{theo:slmodule} {\color{blue} (1)} the exists $m \in \M^c_3$ such that
\beqa
\frac{\tilde p_3}{\text{Disc}(Q_m)} = \frac A D \ , \ \ 
\frac{\tilde p_2}{\text{Disc}(Q_m)} = \frac B D \ , \ \
\frac{\text{Inv}_m} {\text{Disc}(Q_m)} = \frac C D  \ . \nn
\eeqa
By Theorem  \ref{theo:miracle},
\beqa
27 \text{Disc}(Q_m) \tilde p_3^2 + 4 \tilde p_2^3 = -2^4 3^6 \text{Inv}_m^2  \nn
\eeqa
and hence 
\beqa
27 \text{Disc}(Q_m) \Big(\frac A D  \text{Disc}(Q_m)\Big)^2 + 4 \Big(\frac B D  \text{Disc}(Q_m)\Big)^3 =  -2^4 3^6 \Big( \frac C D
 \text{Disc}(Q_m)\Big)^2\ . \nn
\eeqa
Dividing by Disc$(Q_m)^2$ this reduces to
\beqa
\text{Disc}(Q_m)\Bigg( 27\Big(\frac A D  \Big)^2 + 4 \Big(\frac B D  \Big)^3\Bigg) =  -2^4 3^6  \Big( \frac C D \nn
\Big)^2  \
\eeqa
but since by hypothesis 
\beqa
27 D A^2 +4 B^3 = -2^4 3^6 C^2 \  \nn
\eeqa
and $C \ne 0$, we get
\beqa
\text{Disc}(Q_m) = D \ . \nn
\eeqa
It follows that $\zeta_s([m]_{SL(V)}) = (\tilde p_3,\tilde p_2 , \text{Disc}(Q_m), \text{Inv}_m) = (A, B ,D,C)$.
\end{demo}

\end{demo}

\begin{corollary}\label{cor:surj-gene}
The map $\zeta_s :  \M^3_{\text{st}}/SL(V) \to \Gamma_{\mathrm{Eisenstein}}$ is  surjective.
\end{corollary}
\begin{demo}
By Theorem \ref{theo:zeta-psurj} it remains to show that the image of $\zeta_s$ contains
\beqa
\Gamma_{\mathrm{Eisenstein}}\setminus \Gamma^\ast_{\mathrm{Eisenstein}}= \Big\{ (A,B,0,C) \in \Gamma_{\mathrm{Eisenstein}} \Big\} \ . \nn
\eeqa
First note that $\zeta(m)=(A,B,0,C)$ implies that the fundamental cubic $Q_m$ has a multiple root and hence $m$ admits
at least one $v \in V$ such that $m(v,v)$ is proportional to $v$. Motivated by this observation, let us consider
$m\in\M^c_{\text{st}}$ of the form
\beqa
m(e_1,e_1)= a e_1 \ , \ \ m(e_2,e_2) = d e_2 \ , \ \ m(m_1,e_2) = e e_1 + f e_2 \ , \nn
\eeqa
where $a,d,e,f \in \FF$.
Using \eqref{eq:decom} the condition Disc$(Q_m)=0$ is equivalent to $(d-2e)^2(a-2f)^2=0$. Choosing $d=2$ and $e=1$ leads
to the multiplication table
\beqa
m(e_1,e_1)= a e_1 \ , \ \ m(e_2,e_2) = 2 e_2 \ , \ \ m(m_1,e_2) =  e_1 + f e_2 \ , \nn
\eeqa
and we now show that if $(A,B,0,C) \in \Gamma_{\mathrm{Eisenstein}}$ satisfies $B  C \ne 0$ then there exists $m$ of
this form  such that $\zeta(m)=(A,B,0,C)$.
From \eqref{eq:pt3}, \eqref{eq:pt2} and \eqref{eq:Invd} the equation $\zeta(m)=(A,B,0,C)$ is equivalent to
\beqa
\tilde p_3 &=& 72 (a+f)(a-2f)(\epsilon_1\wedge \epsilon_2)^2 = A \ , \nn \\ 
\tilde p_2 &=& - 36 (a-2f)^2 (\epsilon_1\wedge \epsilon_2)^2= B \ , \nn\\\
\text{Inv}_m&=& 4 (a-2f)^3 (\epsilon_1\wedge \epsilon_2)^3= C \ , \nn
\eeqa
where $B$ and $C$ satisfy the Eisenstein constraint
\beqa
 4 (-B)^3 =  (4 \times 27)^2 C^2 \ . \nn
\eeqa
One easily checks that if $(a,f) \in\FF^2$  is any solution of the linear system
\beqa
(a-2f)\epsilon_1 \wedge \epsilon_2 &=& - 9 \frac C B  \ , \nn \\
(a+ f) \epsilon_1 \wedge \epsilon_2 &=&-\frac{AB}{9\times 72 C} \ , \nn
\eeqa
the corresponding multiplication $m$ satisfies $\zeta(m)=(A,B,0,C)$.

To complete the proof of the corollary, it only  remains to show that the image of $\zeta_s$ contains
 points of $\Gamma_{\mathrm{Eisenstein}}$ which are of the form: $(A,0,0,0)$.
(Recall that $B=0$ is equivalent to $C=0$ by the Eisenstein condition.)
For this, consider $m \in \M^c_{\text{st}}$
of the form
\beqa
m(e_1,e_1)=  e_1 \ , \ \ m(e_2,e_2) = c e_1  \ , \ \ m(e_1,e_2) = \frac12 e_2 \ , \nn
\eeqa
$c \in \FF$.
Using \eqref{eq:pt3}, \eqref{eq:pt2} and \eqref{eq:Invd} we obtain
\beqa
\tilde p_3 = 27 c (\epsilon_1 \wedge\epsilon_2)^2 \ , \ \
\tilde p_2 =  \text{Disc}(Q_m) = 0 \ , \text{Inv}_m =0 , \nn
\eeqa
and hence taking
\beqa
c= \frac A {27} \nn \ ,
\eeqa
the corresponding multiplication $m$ satisfies  $\zeta(m)=(A,0,0,0)$.
\end{demo}

The following theorem is the main theorem of the paper.Assuming $\FF$ is algebraically closed  we 
show that $\zeta_s$ is an embedding of the $SL(V)-$moduli space of stable
commutative two-dimensional algebras onto the hypersurface $\Gamma_{\mathrm{Eisenstein}}$ in a four-dimensional vector
space. Furthermore, the ``projectivisation'' of  $\zeta_s$ gives exactly the embedding of the $GL(V)-$moduli space
onto $\mathbb P^2$ of \cite{am}.

The theorem also shows that moduli spaces $-$ both for the action
of $GL(V)$ and $SL(V)$$-$ 
of {\it stable} two-dimensional commutative algebras can be thought of 
as  compactifications of  moduli spaces of {\it generic} two-dimensional commutative algebras.
For instance $\zeta_s(\M^c_{\text{st}}/SL(V))= \Gamma_{\mathrm{Eisenstein}}$ is  a ``compactification'' of
$\zeta_s(\M^c_3/SL(V))$
 where ``the points at infinity'' are given by
$ \Big\{(A,B,0,C) \in 
\Gamma_{\mathrm{Eisenstein}}\Big\}$.

\vfill\eject

\begin{theorem}\label{theo:c-cube}
Suppose $\FF$ is algebraically closed and consider the commutative diagram:

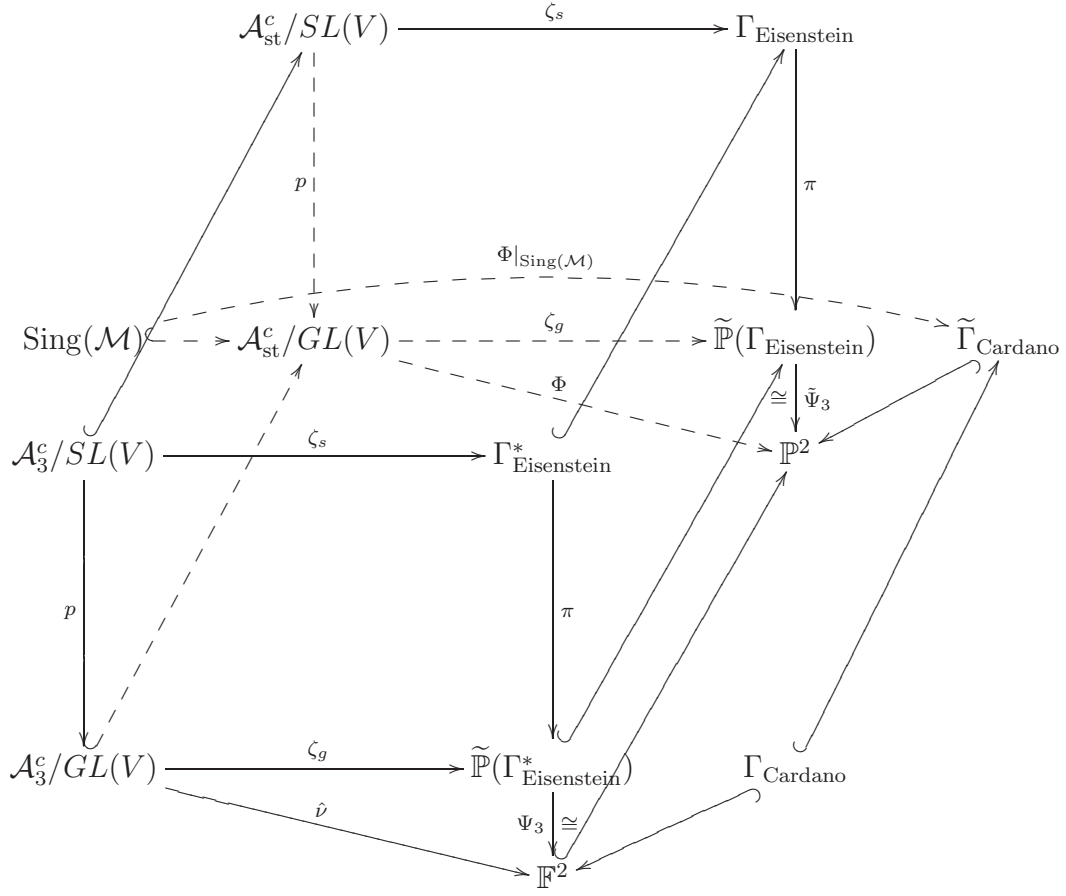
\beqa
\xymatrix{
&\M^c_{\mathrm{st}}/SL(V)\ar[rr]^{\zeta_s}\ar@{-->}[ddd]_p  &&\Gamma_{\mathrm{Eisenstein}}\ar[ddd]^\pi \\{\color{white} toto}\\
{\color{white} toto}\\
\mathrm{Sing}({\cal M})\ar@{^{(}-->}[r]\ar@{-->}@/^2pc/[rrrr]^{\Phi\vert_{\mathrm{Sing}({\cal M})}}&\M^c_{\mathrm{st}}/GL(V)\ar@{-->}[rr]^{\zeta_g} \ar@{-->}[rrd]^{\Phi}  &&\widetilde{\mathbb P}(\Gamma_{\mathrm{Eisenstein}})\ar[d]_{\cong}^{\tilde \Psi_3}
&\widetilde{\Gamma}_{\mathrm{Cardano}}  \ar@{^{(}->}[ld]\\
\M^c_3/SL(V)\ar[rr]^{\zeta_s}\ar@{^{(}->}[uuuur] \ar[ddd]_p&&\Gamma^*_{\mathrm{Eisenstein}} \ar[ddd]^\pi\ar@{^{(}->}[uuuur]&
\mathbb P^2  
\\ {\color{white} toto}\\
{\color{white} toto}\\
\M^c_3/GL(V)\ar[rr]^{\zeta_g}\ar@{^{(}-->}[uuuur]\ar[rrd]^{\hat \nu} && \widetilde{\mathbb P}(\Gamma^*_{\mathrm{Eisenstein}}) \ar[d]_{\Psi_3}^{ \cong}\ar@{^{(}->}[uuuur]  &\Gamma_{\mathrm{Cardano}}\ar@{^{(}->}[ld]\ar@{^{(}->}[uuuur]\\
&& \mathbb F^2\ar@{^{(}->}[uuuur]
} \nn
\eeqa

Then 
\begin{enumerate}
\item $\zeta_s: \M^c_3/SL(V) \to \Gamma^*_{\text{Eisenstein}}$ is a bijection;
\item $\zeta_g : \M^c_3/GL(V) \to \mathbb P(\Gamma^*_{\text{Eisenstein}})$ is a bijection;
\item $\zeta_s: \M^c_{\text{st}}/SL(V) \to \Gamma_{\text{Eisenstein}}$ is a bijection;
\item $\zeta_g : \M^c_{\text{st}}/GL(V) \to  \widetilde{\mathbb P}(\Gamma_{\text{Eisenstein}})$ is a bijection;
\item $\Psi_3:  \widetilde{\mathbb P} (\Gamma^*_{\text{Eisenstein}})\to \FF^2$ is a bijection;
\item $\tilde \Psi_3: \widetilde{\mathbb P} (\Gamma_{\text{Eisenstein}})\to {\mathbb P}^2$ is a bijection;
\end{enumerate}

\end{theorem}
\begin{demo}
Let us remark that the front face of the diagram
involves only generic commutative algebras and the back face only stable commutative algebras.\\

We first prove part 1 and part 2.
Surjectivity for both maps was proved  in Theorem \ref{theo:zeta-psurj}
over an arbitrary field.

To prove the injectivity of $\zeta_s$
let $[m]_{SL(V)}, [m']_{SL(V)} \in \M^c_3/SL(V)$ be such that $\zeta_s([m]_{SL(V)})= \zeta_s([m']_{SL(V)})$. Then by definition,
\beqa
\label{eq:ccond}
(\tilde p_3(m), \tilde p_2(m), \text{Disc}(Q_m), \text{Inv}_m) = (\tilde p_3(m'), \tilde p_2(m'), \text{Disc}(Q_{m'}), \text{Inv}_{m'})
\eeqa
and since $\FF$ has no non-trivial extensions, 
 it follows  from Theorem \ref{theo:slmodule} {\color{blue} (1)},  {\color{blue} (2)} and  {\color{blue} (3)} that 
$[m]_{SL(V)}= [m']_{SL(V)}$.

To prove injectivity of $\zeta_g$, first observe that by definition, $ \widetilde{\mathbb P}(\Gamma^\ast_{\mathrm{Eisenstein}})= U_3$
and $\hat \nu=\Psi_3 \circ \zeta_g$. Hence if $\zeta_g([m]_{GL(V)})= \zeta_g([m']_{GL(V)})$ then $\hat\nu([m]) = \hat\nu([m')]$ and, since
$\FF$ has no non-trivial extensions, it follows from Theorems \ref{theo:inv-gen1-1}, \ref{theo:inv-non-gener} and
\ref{theo:inv-non-gener2} that $[m]_{GL(V)}= [m']_{GL(V)}$.\\

We now prove part 3 and part 4. 
Surjectivity for both maps was proved in Corollary \ref{cor:surj-gene}
over an arbitrary field.

To prove the injectivity of $\zeta_g: \M^c_{\text{st}}/GL(V) \to \widetilde{{\mathbb P}}(\Gamma_{\text{Eisenstein}})$  recall that,
as  already observed in \eqref{eq:id2}, we have 
\beqa
\Phi =\tilde \Psi_3 \circ \zeta_g  \ , \nn
\eeqa
where $\Phi :  \M^c_{\text{st}}/GL(V) \to \mathbb P^2$
is  given by ({\it cf} \cite{am})
\beqa
\Phi([m]) = [\tilde p_3(m),\tilde p_2(m), \text{Disc}(Q_m)]  \ .\nn
\eeqa
It was shown in  \cite{am} that $\Phi$ is injective and hence it follows that both $\zeta_g$ and $\tilde \Psi_3$ are
injective (as we have already shown that $\zeta_g$ is surjective). This proves part (3) of the theorem.

To prove the injectivity of $\zeta_s: \M^c_{\text{st}}/SL(V) \to \Gamma_{\text{Eisenstein}}$ let $m,m' \in \M^c_{\text{st}}$ be
 such that
 \beqa
 \zeta_s([m]_{SL(V)}) = \zeta_s([m']_{SL(V)})\ .\nn
 \eeqa
 Then
 \beqa
 \zeta_g([m]_{GL(V)}) = \zeta_g([m']_{GL(V)})\ \nn
 \eeqa
 and so $\Phi([m]) = \Phi([m'])$. By \cite{am} it follows that there exists $g \in GL(V)$ such that $m ' = g\cdot m$.
 On the other hand $\zeta_s([m]_{SL(V)} =\zeta_s([m']_{SL(V)}$ is equivalent to
 \beqa
 \label{eq:cond}
(\tilde p_2(m),\tilde p_3(m), \text{Disc}(Q_m), \text{Inv}_m) = (\tilde p_2(m'),\tilde p_3(m'), \text{Disc}(Q_{m'}), \text{Inv}_{m'})
 \ .
\eeqa
If Disc$(Q_m)=$ Disc$(Q_{m'})\ne 0$ it follows from Theorem \ref{theo:slmodule} {\color{blue} (1)},  {\color{blue} (2)} and  {\color{blue} (3)} that
$[m]_{SL(V)} = [m']_{SL(V)}$.\\

\noi
If Disc$(Q_m)= $ Disc$(Q_{m'})=0$ 
it follows  (see Appendix  \ref{sec:non-gen})  that $m$ and $m'$ are $GL(V)-$equivalent to one of the algebras in Table 
\ref{tab:1} (two types) or Table \ref{tab:2} (one type).
The condition \eqref{eq:cond} implies that $m$ and $m'$ must be of the same type.  In fact, as we show in the following
lemma, for such algebras $GL(V)-$equivalence implies $SL(V)-$equivalence.

\begin{lemma}
 Let $(e_1,e_2)$ and $(e_1',e_2')$ be bases of $V$. We denote by $(\epsilon_1,\epsilon_2)$ and $(\epsilon_1',\epsilon_2')$
 the corresponding dual bases. Define $\alpha \in \FF^*$ by
$
e_1'\wedge e'_e = \alpha e_1 \wedge e_2$. 
\begin{enumerate}
\item[(i)]  Let $\nu,\nu'\in \FF\setminus\{\frac12\}$ and consider the two commutative algebras $m,m'$ given by
\beqa
\begin{array}{lll}
m(e_1,e_1)= e_1\ , &m(e_2,e_2) =e_2\ ,& m(e_1,e_2)= \frac12 e_1 + \nu e_2\ , \\
m'(e'_1,e'_1)= e'_1\ ,& m(e'_2,e'_2) =e'_2\ ,& m(e'_1,e'_2)= \frac12 e'_1 + \nu' e'_2 \ .
\end{array} \nn
\eeqa
If $m$ and $m'$ are $GL(V)-$equivalent and satisfy \eqref{eq:ccond} then they are $SL(V)-$equivalent.
\item[(ii)] Consider  the two commutative algebras $m,m'$ given by
\beqa
\begin{array}{lll}
m(e_1,e_1)= 0\ , &m(e_2,e_2) =e_2\ ,& m(e_1,e_2)= \frac12 e_1 + e_2 , \\
m'(e'_1,e'_1)= 0\ , &m(e'_2,e'_2) =e'_2\ ,& m(e'_1,e'_2)= \frac12 e'_1 + e'_2 ,
\end{array}\nn
\eeqa
If $m$ and $m'$ are $GL(V)-$equivalent and satisfy \eqref{eq:ccond} then they are $SL(V)-$equivalent.
\item[(iii)]  Let $\lambda,\lambda'\in \FF^*$ and consider the two commutative algebras $m,m'$ given by
 \beqa
\begin{array}{lll}
m(e_1,e_1)= e_1\ , &m(e_2,e_2) =\lambda e_1\ ,& m(e_1,e_2)= \frac12  e_2 , \\
m'(e'_1,e'_1)= e'_1\ , &m(e'_2,e'_2) =\lambda e'_1\ ,& m(e'_1,e'_2)= \frac12  e'_2 ,
\end{array}\nn
\eeqa
If $m$ and $m'$ are $GL(V)-$equivalent and satisfy \eqref{eq:ccond} then they are $SL(V)-$equivalent.
\end{enumerate}
\end{lemma}
\begin{demo}
(i) : Since $m$ and $m'$ are $GL(V)-$equivalent we must have $\nu =\nu'$ (see Theorem \ref{theo:tutu}).
Since by \eqref{eq:ccond}  we have
\beqa
\begin{array}{ll}
\tilde p_3(m) =-18 (2\nu-1)(\nu+1) (\epsilon_1 \wedge \epsilon_2)^2&= \tilde p_3(m')=
-18 (2\nu-1)(\nu+1) (\epsilon'_1 \wedge \epsilon'_2)^2   \nn\\
\tilde p_2(m) =-9 (2\nu-1)^2 (\epsilon_1 \wedge \epsilon_2)^2&= \tilde p_2(m')=
-9 (2\nu-1)^2 (\epsilon'_1 \wedge \epsilon'_2)^2  \ .
\end{array}\nn 
\eeqa
Substituting $ \epsilon'_1 \wedge \epsilon'_2=1/\alpha \epsilon_1 \wedge \epsilon$ and solving the two equations above
gives $\alpha^2=1$.
Since Inv${}_m= -1/2(2\nu-1)^3 (\epsilon_1\wedge\epsilon_2)^3 = $  Inv${}'_m= -1/2(2\nu-1)^3 (\epsilon'_1\wedge\epsilon'_2)^3 $
then $\alpha^3 = 1$. This implies that $\alpha=1$ and 
 the linear mapping $g: V\to V$ given by $g(e_1)=e'_1, \  g(e_2)=e'_2$ is an algebra isomorphism of
determinant one.

(ii): Just as above the equality \eqref{eq:ccond} implies that $\alpha^3=\alpha^2=1$ and hence
the linear mapping
$g: V\to V$ given by $g(e_1)=e'_1, \  g(e_2)=e'_2$ is an algebra isomorphism of
determinant one.

(iii):
Since $\tilde p_3(m) =-27 \lambda (\epsilon_1 \wedge \epsilon_2)^2 =\tilde p_3(m')=-
27 \lambda' (\epsilon'_1 \wedge \epsilon'_2)^2$, we have $\lambda'=\alpha^2 \lambda$. Let $\mu \in \FF^*$ be
the unique square root of $\lambda/\lambda'$ such that $\mu \alpha=1$. It is straightforward to check that
the linear map $g: V\to V$ given
\beqa
g(e_1) =e_1'\ , g(e_2) = \mu e_2' \ , \nn
\eeqa
is an algebra isomorphism of determinant one.
\end{demo}
This completes the proof of part (3) of the Theorem.\\

Finally, observe that part (6) of the theorem implies part (5) of the theorem and that part (6) is a consequence
of the fact that 
\beqa
\tilde \Psi_3 = \Phi \circ \zeta_g^{-1} \  \nn
\eeqa
is the composition of the bijections $\zeta_g^{-1}$ ({\it cf} part (4) above) and $\Phi$ ({\it cf} \cite{am}).

\end{demo}

\section{Disc$(D)$-generic algebras  modulo $SL(V)\times SL(V)$ and binary quadratic forms modulo $SL(V)$}\label{sec:bhar}
Recall that if $m\in {\cal A}^c$,  by definition Disc$(D_m)$ is the discriminant of the
quadratic form $\det L_m$ ({\it cf} Definition \ref{def:D-bha}).
In this section we consider the set of  {Disc$(D)$-generic algebras
\beqa
\M^c_D = \Big\{ m \in \M^c  \ \  \text{s.t.} \ \  \text{Disc}(D_m) \ne 0 \Big\} \ , \nn
\eeqa
modulo the left action of $SL(V) \times SL(V)$ given by
\beqa
{}^{(g_1,g_2)}m(u,v)= g_1 m(g_2^{-1} u, g_2^{-1} v) \ . \nn
\eeqa 
Restricting this action to the diagonally embedded $SL(V)$
\beqa
g \mapsto (g,g)  \nn
\eeqa
we recover the action of  $SL(V)$ of \eqref{eq:GLV}.
It turns out that Disc$(D)$ is not just  $SL(V)$ invariant but in fact  $SL(V) \times SL(V)$ invariant.
We first show that the map which associates to $m$ in $\M^c_D$ the binary quadratic form $D_m$ induces a bijection $\hat D$ 
of the moduli space $\M^c_D/(SL(V)\times SL(V))$ with the moduli space of non-degenerate 
quadratic forms  $S^2(V^*)_{\text{n.d.}}/SL(V)$.
As an application we define a group structure on the elements of $\M^c_D/SL(V)\times SL(V)$ of fixed discriminant $\Delta$,
and show that $\hat D$ induces a group isomorphism between  this group and the  group  of
$SL(V)-$orbits of quadratic forms of discriminant $\Delta$ (with  Gauss composition as the group law).

\subsection{Isomorphism of moduli spaces}

We first show that Disc$(D_m)$ is an $SL(V) \times SL(V)$ invariant.
\begin{proposition}\label{prop:qD}
Let $m \in \M^c$ and $(g_1,g_2) \in SL(V)\times SL(V)$. Then, 
\begin{enumerate}
\item [(i)]
$
D_{{}^{(g_1,g_2)} m} = g_2 \cdot D_m \ .
$
\item[(ii)]
$
\mathrm{Disc} (D_{{}^{(g_1,g_2)} m}) = \mathrm{Disc}(D_m) \ . 
$
\end{enumerate}
\end{proposition}

\begin{demo}
We set $\tilde m= {}^{(g_1,g_2)} m$.
By definition for all $u,v$ in $V$,
\beqa
\tilde m (u,v) = g_1 m (g_2^{-1} u, g_2^{-1} v)   \nn
\eeqa
and thus
\beqa
L_u^{\tilde m}(v) = g_1 L^m_{g_2^{-1}u}(g_2^{-1} v) \ . \nn
\eeqa
Hence
\beqa
L^{\tilde m}_u = g_1\circ L^m_{g_2^{-1}u} \circ g_2^{-1} \ , \nn
\eeqa
and taking determinants 
\beqa
\det(L^{\tilde m}_u) = \det ( L^m_{g_2^{-1}u}) \ .\nn
\eeqa
This implies (see \eqref{eq:DT}) 
\beqa
D_{\tilde m}(u) = D_m(g_2^{-1} u) = g_2 \cdot D_m(u) \ , \ \ \forall u \in V \  . \nn
\eeqa
This proves part (i). Part (ii) follows from this and the fact that,  
as  is well known,  two binary quadratic forms in the same $SL(V)-$orbit have the same discriminant.
\end{demo}
The following corollary is a simple consequence of the proof of the proposition and the fact that $D_m$ is non-degenerate
if $m \in \M^c_D$.
\begin{corollary}\label{cor:xi}
The map $D: \M^c_D \to S^2(V^*)_{\text{n.d.}}$  factors to a  map $\hat D:$ $\M^c_D/SL(V)\times SL(V)\to S^2(V^*)_{\text{n.d.}}/SL(V)$.
If $m \in \M^c_D$ we write $\hat D_{[m]} = [ D_m]$ for the image of the equivalence class  $[m]$ by $\hat D$.
\end{corollary}


\begin{proposition}
\label{prop:hatD}
The map $\hat D: \M^c_D/SL(V)\times SL(V) \to S^2(V^*)_{\text{n.d.}}/SL(V)$ 
 defined in  Corollary \ref{cor:xi} is a bijection.
\end{proposition}

\begin{demo}
To show that $\hat D: \M^c_D/SL(V)\times SL(V) \to S^2(V^*)_{\text{n.d.}}/SL(V)$ is surjective 
it is sufficient to show that $D: \M^c \to S^2(V^*)$ is surjective.
 Let  $(e_1,e_2)$ be a basis of $V$ with dual  basis  $(\varepsilon_1,\varepsilon_2)$ and $m: V\times V \to V$ be given by:
\beqa
m(e_1,e_2) = a e_1 + b e_2\ , \ \
m(e_2,e_2)=  c e_1 + d e_2 \ , \ \
m(e_1,e_2) = e e_1 + f e_2 \ , \ \ 
m(e_2,e_1) = e e_1 + f e_2 \ , \nn
 \eeqa
then by \eqref{eq:TD} 
\beqa
D_m(x e_1 + y e_2) =(af-be) x^2 + (ad-bc) xy + (ed-cf) y^2 \ .\nn
\eeqa
If  $q= \alpha x^2 + 2 \beta xy + \gamma y^2$ then $D_m = q$ {\it iff}
\beqa
\label{eq:abcdef}
af-be =\alpha \ ; \ \ 
 ad-bc= 2 \beta \ ; \ \ 
 ed-cf =\gamma \ .
\eeqa
If $q=0$ then $m \equiv 0$ satisfies $\hat D(m) =q$. If $q\neq 0$ then without loss of generality we can suppose that
$\beta \neq 0$. Then the vectors $v_1 = a e_1 + b e_2$ and $v_2 = c e_1 + d e_2$ form a basis of $V$ and there exist 
constants $A,B$  s.t.
\beqa
e e_1 + f e_2 = A v_1 + B v_2 \ . \nn
\eeqa
Substituting in the system above leads to 
\beqa
A = \frac{\gamma}{2 \beta}  \ \ ; \ \ B= \frac{\alpha}{2 \beta} \ \ ;  \ \ ad -bc =2 \beta \ .\nn 
\eeqa
In other words our system has a solution {\it iff} $ad -bc =2 \beta$ has a solution but this is obvious 
({\it e.g.} $a=2,d=\beta, b=c=0, e = \frac\gamma \beta, f=\frac \alpha 2 $). This completes the proof of the surjectivity of $\hat D$.

To show that $\hat D$ is injective we first prove the following lemma:
\begin{lemma}\label{lem:qD-inj}
Let $q$ be in $S^2(V^*)_{\text{n.d.}}$ and suppose that  $D(m_1) = D(m_2)=q$ for $m_1,m_2$ in $ \M^c_D$. Then there exists
$(g,1) \in SL(V) \times SL(V)$ such that ${}^{(g,1)} m_1 = m_2$.
\end{lemma}

\begin{demo}
It follows from \eqref{eq:abcdef} that  $D_m =q$ {\it iff}
\beqa
af-be=\alpha \ ,  \ \
ad-bc=2 \beta \ , \ \ 
ed-cf=\gamma \ \ . \nn
\eeqa
In particular, the multiplication $m_0$ given by
\beqa
a_0=2 \ , \ \ d_0=\beta \ \ , b_0 = c_0 = 0   \ , e_0 = \frac\gamma \beta \ , \ \ f_0=\frac \alpha 2 \nn
\eeqa
satisfies $D_{m_0}=q$. It is easy to check that if we set
\beqa
g= \begin{pmatrix}
\frac a 2 &\frac c \beta\\\frac b 2 &\frac d \beta
 \end{pmatrix} \ , \nn
\eeqa
then  $\det(g) =1$ and
\beqa
{}^{(g,1)} m_0 = m  \ . \nn
\eeqa
Thus the group $SL(V)\times\{ 1 \} $ acts transitively on the set of solutions of $D_m=q$ which proves the lemma.
\end{demo}

We now prove $\hat D$ is injective.  Suppose that 
\beqa
\hat D_{[m_1]} = \hat D_{[m_2]} \ .\nn
\eeqa
Thus by definition
\beqa
[D_{m_1}] = [D_{m_2}] \ , \nn
\eeqa
from which it follows that there exist $h$ in $SL(V)$ s.t.
\beqa
h \cdot D_{m_1} = D_{m_2} \ . \nn
\eeqa
 By Proposition \ref{prop:qD}{\color{blue}(i)}  this means that 
\beqa
D_{{}^{(1,h)}m_1} = D_{m_2} \ , \nn
\eeqa
and therefore by Lemma \ref{lem:qD-inj} there exists $g$ in $SL(V)$ such that
\beqa
{}^{{}^{(g,1)}}({}^{(1,h)} m_1)= m_2 \ . \nn 
\eeqa 
Hence
\beqa
{}^{(g,h)} m_1 = m_2 \ ,\nn
\eeqa
in other words
\beqa
[m_1] = [m_2] \ \nn
\eeqa
and  $\hat D$ is injective.
\end{demo}

\subsection{Group structure}
The Gauss composition law is a non-trivial group law on 
the set of $SL(V)-$orbits of non-degenerate quadratic forms with fixed  non-zero discriminant.
There is therefore a unique group structure on the set of $SL(V)\times SL(V)-$orbits 
of  Disc$(D)-$generic algebras with fixed discriminant, which makes  $\hat D$ ({\it cf} Proposition \ref{prop:hatD}) into a group
isomorphism. 
We will now give an intrinsic description of this group structure.

Let now $[m],[m']$ be two elements of ${\cal A}^c_D/SL(V)\times SL(V)$ such that Disc$(\hat D_{[m]}) =$Disc$(\hat D_{[m']})$
is not a square in $L^2$. We can always find a basis $\{e_1,e_2\}$ of $V$ which
is neither orthogonal  for $\hat D_{[m]}$ nor for $\hat D_{[m']}$. In this basis we write
\beqa
\hat D_{[m]}(xe_1 + y e_2) &=& \alpha^{\phantom{'}} x^2 +2 \beta^{\phantom{'}} xy + \gamma^{\phantom{'}} y^2\ , \nn\\
\hat D_{[m']}(xe_1 + y e_2) &=& \alpha' x^2 +2 \beta' xy + \gamma' y^2 \ , \nn\\
\text{Disc}(\hat D_{[m]}) = \text{Disc}(\hat D_{[m']})&=& \Delta (\epsilon_1\wedge\epsilon_2)^{\otimes^2} \ ,\nn
\eeqa 
where  $\{\epsilon_1,\epsilon_2\}$ is the dual basis of $\{e_1,e_2\}$,  $\Delta$ is  an element of $\FF^*$ {\it which is not a square in $\FF^*$} and $\beta \beta'\ne 0$.
Note that $\Delta= 4 \beta^2 - 4 \alpha \gamma= 4 \beta'{}^2 -4 \alpha'\gamma'$ and hence 
it follows that $\alpha \gamma \alpha' \gamma'\ne 0$.
By Lemma \ref{lem:qD-inj} the multiplications $m_0$ and $m_0'$ given by
\beqa
&&a_0=2 \ , \ \ d^{\phantom{'}}_0=\beta^{\phantom{'}} \ \ , b^{\phantom{'}}_0 = c^{\phantom{'}}_0 = 0   
\ , e^{\phantom{'}}_0 = \frac{\gamma^{\phantom{'}}} \beta^{\phantom{'}} \ , \ \ f_0=\frac {\alpha^{\phantom{'}}} 2 \nn\\
&&a'_0=2 \ , \ \ d'_0=\beta' \ \ , b'_0 = c'_0 = 0   \ , e'_0 = \frac{\gamma'} {\beta'} \ , \ \ f'_0=\frac {\alpha'} 2 \nn
\eeqa
satisfy $\hat D_{[m_0]} = \hat D_{[m]}$ and $D_{[m'_0]} = \hat D_{[m']}$.
The discriminants of $m_0, m'_0$ are both equal to $\Delta (\epsilon_1\wedge \epsilon_2)^{\otimes^2}$
so we can write 
\beqa
f_0 = \frac{d_0^2 - \frac \Delta 4}{2e_0 d_0} \ , \ \ f'_0 = \frac{d'_0{}^2 - \frac \Delta 4}{2e'_0d'_0} \nn \ 
\eeqa
and this suggests that we define the ``product'' of $[m_0]$ and $[m'_0]$ to be the class of the multiplication $m''_0$ given
 in the basis $\{e_1,e_2\}$  by
\beqa
&&a''_0=\frac{a_0 a'_0}2  \ , \ \ d''_0=d_0 d_0' \ \ , b''_0 =b_0 b'_0\ , \ \ c''_0 = c_0 c'_0   \ , \ \ e''_0 =e_0 e'_0 \ , \ \ 
f''_0=\frac {(d_0 d'_0)^2 -\frac \Delta 4} {2e_0 d_0 e'_0 d'_0} \nn \ .
\eeqa 
By Lemma \ref{lem:qD-inj} one gets
\beqa
\hat D_{[m''_0]}(x e_1 + y e_2) =  \frac{(\beta\beta')^2 -\frac \Delta 4}{\gamma \gamma'} x^2 +2 \beta \beta' xy + \gamma\gamma'y^2 \nn
\eeqa
from which it follows that the discriminant of $m''_0$ is equal to $\Delta(\epsilon_1\wedge \epsilon_2)^{\otimes^2}$.

\begin{proposition}
\label{prop:group-m}
Let $Q_\Delta=S^2(V^*)_\Delta/SL(V)$ denote the set of $SL(V)-$equivalent classes of
quadratic forms of discriminant $\Delta$ and let $\star:Q_\Delta \times Q_\Delta \to Q_\Delta$ be Gauss composition. Then
$\hat D_{[m]}\star \hat D_{[m']}= \hat D_{[m'']}$.
\end{proposition}

\begin{demo}
Set 
\beqa
\FF^*_\Delta= \Big\{a^2 -\Delta b^2 \in \FF^*\  \ \  \text{s.t.} \ \ a,b \in \FF \Big\} \nn
\eeqa
and 
\beqa
G_\Delta = \FF^*/\FF^*_\Delta \ . \nn
\eeqa
It is well know that  $Q_\Delta$ with Gauss composition is isomorphic to $G_\Delta$ under the map $[q] \mapsto [q(v)]$ where
$q \in S^2(V^*)_\Delta$ and $q(v)$ is its  value at any {\it non-isotropic} $v \in V$.
Hence to prove the proposition it is sufficient to check that 
\beqa
\hat D_{[m'']}(e_2)= \hat D_{[m]}(e_2)\hat D_{[m']}(e_2) \ , \nn
\eeqa 
which is immediate since
\beqa
\hat D_{[m_0]}(e_2) = \gamma \ , \ \ 
\hat D_{[m'_0]}(e_2) = \gamma' \ , \ \ 
\hat D_{[m''_0]}(e_2) = \gamma'' =\gamma \gamma'\ . \ \ 
\eeqa
\end{demo}

\begin{remark}
The essential point of our construction was the choice of the basis $\{e_1,e_2\}$ with respect to which 
neither  $\hat D_{[m]}$ nor $\hat D_{[m']}$ was diagonal. Once this  choice has been made $m_0, m_0'$ and $m''_0$ are uniquely defined
but there are many bases of $V$ with the above property.  Our proposition shows that in fact the
$SL(V)\times SL(V)$ equivalence class of $m''$ is independent of this choice.
\end{remark}

\begin{remark}
Proposition \ref{prop:group-m} is an analogue of the isomorphism
\beqa
\mathrm{Cl}(\mathbb Z^2 \otimes \mathrm{Sym}^2 \mathbb Z^2;\Delta) \to \mathrm{Cl}(\mathbb(Sym^2 \mathbb Z^2)^*;\Delta) \ , \nn
\eeqa
of Bhargava \cite{bhar} (p. 225).

\end{remark}

\appendix
\section{Non-generic commutative algebras} \label{sec:non-gen}

\subsection{Stable commutative algebras}
As we pointed out in Section  \ref{sec:res} the $GL(V)-$moduli space  ${\cal M}^c_{\text{st}}$ of {\it stable} two-dimensional
commutative algebras is larger than  the $GL(V)-$moduli space of {\it generic}  two-dimensional
commutative algebras ${\cal M}^c_3$. In fact there is a stratification
\beqa
{\cal M}^c_{\text{st}} = {\cal M}^c_3 \cup  {\cal M}^c_2 \cup {\cal M}^c_1 \ , \nn
\eeqa
where
\beqa
{\cal M}^c_3 &=& \Big\{ m \in {\cal M} \ \ \text{s.t.  \ \ Disc}(Q_m) \ne 0  \Big\} \ ; \nn\\
{\cal M}^c_2 &=& \Big\{ m \in {\cal M} \ \ \text{s.t.  \ \ Disc}(Q_m) = 0 \ \ \text{and} \ \ \tilde p_2(m) \ne 0   \Big\} \ ; \nn\\
{\cal M}^c_1 &=& \Big\{ m \in {\cal M} \ \ \text{s.t.  \ \ Disc}(Q_m) = \tilde p_2(m) =0 \ \ \text{and} \ \ \tilde p_3(m) \ne 0   \Big\}
\ . \nn
\eeqa
By definition the open {\it stratum} is the moduli space of  generic algebras ${\cal M}_3^c$
and is of dimension two. The next {\it stratum}  ${\cal M}^c_2$ is of dimension one and corresponds to algebras whose fundamental
cubic has exactly two distinct roots, and the last {\it stratum}  ${\cal M}^c_1$ is a single
point  and corresponds to algebras whose fundamental  cubic has exactly one  root (see below).

In Section \ref{sec:moduli} we showed that for $[m] \in {\cal M}_3$:
\begin{itemize}
\item if $\hat\nu([m]) \not \in \Gamma_{\text{Cardano}}$ then $\hat\nu([m])$ completely determines the equivalence class of $m$
(even if $\FF$ is not algebraically closed);
\item if $\hat\nu([m])  \in \Gamma_{\text{Cardano}}\setminus\{(0,0)\}$ then $\hat\nu([m])$ together with the equivalence class of
the splitting field (at most quadratic) of the
fundamental cubic $Q_m$ completely determine the equivalence class of $m$;
\item if  $\hat\nu([m]) =(0,0)$  then $\hat\nu([m])$ together with the equivalence class of the splitting field of the
fundamental cubic $Q_m$ completely determine the equivalence class of $m$; 
\item  $\hat\nu([m])  \in \Gamma_{\text{Cardano}} \setminus\{(0,0)\}$ \ {\it iff} \  $[m] \in$
$\Big($Sing$({\cal M}^c_{\text{st}}) \setminus$  Sing(Sing$({\cal M}^c_{\text{st}}))\Big)\cap {\cal M}^c_3$;
\item  $\hat\nu([m]) =(0,0)$ \ {\it iff}\ \   $[m] \in$ Sing(Sing$({\cal M}^c_{\text{st}}))\cap {\cal M}^c_3$;
\item   $\hat\nu([m])  \in \Gamma_{\text{Cardano}}$ \ {\it iff} \  $m$  admits a non-trivial automorphism.
\end{itemize}
In this appendix we will extend all of the above results to ${\cal M}^c_{\text{st}}$.
In order to do this we will give normal forms for algebras in ${\cal M}^c_2$ and ${\cal M}^c_1$ based
on the fact that algebras in ${\cal M}^c_2$ have a basis  of  elements  satisfying $m(x,x)= \lambda x$, and
algebras in ${\cal M}^c_1$ have a one-dimensional subspace of elements satisfying $m(x,x)= \lambda x$.

\begin{theorem}\label{theo:tutu}
In Table \ref{tab:1} (resp. Table \ref{tab:2}) we give the normal forms,  continuous moduli, discrete moduli,
the automorphisms, the values of the  coordinates  $\tilde p_2,\tilde p_3$,
and the values of the invariants  $Q_m,$ {\rm Disc}$(D_m)$, {\rm Inv}$_m$
for ${\cal M}^c_2$ (resp. ${\cal M}^c_1$).  ``Discrete'' moduli  are only  needed when the field
$\FF$ is not algebraically closed.
In these tables we have written $e_1{}^2$ for $m(e_1,e_1)$ etc and $\{\epsilon_1,\epsilon_2\}$ for the basis
dual to $\{e_1,e_2\}$.

\begin{table}
\begin{center}
\begin{tabular}{|c|c|c|c|}
\hline
Normal form&Moduli&Automorphisms&\begin{tabular}{c} Coordinates $\tilde p_2, \tilde p_3$; \\ Invariants $Q_m,$ Disc$(D_m).$
\end{tabular}
\\
\hline
\begin{minipage}{4.1cm}
 \beqa
e_1^2&=&e_1 \nn\\
e_2^2&=&e_2\nn\\
e_1 e_2 &=&\frac12 e_1 + \nu e_2\nn\\ \nn
\eeqa\end{minipage}
&
\begin{minipage}{3cm}
\begin{center}
\bigskip
continuous:\\ \smallskip 
 $\nu \in \FF \setminus\{\frac12\}$.\\ \color{white} marcus\\\color{black}
discrete:\\
none.
\end{center}
\end{minipage}
& $\{1\}$&
\begin{minipage}{6.5cm} \vskip -.4truecm 
\hskip -1.truecm \beqa
&&\hskip -.5truecm \tilde p_2 = - 9 (2\nu -1)^2 (\epsilon_1 \wedge \epsilon_2)^2 \nn \\
&&\hskip -.5truecm \tilde p_3 = - 18 (2\nu -1)(\nu +1) (\epsilon_1 \wedge \epsilon_2)^2 \nn \\
&&\hskip -.5truecm \text{Inv}_m= -\frac12 (2\nu-1)^3 (\epsilon_1 \wedge \epsilon_2)^3 \nn \\
&&\hskip -.5truecm Q_m=(1-2\nu)  x^2 y (e_1 \wedge e_2) \nn\\
&&\hskip -.5truecm \text{Disc}(D_m) =(1-2\nu) (\epsilon_1 \wedge \epsilon_2)^2 \nn
\eeqa
\end{minipage}
\\
\hline
\begin{minipage}{2.4cm}
\beqa
e_1^2&=&0 \nn\\
e_2^2&=&e_2\nn\\
e_1 e_2 &=&\frac12 e_1 +  e_2\nn\\ \nn
\eeqa\end{minipage}&
\begin{minipage}{3cm}
\begin{center}
\bigskip
continuous:\\ 
 point.\\ \color{white} marcus\\\color{black}
discrete:\\
none.
\end{center}\
\end{minipage}
&
 $\{1\}$&
\begin{minipage}{4cm} \vskip -.4truecm 
\beqa
&&\hskip -.5truecm \tilde p_2 = - 36 (\epsilon_1 \wedge \epsilon_2)^2 \nn \\
&&\hskip -.5truecm \tilde p_3 = - 36 (\epsilon_1 \wedge \epsilon_2)^2 \nn \\
&&\hskip -.5truecm \text{Inv}_m= -4(\epsilon_1 \wedge \epsilon_2)^3 \nn \\
&& \hskip -.5truecm Q_m=-2  x^2 y (e_1 \wedge e_2) \nn\\
&& \hskip -.5truecm \text{Disc}(D_m) = 0  \nn
\eeqa
\end{minipage}
\\
\hline
\end{tabular}
\end{center}
\caption{$m\in {\cal M}^c_2$}
\label{tab:1}
\end{table}

\begin{table}
\begin{center}
\begin{tabular}{|c|c|c|c|}
\hline
Normal form&Moduli&Automorphisms&Invariants\\
\hline
\begin{minipage}{3cm}
\beqa
e_1^2&=&e_1 \nn\\
e_2^2&=&\lambda e_1\nn\\
e_1 e_2 &=&\frac12 e_2 \nn\\ \nn
\eeqa\end{minipage}
&
\begin{minipage}{3cm}
\begin{center}
\bigskip
continuous:\\ \smallskip 
point.\\ \color{white} marcus\\\color{black}
discrete:\\\medskip
$[\lambda] \in \FF^\ast/\FF^{\ast 2}$.
\end{center}
\end{minipage}
&$\mathbb Z_2$&
\begin{minipage}{7cm} \vskip -.4truecm 
\beqa
&&\hskip -.5truecm \tilde p_2 = 0 \nn \\
&&\hskip -.5truecm \tilde p_3 = - 27\lambda(\epsilon_1 \wedge \epsilon_2)^2 \nn\\
&&\hskip -.5truecm \text{Inv}_m= 0 \nn \\
&&\hskip -.5truecm Q_m=\lambda   y^3 (e_1 \wedge e_2) \nn\\
&& \hskip -.5truecm \text{Disc}(D_m) = \lambda (\epsilon_1 \wedge \epsilon_2)^2 \nn 
\eeqa
\end{minipage}\\
\hline
\end{tabular}
\end{center}
\caption{$m\in {\cal M}^c_1$}
\label{tab:2}
\end{table}

\end{theorem}

\begin{demo}
If $[m] \in {\cal M}^c_2$, there is a basis of elements satisfying $m(x,x) = k x$ and if $[m] \in {\cal M}^c_1$ there is
a unique one-dimensional nontrivial subalgebra. 
If $ [m] \in {\cal M}^c_2$ it is  relatively straightforward to show that   the normal forms
in Table \ref{tab:1} are the only possible ones even if $\FF$ is not algebraically closed.
If $ [m] \in {\cal M}^c_1$   there is a unique idempotent  $e_1$ and
those multiples of $e_1$ which are perfect squares define $[\lambda] \in \FF^\ast/\FF^{\ast 2}$.
It turns that this class completely characterises  algebras in ${\cal M}^c_1$.
\end{demo}

\begin{remark}
If $\FF$ is algebraically closed ${\cal M}^c_1$ consists of a single point as was shown in \cite{am}.
\end{remark}

 Recall that
\beqa
\widetilde{\Gamma}_{\text{Cardano}}= \Bigg\{ [a,b,c] \in \mathbb P^2 \ \ \text{s.t.} \ \  27 c a^2+ 4 b^3 = 0 \Bigg\} \ , \nn
\eeqa
is the `projective' extension of the `affine' curve ${\Gamma}_{\text{Cardano}}\subset \FF^2$ and that
$\Phi: {\cal M}^c_{\text{st}} \to \mathbb P^2$  is the projective extension of $\hat\nu: {\cal M}^c_3 \to \FF^2$.

It follows  immediately from the Tables \ref{tab:1} and \ref{tab:2}
that the properties given above for  generic algebras extend to stable algebras as follows:
\begin{itemize}
\item if $\Phi([m]) \not \in \widetilde{\Gamma}_{\text{Cardano}}$ then $\Phi([m])$ completely determines the equivalence class of $m$
(even if $\FF$ is not algebraically closed);
\item if $\Phi([m])  \in \widetilde{\Gamma}_{\text{Cardano}} \setminus \{[0,0,1]\}$ then $\Phi([m])$
together with the equivalence class of a (at most) quadratic extension of $\FF$  determines the equivalence class of $m$;
\item if  $\Phi(m) =[0,0,1]$  then $\Phi([m])$ together with the equivalence class of the splitting field of the
fundamental cubic $Q_m$ completely determine the equivalence class of $m$; 
\item  $\Phi([m])  \in \widetilde{\Gamma}_{\text{Cardano}} \setminus\{[0,0,1]\}$ \ {\it iff} \  $[m] \in$
$\Big($Sing$({\cal M}^c_{\text{st}}) \setminus$  Sing(Sing$({\cal M}^c_{\text{st}}))\Big)$;
\item  $\Phi([m]) =[0,0,1]$ \ {\it iff}\ \   $[m] \in$ Sing(Sing$({\cal M}^c_{\text{st}}))$;
\item   $\Phi([m])  \in \widetilde{\Gamma}_{\text{Cardano}}$ \ {\it iff} \  $m$  admits a non-trivial automorphism.
\end{itemize}

\subsection{Non-stable commutative algebras}
For completeness we give the  classification of non-stable  algebras which by definition are   those algebras
on which   all scalar-valued invariant polynomials vanish.
\begin{theorem}
In Table \ref{tab:3}  we give the normal forms, moduli,
 automorphisms, 
and the value of the fundamental cubic  $Q_m$ for non-stable commutative algebras.
In these tables we have written $e_1{}^2$ for $m(e_1,e_1)$ etc and $(\epsilon_1,\epsilon_2)$ for the basis
dual to $(e_1,e_2)$.
\begin{table}\
\begin{center}
\begin{tabular}{|c|c|c|c|}
\hline
Normal form&Moduli&Automorphisms& Fundamental cubic\\
\hline
\begin{minipage}{4.1cm}
 \beqa
e_1^2&=&e_1 \nn\\
e_2^2&=&0\nn\\
e_1 e_2 &=&\nu e_2\nn\\ \nn
\eeqa\end{minipage}
&
\begin{minipage}{3cm}
\begin{center}
\bigskip
continuous: \\ \smallskip 
$\nu \in \FF \setminus\{\frac12\}$.\\ \color{white} marcus\\\color{black}
discrete:\\\medskip
none.
\end{center}
\end{minipage}
&$\FF^*$& $Q_m = (1-2\nu) y x^2 (e_1\wedge e_2)^2$\\
\hline
\begin{minipage}{4.1cm}
 \beqa
e_1^2&=&0 \nn\\
e_2^2&=&0\nn\\
e_1 e_2 &=& e_2\nn\\ \nn
\eeqa\end{minipage}
&
\begin{minipage}{3cm}
\begin{center}
\bigskip
continuous: \\ \smallskip 
Point.\\ \color{white} marcus\\\color{black}
discrete:\\\medskip
none.
\end{center}
\end{minipage}
&$\FF^*$& $Q_m = -2 y x^2 (e_1\wedge e_2)^2$\\
\hline
\begin{minipage}{4.1cm}
 \beqa
e_1^2&=&0 \nn\\
e_2^2&=&e_1 + \delta e_2\nn\\
e_1 e_2 &=& \frac12 \delta e_1 \nn\\ \nn
\eeqa\end{minipage}
&
\begin{minipage}{3cm}
\begin{center}
\bigskip
continuous: \\ \smallskip 
$\delta \in \{0,1\}$.\\ \color{white} marcus\\\color{black}
discrete:\\\medskip
none.
\end{center}
\end{minipage}
&$\FF^*$& $Q_m = y^3 (e_1\wedge e_2)^2$\\
\hline
\begin{minipage}{4.1cm}
 \beqa
e_1^2&=&e_1 \nn\\
e_2^2&=&0\nn\\
e_1 e_2 &=& \frac12 e_2\nn\\ \nn
\eeqa\end{minipage}
&
\begin{minipage}{3cm}
\begin{center}
\bigskip
continuous: \\ \smallskip 
Point.\\ \color{white} marcus\\\color{black}
discrete:\\\medskip
none.
\end{center}
\end{minipage}
&$\FF^*$& $Q_m =0$\\
\hline
\end{tabular}
\end{center}
\caption{Non-stable algebras}
\label{tab:3}
\end{table}
\end{theorem}
\begin{demo} The proof of this theorem is very similar to the proof of the theorem above.
\end{demo}
\begin{remark}
It is interesting to observe that non-stable commutative algebras are characterised  amongst all commutative
algebras as being those whose automorphism group is isomorphic to $\FF^*$. Similarly stable commutative
algebras are characterised as those whose automorphism group is finite.
\end{remark}

\begin{remark}
Note that at the end of the article \cite{am} a description of the moduli space of non-stable not
necessarily commutative algebras is given.
\end{remark}

\phantomsection
\newpage
\bibliographystyle{utphys}
\bibliography{TPI}

\end{document}